\documentclass[manuscript, nonacm, natbib=false]{acmart}

\author{Sudhanva Vijaya Kulkarni}
\affiliation{%
  \institution{Arith Inc, USA}
  \city{Berkeley}
  \state{CA}
  \country{USA}
}
\email{sudhanvakulkarni@berkeley.edu}

\author{James W. Demmel}
\affiliation{%
  \institution{University of California, Berkeley}
  \city{Berkeley}
  \state{CA}
  \country{USA}
}
\email{demmel@berkeley.edu}

\author{Laura Grigori}
\affiliation{%
  \institution{PSI Center for Scientific Computing, Theory and Data and Institute of Mathematics, EPF Lausanne, Switzerland}
  \city{Lausanne}
  \country{Switzerland}
}
\email{laura.grigori@epfl.ch}

\acmJournal{TOMS}
\acmVolume{0}
\acmNumber{0}
\acmArticle{0}
\acmYear{2026}
\acmMonth{1}
\acmDOI{10.1145/xxxxxxx.xxxxxxx}
\copyrightyear{2026}

\usepackage{subcaption}
\usepackage{amsmath}
\usepackage{amssymb,amsfonts}
\usepackage{graphicx}
\usepackage{textcomp}
\usepackage{xcolor}
\usepackage{algorithm}
\usepackage{algpseudocode}
\usepackage{tikz}
\usetikzlibrary{calc,shapes.geometric,arrows,positioning}
\usepackage{amsthm}
\usepackage{booktabs}
\usepackage{array}
\usepackage{graphicx}
\usepackage{subcaption}
\usepackage[utf8]{inputenc}
\usepackage[style=acmnumeric]{biblatex}
\allowdisplaybreaks

\newtheorem{theorem}{Theorem}[section]
\newtheorem{lemma}[theorem]{Lemma}
\newtheorem{corollary}{Corollary}[theorem]

\renewcommand{\arraystretch}{1.4}
\newcommand{\eh}{\epsilon_h}
\newcommand{\gfp}{\gamma_{\mathrm{fp}}}

\newif\ifshownotes
\shownotesfalse    

\newcommand{\ignore}[1]{}

\newcommand{\jwdnote}[1]{\ifshownotes{\color{red}{[JWD: #1]}}\else\expandafter\ignorespaces\fi}
\newcommand{\jwdnotenew}[1]{\ifshownotes{\color{red}{[JWD: #1]}}\else\expandafter\ignorespaces\fi}
\newcommand{\jwdnotenewnew}[1]{\ifshownotes{\color{red}{[JWD: #1]}}\else\expandafter\ignorespaces\fi}
\newcommand{\svknote}[1]{\ifshownotes{\color{blue}{[SVK: #1]}}\else\expandafter\ignorespaces\fi}
\newcommand{\svknotenew}[1]{\ifshownotes{\color{blue}{[SVK: #1]}}\else\expandafter\ignorespaces\fi}
\newcommand{\svknotenewnew}[1]{\ifshownotes{\color{blue}{[SVK: #1]}}\else\expandafter\ignorespaces\fi}
\newcommand{\lgnote}[1]{\ifshownotes{\color{magenta}{[LG: #1]}}\else\expandafter\ignorespaces\fi}
\tikzstyle{inputblock} = [draw, thick, minimum width=1.8cm, minimum height=0.6cm, align=center, fill=blue!20]
\tikzstyle{matrixblock} = [draw, thick, minimum width=2.2cm, minimum height=0.6cm, align=center, fill=green!30]
\tikzstyle{outputblock} = [draw, thick, minimum width=2.5cm, minimum height=0.6cm, align=center, fill=yellow!30]
\tikzstyle{block} = [draw, thick, rounded corners=5pt, minimum width=2.0cm, minimum height=1cm, align=center, fill=blue!20]
\tikzstyle{arrow} = [thick, ->, >=stealth]

\begin{document}
\title{Adaptive Precision Cholesky with Error Guarantees}

\begin{abstract}
Recent hardware accelerators increasingly support low-precision numerics such as 8-bit floats and ints that are used to speed up machine learning applications. In this paper, we show that a Cholesky factorization in which both GEMM and TRSM are implemented using the fixed-point Ozaki 1 scheme \cite{Ootomo_Ozaki_Yokota_2024} satisfies the same backward error bounds as the standard floating-point Cholesky factorization \cite{demmel1989cholesky}. We also show that a Cholesky factorization computed using a low-precision MicroScaled GEMM \cite{ocp2023mx} satisfies the same bound as \cite{demmel1989cholesky}. We then use our error analyses to construct an adaptive precision Cholesky factorization that can exploit low-precision formats while respecting a prescribed element-wise error bound threshold, as long as the factorization is successful. We test our factorization in two settings. First, we use it to choose the number of Ozaki 1 slices at each iteration of
the Cholesky factorization, with SYRK/GEMM and TRSM 
emulated in low-precision (INT8) fixed-point arithmetic to attain double-precision accuracy. Second, we use it to build a low-precision Cholesky preconditioner for the conjugate gradient method (PCG), in which SYRK/GEMM take low-precision inputs and accumulate in high precision, and the input format is selected adaptively for each panel. If the factorization breaks down due to loss of positive definiteness, we use diagonal perturbation to ensure we still construct a valid preconditioner, but lose error bound guarantees.
We test the adaptive precision Ozaki 1-based factorization on an NVIDIA GH200, and compare it against a right-looking Cholesky implemented using two kinds of BLAS: first using cuBLAS' Automatic Dynamic Precision (ADP), and then using a fixed slice count of 7 for all Ozaki 1 BLAS operations. We report speedups of up to 29\% 
\jwdnotenewnew{Later sections refer to a 15\% speedup.} \svknote{Fixed}
over the latter and up to 69\% over the former, while still maintaining the same error guarantees.
For PCG, we show that with a diagonal perturbation, FP8 factorizations form effective preconditioners on the 13 tested SuiteSparse matrices, reducing the number of PCG iterations by up to  $12\times$ over CG on the same matrix diagonally scaled by powers of two. For synthetic matrices with arithmetic spectra, and geometric spectra with condition numbers $\kappa \le 10^4$, FP8 is able to bring down the number of CG iterations. \jwdnotenewnew{Incomplete sentence.}
\jwdnote{Anything about speedups with respect to
standard floating point?} \svknote{If this refers to speedup of emulation over floating point, there is no speedup on H200 as we discussed in our meeting}
\svknote{TODO:  I don't like the last sentence since it doesn't clearly convey the tradeoff we are characterizing}
\jwdnote{To be updated when all results are finalized.}
\end{abstract}



\maketitle

\raggedbottom

\section{Introduction and Related Work}
Many recent studies have shown that mixed-precision hardware can significantly accelerate the solution of linear systems (see \cite{Kashi2024MixedPrecSurvey} for a survey). We consider two approaches that have been used to achieve these speedups. The first approach combines a low-precision factorization with iterative refinement to recover a higher accuracy solution of the linear system. For example Haidar et al.
\cite{haidar-mixed-prec} achieve significant speedups by combining a mixed-precision LU factorization that exploits half-precision NVIDIA Tensor Cores with GMRES-based iterative refinement (GMRES-IR). In addition, GMRES-IR was accelerated using multiple precisions in  \cite{five-prec-GMRES}. The strategy in both these works is to use an FP16/TF32/BF16 LU factorization and then use it to precondition the "correction" equation, which is solved with GMRES. Iterative refinement is required to compensate for the fact that FP16 and BF16 both have many fewer mantissa bits than FP64. This implies a larger machine epsilon and hence a less accurate solution. In this work, since we are only concerned with symmetric positive definite systems, we use the mixed-precision factorization as a preconditioner for a preconditioned conjugate gradient (PCG) rather than GMRES. Mixed-precision PCG with low-precision incomplete Cholesky preconditioners has been analyzed in \cite{Bake2025MixedPCG}. Low-precision formats like FP16 also provide a much smaller range of representable numbers. Methods to mitigate the effect of the smaller range of available FP16 numbers have been studied in \cite{squeezing_fp16}. Sahraneshinsamani et al.~\cite{Sahraneshinsamani2024} show how pre-pivoting with low-precision GEMMs can  eliminate the pivoting bottleneck on GPU-like architectures. Other approaches include using custom data structures 
such as those described  in \cite{Amestoy2022Mixed}, where a mixed-precision SVD is used to store different blocks of a Block-Low-Rank matrix in different precisions. This has the effect of storing "less significant" (smaller singular value) 
 parts of the blocks in lower precision. Another example is \cite{mixed-prec-chol}, where off-diagonal blocks beyond some user-specified distance away from the diagonal
 are stored in FP32, and the rest are stored in FP64. A similar idea is applied in \cite{Ren2024Cholesky} where the precision to which the off-diagonal blocks are rounded is decided by the magnitude of  their Frobenius norm (blocks with a smaller Frobenius norm are stored in lower precision). This idea is further extended in \cite{MXCholesky} while also applying MicroScaling. The Cholesky factorizations described in \cite{mixed-prec-chol, Ren2024Cholesky}
have the problem of not providing error-bound guarantees, and not accounting for problems such as factorization failure and underflow, which become common in the world of low-precision floats.
The concurrent work \cite{Cao2026RightPrecision} solves this issue by picking the precision of each tile using the ratio between the norm of the tile and the norm of the input matrix. Our work differs from that work in that it also considers emulated GEMMs and can provide either an element-wise or norm-wise error bound rather than only a norm-wise bound. Our algorithm also uses a different precision selection criterion.
 More recently, FP8/6/4 formats have also become popular due to their utility in machine learning
\cite{micikevicius2022fp8formatsdeeplearning}. Many FP8 (and shorter) formats are available \cite{ocp2023mx}, \cite{noune2022numerical}, \cite{tesla2023dojo}, \cite{10.14529/jsfi170206}, \cite{luo2024ascendhifloat8formatdeep}, including the emerging IEEE 3109 standard for low-precision floats \cite{p3109_public_repo}. 
When using these \jwdnotenewnew{three?} low precisions, the Cholesky factorization can break down. We avoid this by using diagonal perturbation, but lose error guaranties in the process.

The second approach relies on emulation, which splits a higher precision GEMM into a sum of many scaled low-precision fixed-point or floating-point GEMMs. Ozaki schemes \cite{Ootomo_Ozaki_Yokota_2024, ozaki2025ozakischemeiigemmoriented} use lower-precision GEMMs to decompose a floating-point matrix into a sum of matrices whose products can be calculated exactly using low precision and then accumulated using higher precision. This allows FP64 matrix products to be computed using INT8 Tensor Cores, 
albeit with a weaker roundoff error bound \cite{demmel2026gradeaccuracyblas}. The error bound of a GEMM implementation can receive a grade $A, B$ or $C$, where the meaning of these grades is displayed in Table~\ref{tab:gemm-grades}.
In this paper, we show that using Ozaki (or any "grade B" GEMM as defined in \cite{demmel2026gradeaccuracyblas}) gives you the same backward error guarantees as the regular floating-point Cholesky \cite{demmel1989cholesky}.
Hence, when using emulation, we can simply replace DGEMM or DTRSM with their emulated equivalents, and we do not need an extra refinement phase. Luszczek et al.~\cite{emulatedLU} used the Ozaki scheme for the BLAS in an LU decomposition, but did not present an error analysis. Kouya~\cite{Kouya2026OzakiLU} also uses Ozaki 2 in an LU decomposition. The Ozaki scheme has also been extended to FP8 \cite{Mukunoki2025FP8Ozaki} and FP4 \cite{Hayashi2026FP4Ozaki} tensor cores.
NVIDIA recently added a fixed-point INT8 emulation of DGEMM \cite{Schwarz_2026} to cuBLAS. Other implementations include ozIMMU 
\cite{ootomo2024dgemmintegermatrixmultiplication} 
and also the Ozaki 2 implementation \cite{ozaki2025ozakischemeiigemmoriented}. We use the cuBLAS implementation for our experiments.

\begin{table}[h]
\centering
\caption{GEMM error bound grades for $\hat{C} = \text{fl}(AB)$, where $\epsilon$ is machine epsilon and $f_A, f_B, f_C$ are at most low degree polynomial functions of $n$, where for simplicity we assume $A$, $B$ and $C$ are
$n \times n$ matrices.}
\label{tab:gemm-grades}
\begin{tabular}{c l}
\hline
Grade & Error Bound \\
\hline
A & $|\hat{C} - C| \le f_A(n)\,\epsilon\,|A|\,|B|$ \\
B & $|\hat{C}-C|_{ij} \le f_B(n)\,\epsilon\,\|A(i,:)\|\,\|B(:,j)\|$ \\
C & $\|\hat{C}-C\| \le f_C(n)\,\epsilon\,\|A\|\,\|B\|$ \\
\hline
\end{tabular}
\end{table}

\begin{table}[h]
\centering
\caption{TRSM error bound grades for solving  $LX = B$, with computed solution $\hat{X}$, where $\epsilon$ is machine epsilon and $g_A, g_B, g_C$ are at most low degree polynomial functions of $n$,
where for simplicity we assume $L$, $X$ and $B$ are $n \times n$ matrices.} 
\label{tab:trsm-grades}
\begin{tabular}{c l}
\hline
Grade & Error Bound \\
\hline
A & $|L\hat{X}-B| \le g_A(n)\,\epsilon\,|L|\,|\hat{X}|$ \\
B & $|L\hat{X}-B|_{ij} \le g_B(n)\,\epsilon\,\|L(i,:)\|\,\|\hat{X}(:,j)\|$ \\
C & $\|L\hat{X}-B\| \le g_C(n)\,\epsilon\,\|L\|\,\|\hat{X}\|$ \\
\hline
\end{tabular}
\end{table}


\section{Contributions}
Our contributions may be summarized as follows.
\jwdnote{Add sentences like "This yields speedups of up to X
compared to Y." where appropriate,
or some other numerical improvement
metric.}

\begin{itemize}
    \item Our main contribution is an adaptive-precision Cholesky factorization that selects the number of Ozaki 1 slices, or the floating-point precision (FP8, FP16, etc.) for GEMM/SYRK and  the number of Ozaki 1 slices for TRSM per iteration of Cholesky adaptively to reach a prescribed accuracy.
    \jwdnotenewnew{adaptively
    with what goal?}
    \svknotenewnew{Fixed}
     
    \jwdnotenewnew{Ditto.}
    In the appendix, we also give an adaptive algorithm for Ozaki 2 that selects the number of moduli for a prescribed accuracy, based on only a first-order error analysis (we keep only the dominant first-order
   \jwdnotenewnew{$t$ undefined.} terms of the GEMM error bound).

    \item The adaptive algorithm is enabled by our error analysis: we present an error analysis of Cholesky using emulated BLAS for GEMM/SYRK and TRSM, and an error analysis of Cholesky that uses mixed-precision floating-point Tensor Cores for GEMM/SYRK. 
    We show that an emulated Cholesky obeys the same classical backward 
    error bound as the classical floating-point Cholesky.

    \item We evaluate the emulation algorithm on NVIDIA GH200 GPUs using cuBLAS FP64 fixed-point emulation. We compare 3 strategies, using 1) our adaptive algorithm, 2) static slice counts (7 and 8 slices), and 3) cuBLAS' automatic dynamic precision scheme \cite{Schwarz_2026} which decides the number of slices by looking at each GEMM without the context of Cholesky factorization. Our adaptive algorithm is up to 29\% faster than the static 7-slice strategy and up to 69\% faster than cuBLAS ADP while meeting the same backward error bound. 
    Our success criterion for all three implementations is the double precision backward error bound described in \cite{demmel1989cholesky,Higham_main_book}. This means that if Cholesky returns a factor $L$, we check that $|a_{ij} - \sum_k l_{ik}l_{jk}| \le n\epsilon_{64}\sqrt{a_{ii}a_{jj}}$, where the subscript $ij$ denotes the $(i,j)$-th element of the corresponding matrix. 
    We also evaluate the Ozaki 2 variant on synthetic matrices against fixed 12-, 13- and 14-moduli implementations. We test on synthetic matrices and matrices from SuiteSparse 
    \cite{SuiteSparse2011,Kolodziej2019}.
    \item We use the mixed-precision floating-point version of the algorithm to generate preconditioners (for PCG) with varying tolerances, using a diagonal perturbation to handle factorization failures, and evaluate their effectiveness in reducing the number of CG iterations. If the algorithm is given a higher tolerance, the adaptive Cholesky will use more low-precision GEMMs. On the 13 SuiteSparse matrices we tested, and on synthetic matrices with arithmetic and random spectra, an FP8 preconditioner reduces the number of PCG iterations by up to $12\times$ compared to CG on the power-of-two equilibrated matrix. 
    \jwdnotenewnew{``but it fails'' sounds like you are
    referring to the 1 out of 13 SuiteSparse matrices that failed, but the rest of the
    sentence refers to ``matrices'', so more than 1.} \svknotenewnew{Thanks, how about something like this?}
    On synthetic matrices with geometric spectra, the FP8 preconditioner improves the convergence of CG for $\kappa \le 10^4$. However, the convergence of CG is not improved much on synthetic matrices with $\kappa \ge 10^8$.
\end{itemize}

Section~\ref{sec:bg} proves the lemmas and theorems required for the error analysis of a mixed-precision Cholesky. Section~\ref{sec:mp-chol} carries out the error analysis, and Section~\ref{sec:adaptive_alg} states and proves the correctness of our adaptive precision algorithm. Our experiments are presented in Section~\ref{sec:experiments}.

\section{Background}\label{sec:bg}
We first introduce the notation we use throughout the rest of the paper. We use lower case letters for scalars and upper case letters for matrices and submatrices. Additionally, lower case indices are used to denote an element of a matrix and upper case indices are used to denote submatrices. Hence, $a_{ij}$ is the $(i,j)$-th element of the matrix $A$, while $A_{IJ}$ is the $(I,J)$-th submatrix; also, $l_{ik}$ is an entry of $L$ and $L_{IK}$ a subblock. Subscripts can be used for both indices (as $a_{ij}$) and labels (as $\epsilon_h$, $\epsilon_l$ to denote high and low-precision machine epsilon, respectively). We sometimes refer to the precision using the machine epsilon when it makes sense. Iteration/panel/slice number is denoted with a parenthesized superscript (e.g. $\epsilon^{(J)}$ is the precision chosen at iteration $J$). We also use MATLAB notation: $A[i_1:i_2,j_1:j_2]$ denotes the submatrix with rows $i_1$ to $i_2$ and columns $j_1$ to $j_2$. $A[i_1:i_2,:]$ refers to rows $i_1$ to $i_2$ and all columns. We use a similar notation for $A[:, j_1:j_2]$. Finally, $A[i,:]$ is the row $i$ and $A[:,j]$ is the column $j$. 
\subsection{Floating-Point Error Analysis}

Floating-point computations are generally not exact. 
Errors are introduced when a higher precision floating-point number is rounded to a lower precision one. In addition to roundoff errors, there are also underflow and overflow errors. For this paper, we ignore overflow errors and high-precision underflow errors (which could occur during accumulation, but are very small and unlikely). The roundoff model used for high precision binary operations is as follows.
    For a binary operation \( \circ \) performed in (high) precision \( \mathcal{H}\), we have
    \begin{equation} \label{eq:flop_error}
\text{fl}(a \circ b) = (a \circ b)(1 + \delta)
\end{equation}
    with \( |\delta| < \epsilon_h\), where $\epsilon_h$ is the machine epsilon. 
Note that we have ignored the underflow since it is a high-precision operator.
If a number \( a\) is rounded to some lower precision \( \mathcal{L} \) with gradual underflow, then
    \begin{equation}
        \label{eq:round error}
        \text{fl}(a) = a(1 + \delta) + \eta
    \end{equation}
    with \(|\delta| <  \epsilon_l \), \( |\eta| < \mu \epsilon_l\), and $\eta \cdot \delta = 0$, 
    where  \( \mu = 2^{e_{min}}\) and $e_{min}$
    is the smallest exponent in the format \cite{IEEE754-2019}. We always assume RoundToNearestEven for the rounding mode unless specified otherwise. We also assume gradual underflow for the rest of the document.

\subsection{Emulation with Ozaki 1}
\label{sec:ozaki}

The Ozaki scheme emulates a high-precision matrix product using a sequence of
low-precision products~\cite{Ootomo_Ozaki_Yokota_2024,ozaki2025ozakischemeiigemmoriented,Uchino_2025}. In this paper, we consider the Ozaki 1 algorithm from \cite{Ootomo_Ozaki_Yokota_2024} 
applied to FP64 matrices with INT8/UINT8 slices. Consider two FP64 matrices $A,B$, divided into $p \in \mathbb{N}$ slices, where a slice represents a
short set of consecutive mantissa bits
in each matrix entry.
We can write
\begin{equation}
    A = \sum_{q=1}^{p} A^{(q)} + R_A^{(p)}, \qquad
    B = \sum_{q=1}^{p} B^{(q)} + R_B^{(p)},
    \label{eq:ozaki-split}
\end{equation}
where each $A^{(q)}, B^{(q)}$ is a scaled low-precision INT8/UINT8 matrix and
$R_A^{(p)} = A - \sum_{q=1}^pA^{(q)}$, $R_B^{(p)} = B - \sum_{q=1}^pB^{(q)}$ are the residuals. The matrix product using $p$ slices is approximated by
\begin{equation}
    \widetilde{C}^{(p)} \;:=\;
    \sum_{q=1}^p \sum_{r=1}^{p - q + 1} A^{(q)}B^{(r)}
    \label{eq:ozaki-approx}
\end{equation}
which contains $\tfrac{1}{2}p(p+1)$ low-precision GEMMs. Each $A^{(q)} B^{(r)}$ is
computed exactly.
We let $\hat{C}^{(p)}$ denote the
final result of accumulating all the
products $A^{(q)}B^{(r)}$ using FP64. Then, for the uniform INT8 slice case described in \cite{Uchino_2025}, this gives the following error bound
\begin{equation}\label{eqn_ChatvsCtilde}
\begin{aligned}
| \hat{C}^{(p)} - \widetilde{C}^{(p)} |
&\leq \Bigl(\frac{1}{2}p(p+1) - 1\Bigr) \epsilon_{64}
\sum_{q=1}^p \sum_{r=1}^{p-q+1}|A^{(q)} B^{(r)}| \\
&\leq \Bigl(\frac{1}{2}p(p+1) - 1\Bigr) \epsilon_{64}
\sum_{q=1}^p \sum_{r=1}^{p-q+1}|A^{(q)}| |B^{(r)}| \le \Bigl(\frac{1}{2}p(p+1) - 1\Bigr) \epsilon_{64}|A||B|
\end{aligned}
\end{equation}

Current implementations of Ozaki 1 FP64 emulation with integers either use INT8 for all slices \cite{Uchino_2025}, or use INT8 for the first slice (to capture the sign bit) and UINT8 for the remaining slices \cite{Schwarz_2026}. The latter approach can capture more mantissa bits while using the same number of slices.  
An error analysis for the uniform INT8 case was presented in \cite{Uchino_2025}. Here, we partially extend that analysis to the cuBLAS implementation described in \cite{Schwarz_2026}.  
Let $d : \mathbb{N} \to \mathbb{N}$ denote a slice width function, where
$d(q)$ is the number of mantissa bits represented by the $q$-th slice. Define $d(0) = 0$. 
Then, the effective mantissa width of a $q$-slice emulation is
\begin{equation}
    w(q) = \sum_{r=0}^{q} d(r) .
\end{equation}
 Two examples are:
\begin{itemize}
    \item {Uniform slices.} If every slice is INT8, then
        $d(q) = 7$ for $q > 0$, and $d(0) = 0$. So $w(q) = 7q$.
    \item {Signed leading slice.} If the leading slice is stored as signed INT8 and subsequent slices as unsigned UINT8, then $d(0) = w(0) = 0$, $d(1) = w(1)  = 7$ and $d(q) = 8$ for $q \ge 2$, giving
        $w(q) = 7 + 8(q-1) = 8q-1$.
\end{itemize}



\subsubsection{Error analysis.}
 Consider multiplying $A,B$ using a $p$-slice Ozaki 1 scheme. Define the vectors $g$ and $f$ with $g_i = \text{ufp}\!\left(\max_j |a_{ij}|\right)$ and $f_j = \text{ufp}\!\left(\max_i |b_{ij}|\right)$. Here $\text{ufp}$ (unit in the first place) is defined as $\text{ufp}(y) = 2^{\lfloor \log_2 |y| \rfloor} \text{ for } y \ne 0, \text{ and } \text{ufp}(0) := 0$.
The error in an Ozaki 1 GEMM consists of 2 parts, the truncation error and the accumulation error.
First, we bound the truncation error $|AB - \widetilde{C}^{(p)}|$
due to splitting the mantissa bits. This requires bounding $|A^{(q)}|, |R_A^{(q)}|, |B^{(q)}|, |R_B^{(q)}|$, for $1 \le q \le p$. These can be bounded above using bounds similar to those in \cite{Uchino_2025}:
\begin{equation}
\begin{aligned}
    |A^{(q)}| &\le 2^{-w(q-1)+1} g \mathbf{1}^T, &
    |R_A^{(q)}| &\le 2^{-w(q) + 1} g \mathbf{1}^T, \\
    |B^{(q)}| &\le 2^{-w(q-1)+1} \mathbf{1} f^T, &
    |R_B^{(q)}| &\le 2^{-w(q) + 1} \mathbf{1} f^T
\end{aligned}
    \label{eq:slice-bound-A}
\end{equation}
for every $1 \le q \le p$. Also, $|A| \le 2\cdot g\mathbf{1}^T$, and $|B| \le 2\cdot \mathbf{1}f^T$, where $\mathbf{1}$ is the column vector of all 1s. 


From~\eqref{eq:ozaki-split}, we have $AB - \widetilde{C}^{(p)} = (\sum_{q=1}^pA^{(q)} + R_A^{(p)})(\sum_{q=1}^pB^{(q)} + R_B^{(p)}) - \sum_{q=1}^p\sum_{r=1}^{p-q + 1}A^{(q)}B^{(r)}$. To transform the residual into a more convenient form, write $\sum_{q=1}^p\sum_{r=1}^{p-q + 1}A^{(q)}B^{(r)} = \sum_{q=1}^pA^{(q)}(\sum_{r=1}^{p-q+1}B^{(r)})$.
Now, since $B = \sum_{r=1}^pB^{(r)} + R_B^{(p)}$, we can write $\sum_{q=1}^pA^{(q)}(\sum_{r=1}^{p-q+1}B^{(r)}) = \sum_{q=1}^pA^{(q)}(B - R_B^{(p - q + 1)}) =  \left (\sum_{q=1}^pA^{(q)} \right)B - \sum_{q=1}^pA^{(q)}R_B^{(p - q + 1)} = \left ( A - R_A^{(p)} \right )B - \sum_{q=1}^pA^{(q)}R_B^{(p - q + 1)}$.
Then by the triangle inequality,
\begin{align*}
    \left| AB - \widetilde{C}^{(p)} \right|
    &\le \sum_{q=1}^{p} |A^{(q)}|\,|R_B^{(p-q+1)}| \;+\; |R_A^{(p)}|\,|B| \le \sum_{q=1}^{p} 2^{-w(q-1)+1}\, g \mathbf{1}^T \cdot 2^{-w(p-q+1)+1}\, \mathbf{1} f^T \;+\; 2 \cdot 2^{-w(p)+1}\, g \mathbf{1}^T \cdot \mathbf{1} f^T.
\end{align*}
Using $\mathbf{1}^T \mathbf{1} = n$, each term in the first sum equals
$4n\, 2^{-(w(q-1) + w(p-q+1))}\, g f^T$.  Assuming $p > 2$, for $q = 1$, $w(0) + w(p) = 8p - 1$, and for $q > 1$,  $w(q-1) + w(p - q + 1) = 
8(q-1)-1 + 8(p-q+1)-1 = 8p-2$, so
\begin{equation}
\begin{aligned} \label{Ozaki-trunc-err}
   \left| AB - \widetilde{C}^{(p)} \right|
    &\;\le\; 4n\, p\, 2^{-8p + 2}\, g f^T \;+\; 2n\, 2^{-8p + 2}\, g f^T
    \le \; 4\left (p+\frac{1}{2} \right )\, n\, 2^{-8p + 2}\, g f^T.
\end{aligned}
\end{equation}

For the accumulation error, we assume that the cuBLAS Ozaki 1 has the same accumulation error as \eqref{eqn_ChatvsCtilde}. Adding together truncation and accumulation errors (from \eqref{Ozaki-trunc-err} and \eqref{eqn_ChatvsCtilde} respectively) gives us
\begin{equation}
\begin{aligned}
|AB - \hat{C}^{(p)}| &\leq
|AB - \widetilde{C}^{(p)}| +
|\widetilde{C}^{(p)} - \hat{C}^{(p)}| \leq 4\left (p+\frac{1}{2}\right )n2^{-8p + 2}gf^T + \left (\frac{1}{2}(p(p+1)) - 1\right )\epsilon_{64}|A||B|.
\end{aligned}
\end{equation}
For convenience, we define
$\Theta(p, n) = 4\left ( p + \frac{1}{2}\right )n2^{-8p + 2}$ and $\Gamma(p) = \frac{1}{2}p(p+1) - 1$, and write our error bound as
\begin{equation}\label{eqn:GEMMerrbnd}
    \left| AB - \hat{C}^{(p)} \right| \leq \Theta(p,n)gf^T + \Gamma(p)\epsilon_{64}|A||B|.
\end{equation}

\subsection{TRSM with Ozaki 1}
We consider solving the system $LX = B$ where $L \in \mathbb{R}^{n\times n}$ is a lower triangular matrix and $B \in \mathbb{R}^{n\times m}$ is a general rectangular matrix. We want to solve for $X$ which has the same dimensions as $B$. There are many possible implementations of TRSM, here we use the textbook blocked implementation of TRSM (shown in Algorithm~\ref{alg:blocked-trsm}, using block size $n_b$) for simplicity. Off-diagonal updates are computed using emulated GEMM, and the smaller TRSM with the $n_b \times n_b$ diagonal block $L_{II}$
is done in regular FP64 floating point. We obtain the following theorem.
\begin{theorem} \label{thm:trsm}
    Let $L \in \mathbb{R}^{n\times n}$, $B \in \mathbb{R}^{n\times m}$, and $X \in \mathbb{R}^{n\times m}$.
    We solve $LX = B$ using a blocked TRSM (with block size $n_b$) where off-diagonal updates are emulated with Ozaki 1 using $p$ slices. Then, for block row and column indices $I,J$, the computed solution $\hat{X}$ satisfies:
    \begin{equation}
        |L\hat{X} - B|_{IJ} \le n_b\epsilon_{64}|L_{II}||\hat{X}_{IJ}| + \Theta(p,n)[\rho \chi^T]_{IJ} + (\Gamma(p) + I)\epsilon_{64}\Bigl(\sum_{K=1}^{I - 1}|L_{IK}||\hat{X}_{KJ}| + |B_{IJ}|\Bigr)
    \end{equation}
    \jwdnote{Isn't the first term covered
    by the 3rd term? Maybe clearer for
    limit of sum in 3rd term not to be
    the upper bound $n/n_b$?} \svknote{It is covered only if $\frac{n}{n_b} \ge n_b$, but I agree that the bound will be correct if the third term only goes till $n/n_b - 1$}
    where $\rho \in \mathbb{R}^n$
    is a column vector with $\rho_i = \max_{j=1}^n |l_{ij}|$
    , and $\chi \in \mathbb{R}^m$ is a column vector
    with $\chi_j = \max_{i=1}^n |\hat{x}_{ij}|$, where $i,j$ are row and column indices,
    and $[\rho \chi^T]_{IJ}$ refers to
    the $I,J$-th subblock of the $n \times m$
    matrix $\rho \chi^T$.
    Also, $\Theta (p,n) = 4 \left ( p + \frac{1}{2}\right )n2^{-8p +2}$, and $\Gamma (p) = \frac{1}{2}p(p+1) - 1$.
\end{theorem}

\begin{algorithm}[H]
\caption{Blocked Right-Looking TRSM: Solve $LX = B$}
\label{alg:blocked-trsm}
\begin{algorithmic}[1]
    \Require Lower-triangular $L \in \mathbb{R}^{n \times n}$,
             $B \in \mathbb{R}^{n \times m}$,
             block size $n_b$
    \Ensure  $X \in \mathbb{R}^{n \times m}$ 
             such that $LX = B$
    \Statex Partition $L$ and $B$ into $n/n_b \times m/n_b$ blocks of size $n_b \times n_b$
        with $L_{IJ}, B_{IJ}$ the $(I,J)$ 
        submatrices

    \For{$J = 1, \dots, m/n_b$}
    \For{$I = 1, \ldots, n/n_b$}
        \For{$K = 1, \dots, I - 1$}
            \State $B_{IJ} \gets B_{IJ} - L_{IK}\, X_{KJ}$
                \Comment{GEMM update of trailing panels} 
        \EndFor
        \State $X_{IJ} \gets L^{-1}_{II}B_{IJ}$
            \Comment{Unblocked TRSM using FP64}
        \EndFor
    \EndFor
\end{algorithmic}
\end{algorithm}
\begin{proof}
For simplicity, assume $n_b|n$ and $n_b|m$. 
We need to do a floating-point TRSM on the updated block $\hat{B}_{IJ} = B_{IJ} - \sum_{K=1}^{I-1} L_{IK}\hat{X}_{KJ} + E_{IJ}$, where $E_{IJ}$ is the error from performing the $I - 1$ GEMMs with Ozaki 1 using $p$ slices. 
For each $L_{IK}, \hat{X}_{KJ}$, denote their column vectors
of row and column ufp maxnorms respectively as $g^{(I, K)}, f^{(J, K)}$. 
\jwdnote{Confusing notation, 
since $g^{(I)}$ also 
depends on $K$. Maybe use $g^{(I,K)}$? Ditto for $f^{(J)}$. Looks like summands in first sums in (11) don't depend on $K$.}
Note that $g^{(I,K)}$ and $L_{IK}$ are 0 in all elements for blocks with $I < K$. Then,
\jwdnote{Compressed formatting of (11). Change (12), (13) and (14) analogously.} \svknote{Thank you, changed the formatting}
\jwdnote{Need to propagate changing
$g^{(I)}$ to $g^{(I,K)}$, ditto for $f^{(J)}$.}
\begin{equation}\label{eqn:Thm3.1}
\begin{aligned}
    |E_{IJ}| &\le \sum_{K=1}^{I - 1} \Theta (p, n_b)g^{(I,K)}f^{(J,K)T} 
    + (\Gamma (p) +  I)\epsilon_{64}\Bigl(\sum_{K=1}^{I - 1}|L_{IK}||\hat{X}_{KJ}| + |B_{IJ}|\Bigr) \\
    &\le \sum_{K=1}^{n/n_b} \Theta (p, n_b)g^{(I,K)}f^{(J,K)T} 
    + (\Gamma (p)  + I)\epsilon_{64}\Bigl(\sum_{K=1}^{I - 1}|L_{IK}||\hat{X}_{KJ}| + |B_{IJ}|\Bigr),
\end{aligned}
\end{equation}
\jwdnote{Where does factor $n/n_b$ come from?
Not in error analysis in sec. 3.2.} \svknote{It comes from there being at most $n/n_b$ GEMMs to update each diagonal block, and each of those GEMMs is accumulated by accumulating $\Gamma(p)$ INT8 matrices in FP64. Will add some text about this}
\jwdnote{But aren't you just adding the
bound in (9) $n/n_b$ times, so the factor
is implicit in the sum?} \svknote{No I don't think we can do that, since we are adding $\Gamma (p) n/n_b$ terms in FP64. If I am adding the result of 2 Cholesky GEMMs, (results of which are C and D), I am adding C\_1 + C\_2 + ... + C\_k + D\_1 + D\_2 + ... + D\_k, so my total roundoff will be 2*k. Does this make sense?} where we increase the upper limit on the first sum from $I - 1$ to $n/n_b$, which is an upper bound since all extra terms are positive. This is also a more convenient form for the inequality, since it includes the diagonal block. The first term on the right hand side is just the sum of all truncation errors. The second term is a collection of all the accumulation errors. To see why we get an extra $\Gamma (p) + I$ factor for the accumulation error,
\jwdnote{$\frac{n}{n_b}$ factor? Ditto later.}
notice that every $p$-slice Ozaki 1 gives us $\Gamma (p) + 1$ FP64 product matrices (which have been multiplied exactly in INT8) to accumulate. After accumulating each of these, we need to sum up the $I - 1$ $L_{IK}\hat{X}_{KJ}$ products that we need to subtract from $B_{IJ}$. Assuming accumulation is done in a nested manner with the $\Gamma (p) + 1$ terms accumulated first, then the $I - 1$ partial sums accumulated after, we get at most $\Gamma (p) + I$ roundoff terms on each element of the result matrix $B_{IJ} - \sum_{K=1}^{I-1}L_{IK}\hat{X}_{KJ}$, thus giving us the required factor. Next, each entry of $g^{(I,K)}, f^{(J,K)}$ is bounded above by the corresponding entry of 
$\rho, \chi$ respectively (since the former is max over submatrices). 
Thus, we can rewrite (\ref{eqn:Thm3.1}) as 
\begin{equation}
\begin{aligned}
       |E_{IJ}| &\le \Theta (p, n) [\rho \chi^T]_{IJ}  + (\Gamma (p) + I)\epsilon_{64}\Bigl(\sum_{K=1}^{I - 1}|L_{IK}||\hat{X}_{KJ}| + |B_{IJ}|\Bigr).
\end{aligned}
\end{equation}
 When we perform the $n_b \times n_b$ TRSM with $L_{II}$, we know from \cite{Higham_main_book} that for every column $j$, 
$\hat{X}_{IJ}$ satisfies
\begin{align*}
(L_{II} + \Delta L_{II}^{(j)})\hat{X}_{IJ}[:,j] = \hat{B}_{IJ}[:,j] \implies \hat{B}_{IJ}[:,j] - L_{II}\hat{X}_{IJ}[:,j] = \Delta L_{II}^{(j)} \hat{X}_{IJ}[:,j],
\end{align*}
with $|\Delta L_{II}^{(j)}| \le n_b\epsilon_{64} |L_{II}|$ for each column $0 < j \le n_b$. So, the inequality per column is $|\hat{B}_{IJ}[:,j] - L_{II}\hat{X}_{IJ}[:,j]| \le n_b\epsilon_{64}|L_{II}| |\hat{X}_{IJ}[:,j]|$. Since this holds for every column, we can say $|\hat{B}_{IJ} - L_{II}\hat{X}_{IJ}| \le n_b\epsilon_{64}|L_{II}| |\hat{X}_{IJ}|$. Thus
, $|B_{IJ} - \sum_{K=1}^{I}L_{IK}\hat{X}_{KJ}| \le |(B_{IJ} - \sum_{K=1}^{I-1}L_{IK}\hat{X}_{KJ}) - \hat{B}_{IJ}| + |\hat{B}_{IJ} - L_{II}\hat{X}_{IJ}|$ can be written as
\begin{equation} \label{b - lk}
\Bigl|B_{IJ} - \sum_{K=1}^{I}L_{IK}\hat{X}_{KJ}\Bigr| \le n_b\epsilon_{64}|L_{II}||\hat{X}_{IJ}| + \Theta (p, n) [\rho \chi^T]_{IJ} + (\Gamma (p) + I)\epsilon_{64}\Bigl(\sum_{K=1}^{I - 1}|L_{IK}||\hat{X}_{KJ}| + |B_{IJ}|\Bigr).
\end{equation}
\end{proof}

For simplicity, we relax the bound in (\ref{b - lk}) as 
\begin{equation} \label{simplified b - lk}
    \begin{aligned}
        \Bigl|B_{IJ} - \sum_{K=1}^{I}L_{IK}\hat{X}_{KJ}\Bigr| &\le \Theta (p,n)[\rho \chi^T]_{IJ} + \epsilon_{64}\max\left(n_b, \Gamma (p) + I\right)\Bigl(\sum_{K=1}^{I}|L_{IK}||\hat{X}_{KJ}| + |B_{IJ}|\Bigr).
\end{aligned}
\end{equation}
Since  $I \le \frac{n}{n_b}$, we obtain
\begin{align}
        \Bigl|B_{IJ} - \sum_{K=1}^{I}L_{IK}\hat{X}_{KJ}\Bigr| &\le \Theta (p,n)[\rho \chi^T]_{IJ} + \epsilon_{64}\max\left(n_b, \Gamma (p) + \frac{n}{n_b}\right)\Bigl(\sum_{K=1}^{n/n_b}|L_{IK}||\hat{X}_{KJ}| + |B_{IJ}|\Bigr) \notag
        \\
        \implies\quad
        |B - L\hat{X}| &\le \Theta (p,n)\rho \chi^T + \epsilon_{64}\max\left(n_b, \Gamma (p) + \frac{n}{n_b}\right)(|L||\hat{X}| + |B|) \label{simplified b - lk full matrix}.
\end{align}

If we would like to remove $B$ from the right hand side, define $R = B - L\hat{X}$. Then, $B = R + L\hat{X}$, and we can rewrite (\ref{simplified b - lk full matrix}) as
\begin{align}
    |R| &\le \Theta (p,n)\rho \chi^T + \epsilon_{64}\max\left(n_b, \Gamma (p) + \frac{n}{n_b}\right)(2|L||\hat{X}| + |R|) \notag
    \\
    \implies\quad
    |R| &\le \frac{1}{1 - \epsilon_{64}\max\left(n_b, \Gamma (p) + \frac{n}{n_b}\right)} \left ( \Theta (p,n)\rho \chi^T + \epsilon_{64}\max\left(n_b, \Gamma (p) + \frac{n}{n_b}\right)(2|L||\hat{X}|) \right )
    \label{b - lk final}.
\end{align}

\begin{corollary}
Blocked TRSM with grade ``B'' GEMM also gets ``B''. 
\cite{demmel2026gradeaccuracyblas}
\jwdnote{Missing reference} \svknote{Will cite the arXiv since it seems like it will be done before this}
\end{corollary}
\begin{proof}
    We can easily rewrite (\ref{b - lk final}) \jwdnote{Missing reference} \svknote{fixed}
    in terms of row and column norms by noting that $|L||\hat{X}| \le n\rho \chi^T$. Therefore
\begin{equation*}
    \begin{aligned}
        |R| &\le \frac{1}{1 - \epsilon_{64}\max\left(n_b, \Gamma (p) + \frac{n}{n_b}\right)} \left(\Theta (p,n)\,\rho \chi^T \;+\; \epsilon_{64}\max\left(n_b, \Gamma (p) + \frac{n}{n_b}\right)\left(2|L||\hat{X}|\right)\right)
        \\[4pt]
        &\le \frac{1}{1 - \epsilon_{64}\max\left(n_b, \Gamma (p) + \frac{n}{n_b}\right)} \left(\Theta (p,n) + 2n\epsilon_{64}\max\left(n_b, \Gamma (p) + \frac{n}{n_b}\right)\right)\left(\rho \chi^T\right).
    \end{aligned}
\end{equation*}
Although this proof is done only for Ozaki 1, it generalizes to a generic grade B GEMM by changing $\Theta$ and $\Gamma$ appropriately.
\end{proof}

\subsection{MicroScaling}\label{sec:mx}
MicroScaling is a strategy used to deal with the limited range of numbers in low-precision floats \cite{ocp2023mx}.  
Given a tensor $A$ that is stored in high-precision (say FP32), rather than rounding naively to lower-precision (say FP8), we can divide the tensor into blocks and scale each block before rounding. The scale factors can be stored with their own precision $scal$, which will be FP32 or E8M0 for our experiments. Throughout the paper, we assume that there is no overflow or underflow for the scale factors.  We refer to the post-scaled-and-rounded elements of each block as private elements. For matrices, these blocks can be either a row/column, a sub-row/column, or a submatrix. A simple implementation of MicroScaling on a vector $V$ is shown in Algorithm \ref{alg:microscaling}. In this algorithm, the ComputeScale function is a helper that computes and returns $\theta_\ell = \text{absmax}(V^{(\ell )})/\beta$ rounded to precision $scal$, where $\beta$ can be any number between 1 and the largest finite number in the target format (to avoid overflow), and $V^{(\ell)}$ is an input block. We refer to $\beta$ as the "microscaling target" since it tells the algorithm what low-precision value the largest number in the block should map to. When rounding the scale factor, $\theta_\ell$ \jwdnote{Add "$\theta_\ell$" to be clear.}, we use RoundAwayFromZero. This ensures that $v^{(\ell)}_i/\theta_{\ell}$ does not overflow if $\beta$ is within $(1 + \epsilon_{scal})$ of the largest finite number of the format. Here $\epsilon_{scal}$ denotes the machine epsilon of the type the scale factor is stored in.  However, if $\beta$ is a smaller number, any rounding mode should be fine.

\begin{algorithm}
\caption{MicroScaling (MX) Quantization of a vector $V$}
\label{alg:microscaling}
\begin{algorithmic}[1]
\Require Vector $V \in \mathbb{R}^n$ in precision $h$, Block size $b_{\mathrm{mx}} \leq n$, target low-precision $l$, microscaling target $\beta$ 
\Ensure Quantized vector $\hat{V}$ in precision $l$ and scale vector $\theta$ in precision $scal$

\State Partition $V$ into contiguous blocks $V^{(1)}, V^{(2)}, \dots, V^{(\lceil n/b_{\mathrm{mx}} \rceil)}$, where $V^{(\ell)}$ contains $v_i$ for $(\ell-1)b_{\mathrm{mx}} < i \le \min(\ell b_{\mathrm{mx}}, n)$
\For{$\ell = 1, \ldots, \lceil n/b_{\mathrm{mx}} \rceil$}
    \State $\theta_\ell \gets \Call{ComputeScale}{V^{(\ell)}, \beta}$
    \For{each element $v_i \in V^{(\ell)}$}
        \State $u_i \gets v_i / \theta_\ell$\label{line:mx-scale}
        \Comment{In precision $h$}
        \State $\hat{v}_i \gets \Call{RoundToLowPrecision}{u_i, l}$
        \Comment{With RoundToNearestEven}
    \EndFor
\EndFor
\State \Return $\hat{V}$, $\theta$
\end{algorithmic}
\end{algorithm}

\begin{lemma} \label{MX_underflow lemma}
    Consider some vector \( V\) in some higher precision \( h \). If we run Microscaling in Algorithm~\ref{alg:microscaling} with some block size
    \( b_{\mathrm{mx}} \le \text{len}(V)\),\jwdnote{Specify the kind of rounding,
    probably ToNearest. Best to do so in Sec. 3.1. Does tie-breaking matter? You refer to Higham's book, but I note that his definitions of emax
    and emin differ from the 754 standard.} \svknote{Changed 3.1 to use IEEE 754 instead of Higham's book, and I don't think tie-breaking matters here.} and obtain a vector of private elements \( \hat{V}\) (in precision $l$), and  a vector of scale factors $\theta$ (in precision $scal$), then \( \forall i \in \{1, \ldots , \text{len}(V)\}\),
    \begin{equation}
    \displaybreak
         \hat{v}_i \theta_{\lceil i/b_{\mathrm{mx}} \rceil} = v_i(1 + \delta_h)(1 + \delta_l) + \theta_{\lceil i/b_{\mathrm{mx}} \rceil}\eta_l,
    \end{equation}
    where $\delta_l \eta_l = 0$, and \(|\delta_l| < \epsilon_{l} \), \(|\delta_h| < \epsilon_{h} \), (and $\delta_h = 0$ if the division is exact), \( |\eta_l| < \mu \epsilon_{l}\),
    where  \( \mu = 2^{e_{min}}\) and $e_{min}$
    is the smallest exponent in the lower-precision format (alternatively, $\mu$ is the smallest normalized number). \jwdnote{Are you assuming there is
    no underflow in line 5 of Alg 2?} \svknote{Yes. Will specify}
    So, if there is underflow, then $|v_i(1 + \delta_h)| \le \mu \theta_{\lceil i/b_{\mathrm{mx}} \rceil}$ 
    \jwdnote{Why the factor $\epsilon_l$?} \svknote{$\mu$ is the smallest normal number, since I assume gradual underflow,  the smallest number in the format would be $\mu \epsilon_l$, then add the scale factor because we scale.}
    \jwdnote{But ``if there is underflow''
    implies the result is less than the
    underflow threshold $\mu$, not the
    smallest subnormal $\mu \epsilon_l$.} \svknote{Ah, you're right. Sorry, somehow I propagated a pretty big misunderstanding about the underflow model throughout the paper. Changing it. The final error bounds themselves shouldn't change much however}
    and  $|\hat{v}_i\theta_{\lceil i/b_{\mathrm{mx}} \rceil} - v_i| \le \mu \epsilon_l \theta_{\lceil i/b_{\mathrm{mx}} \rceil} + \epsilon_h|v_i|$ , 
    \jwdnote{Are you assuming $\delta_h = 0$?} \svknote{Yes I assumed $\delta_h \approx 0$. Thanks for catching this.}
    else $|\hat{v}_i\theta_{\lceil i/b_{\mathrm{mx}} \rceil} - v_i| \le ((1 + \epsilon_h)(1 + \epsilon_l) - 1) |v_i| \approx \epsilon_{l}|v_i|$.
\end{lemma}

\begin{proof}
    Following steps 5,6 of Algorithm \ref{alg:microscaling}, we can observe that 
    \begin{align*}
        \hat{v}_i = \frac{v_i}{\theta_{\lceil i/b_{\mathrm{mx}} \rceil}}(1 + \delta_h)(1 + \delta_l) + \eta_l,
    \end{align*}
    where $\eta_l \delta_l = 0$. Multiplying by $\theta_{\lceil i/b_{\mathrm{mx}} \rceil}$ on both sides gives us the desired result.
\end{proof}
There are three important vectors in this context. The first is the input vector $V$, which is stored in some high precision such as FP32. The second is the low-precision $\hat{V}$ which, along with the scale factors $\theta$, is able to approximate the original $V$. The element $v_i$ is approximated by $\hat{v}_i\theta_{\lceil i/b_{\mathrm{mx}} \rceil}$. We also define a third vector, $\bar{V}$ with $\Bar{v}_i = \hat{v}_i\theta_{\lceil i/b_{\mathrm{mx}} \rceil}$. Thus, $|\Bar{v}_i - v_i|$ represents the quantization error we get when we microscale $V$. It is important to note that $\Bar{V}$ does not include any rounding error that we get from multiplying the scale factors with the low-precision private elements. We have decided to count this roundoff term in the dot product, since the dequantization of these MX vectors occurs during the dot product. Furthermore, we assume that the scaling of the high-precision elements in step~\ref{line:mx-scale} is done via a division, and so we have only $(1 + \delta_h)$ for our roundoff term. If scaling is implemented via a cheaper reciprocal and multiply technique, the roundoff term would be $(1 + \delta_{h1})(1 + \delta_{h2})$ (where $|\delta_{h1}|,|\delta_{h2}|$ are bounded by $\epsilon_h$). Using a power of 2 scale factor would make the scaling exact and so we would not need a $\delta_h$ term at all.
\jwdnote{But maybe underflow?} \svknote{I have assumed no underflow in high precision. Since the scaling happens in high precision, we only say $\delta_h = 0$}

The same lemma can be extended to matrices with any arbitrary MicroScaling layout.
 For MicroScaling in our experiments, we try 2 kinds of scalings. The first one uses one shared FP32 scale factor for a full row/col of FP8 (E4M3 - 1 sign bit, 4 exponent bits, 3 mantissa bits) numbers. The second is OCP MX FP8, which is a shared E8M0 scale factor for a contiguous block of 32 private E4M3 elements.

 \subsubsection{MX GEMM} 
 We consider the problem of multiplying two MicroScaled matrices represented by $\hat{A}, S_A$ and $\hat{B}, S_B$, where $\hat{A}$ has scales $S_A$ for contiguous sub-rows and $\hat{B}$ has scales $S_B$ for contiguous sub-columns. The algorithm for doing this multiplication is written in Algorithm \ref{alg:microscaling GEMM}. \jwdnote{Since you are not actually multiplying $\hat{A}$ and $\hat{B}$, but
 the matrices they represent along with their
 scale factors, would it be clearer to define
  ``microscaled matrix'' as a pair
  $(S_A,\hat{A})$ implicitly representing
  their ``product''? And $(\hat{B},S_B)$ for
  row scaling? Maybe too much notation, but you get the idea.}\svknote{That makes sense to do}\jwdnote{Or maybe just ``multiply the microscaled matrices represented
  by $\hat{A}$ and $S_A$, and by
  $\hat{B}$ and $S_B$'', up to you.}\svknote{How about something like the above?}\jwdnote{Looks good.}
The matrices $S_A, S_B$ are of dimension $m \times \lceil \frac{k}{b_\mathrm{mx}}\rceil$, $\lceil \frac{k}{b_\mathrm{mx}}\rceil\times n$, respectively. Both matrices $\hat{A},\hat{B}$ are assumed to have the same block size $b_{\mathrm{mx}}$ so that each partial dot product of length $b_{\mathrm{mx}}$ uses exactly one scale factor from each $S_A$ and $S_B$. Algorithm \ref{alg:microscaling GEMM} describes how we can compute the higher-precision product $C$. There are other possible precision configurations ($C$ can be low-precision and microscaled, $A,B$ can have different block sizes, etc), but we limit our analysis to this case since it is the only relevant case for our mixed-precision Cholesky factorization (described later in Algorithm~\ref{A_chol}). 

 \begin{algorithm}
\caption{GEMM of two microscaled matrices $\hat{A},\hat{B}$ with scale factors $S_A,S_B$, and microscaling block size $b_{\mathrm{mx}}$}
\label{alg:microscaling GEMM}
\begin{algorithmic}[1]
\Require Matrices $\hat{A} \in \mathbb{R}^{m \times k}, \hat{B} \in \mathbb{R}^{k \times n}$ in precision $l$, Block size $b_{\mathrm{mx}} \leq k$, scale factors $S_A \in \mathbb{R}^{m \times \lceil \frac{k}{b_{\mathrm{mx}}} \rceil}, S_B \in \mathbb{R}^{\lceil \frac{k}{b_{\mathrm{mx}}} \rceil \times n}$, (high) computation precision $h$
\Ensure Output/Accumulation matrix $C \in \mathbb{R}^{m \times n}$ in precision $h$
\For{$i = 1, \ldots, m$}
    \For{$j = 1, \ldots, n$}
    \State $c_{ij} \gets 0$
        \For{$\ell = 1, \ldots, \lceil k/b_{\mathrm{mx}} \rceil$}
            \State $p_\ell \gets 0$
            \Comment{Block-local accumulator, precision $h$}
            \For{$k' = (\ell-1)b_{\mathrm{mx}} + 1, \ldots, \min(\ell \, b_{\mathrm{mx}}, k)$}
                \State\label{line:mx-fma} $p_\ell \gets p_\ell + \hat{a}_{i,k'} \cdot \hat{b}_{k',j}$
                \Comment{Inputs in precision $l$, accumulation in precision $h$}
            \EndFor
            \State\label{line:mx-acc} $c_{ij} \gets c_{ij} + (s_A)_{i,\ell} \, (s_B)_{\ell,j} \, p_\ell$
            \Comment{In precision $h$}
        \EndFor
    \EndFor
\EndFor
\State \Return $C$
\end{algorithmic}
\end{algorithm}
 
 To write the error bounds in matrix form, we define a broadcast matrix that is an expanded version of the matrix of scale factors. The new larger matrix has the same shape as the low-precision matrix of private elements, thus the scaling can be applied entrywise. Let $N = \lceil k/b_{\mathrm{mx}} \rceil$ be the number of blocks along the contraction dimension and let $\mathcal{K}_\ell = \{k' : (\ell-1)b_{\mathrm{mx}} < k' \le \min(\ell b_{\mathrm{mx}}, k)\}$ be the contraction indices in block $\ell$. Then, for each $k' \in \mathcal{K}_\ell$, scale factors $S_A \in \mathbb{R}^{m \times N}$ and $S_B \in \mathbb{R}^{N \times n}$, and matrices $A \in \mathbb{R}^{m \times k}$ and $B \in \mathbb{R}^{k \times n}$, 
 \jwdnote{Why switch notation from $A$ and
 $B$ to $X$ and $Y$? Already a lot of
 notation.} \svknote{True, changed.}
 define the broadcasts $\Bar{S}_A \in \mathbb{R}^{m \times k}$ and $\Bar{S}_B \in \mathbb{R}^{k \times n}$ by
 \begin{equation} \label{eq:blockscale}
     (\Bar{s}_A)_{ik'} = (s_A)_{i,\lceil k'/b_{\mathrm{mx}} \rceil}, \qquad (\Bar{s}_B)_{k'j} = (s_B)_{\lceil k'/b_{\mathrm{mx}} \rceil,j}.
 \end{equation}
Thus, each replicates its scale factors across the $b_{\mathrm{mx}}$ entries of its block, matching the shapes of $A$ and $B$, respectively. The application of scale factors is then just a Hadamard product: $\Bar{S}_A \odot A$ scales $A$ by rows and $B \odot \Bar{S}_B$ scales $B$ by columns. Since Hadamard products are commutative, the ordering of $\Bar{S}_A$ and $\Bar{S}_B$ has no effect here; we nevertheless place $\Bar{S}_A$ on the left and $\Bar{S}_B$ on the right for notational consistency with scaling by diagonal matrices. Now, let $D_A \in \mathbb{R}^{m \times m}$ and $D_B \in \mathbb{R}^{n \times n}$ be the diagonal matrices that hold the largest scale factor of each row of $A$ and each column of $B$,
 \begin{equation} \label{eq:Dscales}
     (d_A)_{ii} = \max_{1 \le \ell \le N} (s_A)_{i\ell}, \qquad (d_B)_{jj} = \max_{1 \le \ell \le N} (s_B)_{\ell j} .
 \end{equation}
 Since $(s_A)_{i,\lceil k'/b_{\mathrm{mx}} \rceil} \le (d_A)_{ii}$ and $(s_B)_{\lceil k'/b_{\mathrm{mx}} \rceil,j} \le (d_B)_{jj}$ for every $k'$, this gives
 \begin{equation} \label{eq:Dbound}
     |\bar{S}_A \odot A| \le D_A|A|, \qquad |B \odot \bar{S}_B| \le |B|D_B ,
 \end{equation}
 with $D_B$ on the right since $B$ is scaled by columns. We can now state the matrix version of Lemma \ref{MX_underflow lemma}.
 \jwdnote{Do you mean Lemma 3.2? Good to use labels, so names change consistently.}

 \begin{lemma} \label{lem:MXmatrix}
     Let $A$ be microscaled along its rows and $B$ along its columns, giving $\hat{A}, S_A$ and $\hat{B}, S_B$, and write $\bar{A} = \bar{S}_A \odot \hat{A}$ and $\bar{B} = \hat{B} \odot \bar{S}_B$ for the exactly dequantized $\hat{A}, S_A$ and  $\hat{B}, S_B$ (these are defined analogously to $\bar{V}$, and don't include the rounding error we get from multiplying the scale factor with the private element). \jwdnote{Confusing: what is the difference between $A$ and $\bar{A}$, and between $B$ and $\bar{B}$, given the
     word ``exactly''?} \svknote{How does the new statement sound?}
     Define $A'$ as $a'_{ik'} = a_{ik'}$ 
     \jwdnote{Is $a_{ik}$ an entry of $A$ or
     $\bar{A}$?} \svknote{A. Entries of $\Bar{A}$ are denoted $\bar{a}_{ij}$}
     when $\hat{a}_{ik'}$ doesn't underflow and 0 otherwise, and ditto for $B'$. Let $M_A, M_B$ have complementary sparsity patterns to $A', B'$ with every nonzero entry equal to $\mu$ (smallest positive normal number
     \jwdnote{Smallest positive normal number?}
     in precision $l)$ 
     \jwdnote{$\epsilon_l$ seems wrong.}, 
     and set $A'' = |A'| + D_AM_A$, $B'' = |B'| + M_BD_B$. Since all entries are positive, $|A''| = A''$ and $|B''| = B''$. Then, with $\phi = (1+\epsilon_h)(1+\epsilon_l) - 1 = \epsilon_l + \epsilon_h + \epsilon_h\epsilon_l$,
     \begin{equation} \label{eq:MXmatrix}
         |\bar{A} - A| \le \phi|A''|, \qquad |\bar{B} - B| \le \phi|B''| .
     \end{equation}
 \end{lemma}
 \begin{proof}
     Given row $i$, entry $(i,k')$ lies in block $\lceil k'/b_{\mathrm{mx}} \rceil$ and carries the scale factor $(\bar{s}_A)_{ik'}$, the broadcasted entry from \eqref{eq:blockscale}, so Lemma \ref{MX_underflow lemma} gives $\bar{a}_{ik'} = a_{ik'}(1+\delta_h)(1+\delta_l) + (\bar{s}_A)_{ik'}\eta_l$ with $\delta_l\eta_l = 0$. Now, if $\hat{a}_{ik'}$ did not underflow, then $\eta_l = 0$ and $|\bar{a}_{ik'} - a_{ik'}| \le ((1+\epsilon_h)(1+\epsilon_l) - 1)|a_{ik'}| = \phi|a_{ik'}|$. If it did underflow, then $\delta_l = 0$, and $|a_{ik'}| < (1 + \epsilon_h)\mu (\bar{s}_A)_{ik'}$. Using this upper bound on $|a_{ik'}|$, we can write $\bar{a}_{ik'} - a_{ik'} = \delta_h a_{ik'} + (\bar{s}_A)_{ik'}\eta_l \implies |\bar{a}_{ik'} - a_{ik'}| \le \epsilon_h (1 + \epsilon_h)\mu (\bar{s}_A)_{ik'} + \epsilon_l \mu (\bar{s}_A)_{ik'} \le \epsilon_h \epsilon_l \mu (\bar{s}_A)_{ik'} +  \epsilon_h \mu (\bar{s}_A)_{ik'}+ \epsilon_l \mu (\bar{s}_A)_{ik'} = \phi \mu (\bar{s}_A)_{ik'}$, \jwdnote{Add "$=\phi \mu (\bar{s}_A)_{ik'}$" for clarity?} \svknote{Makes sense. Added it in.}
     where we used $\epsilon_h^2 \le \epsilon_h \epsilon_l$ for the last bound. The two cases have disjoint supports, so by \eqref{eq:Dbound},
     \begin{equation*}
         |\bar{A} - A| \le \phi|A'| + \phi(\bar{S}_A \odot M_A) \le \phi(|A'| + D_AM_A) = \phi|A''| .
     \end{equation*}
     \jwdnote{Shouldn't $\epsilon_l D_A M_A$ be $D_A M_A$, since $M_A$ already includes $\epsilon_l$?} \svknote{Thanks, changed the definition to be just $\mu$}
     The argument for $B$ is identical, with $D_B$ multiplying from the right.
 \end{proof}

Now, we can state the error bound that we obtain when multiplying two microscaled matrices using Algorithm \ref{alg:microscaling GEMM}. 
 \jwdnote{Of what?}

 \begin{lemma} \label{lem:MXGEMM-DS}
     With the notation of Lemma \ref{lem:MXmatrix}, the output $\hat{C}$ of Algorithm \ref{alg:microscaling GEMM} satisfies
     \begin{equation} \label{eq:MXGEMM}
     \begin{aligned}
         |\hat{C} - AB| &\le \phi(|A''||B| + |A||B''| + \phi|A''||B''|)  + (k + 2)\epsilon_h(|A| + \phi|A''|)(|B| + \phi|B''|) .
     \end{aligned}
     \end{equation}
 \end{lemma}
 \begin{proof}
     Algorithm \ref{alg:microscaling GEMM} applies scale factors from $S_A, S_B$ to the length $b_{\mathrm{mx}}$ dot products of blocks of $\hat{A}, \hat{B}$. These dot products have input precision $l$ and accumulation and output precision $h$. The exact expression for the $(i,j)$ entry of $C$ is 
     \jwdnote{Missing equal sign in middle?} \svknote{Yes, added it in. Thank you}
     \begin{equation} \label{eq:DSexact}
         c_{ij} = \sum_{\ell = 1}^{N}(s_A)_{i\ell}(s_B)_{\ell j}\sum_{k' \in \mathcal{K}_\ell}\hat{a}_{ik'}\hat{b}_{k'j} =\sum_{\ell = 1}^{N}\sum_{k' \in \mathcal{K}_\ell}(s_A)_{i\ell}\hat{a}_{ik'} (s_B)_{\ell j}\hat{b}_{k'j}.
     \end{equation}
     We can turn the nested sum into just one sum by applying each of the scale factors to their corresponding private elements,
     \begin{equation*}
         c_{ij} = \sum_{\ell = 1}^{N}\sum_{k' \in \mathcal{K}_\ell}\bar{a}_{ik'}\bar{b}_{k'j} = \sum_{k' = 1}^{k}\bar{a}_{ik'}\bar{b}_{k'j}.
     \end{equation*}
     Counting the rounding error terms carried by the summand $\bar{a}_{ik'}\bar{b}_{k'j}$, we get two from applying scale factors, at most $b_{\mathrm{mx}} - 1$ from the FMAs forming the block dot product in line~\ref{line:mx-fma}, and at most $N - 1$ from accumulating the $N$ block results in line~\ref{line:mx-acc}. Since $b_{\mathrm{mx}} + N - 1 \le k$ for every $1 \le b_{\mathrm{mx}} \le k$, this is at most $k + 2$ in total, so
     \begin{equation} \label{eq:MXentry}
         \hat{c}_{ij} = \sum_{k'=1}^{k}\bar{a}_{ik'}\bar{b}_{k'j}(1 + \zeta_{k'}), \qquad |\zeta_{k'}| \le (k+2)\epsilon_h .
     \end{equation}
     Hence, $|\hat{c}_{ij} - (\bar{A}\bar{B})_{ij}| \le (k+2)\epsilon_h\sum_{k'=1}^{k}|\bar{a}_{ik'}||\bar{b}_{k'j}|$, which in matrix form reads
     \begin{equation} \label{eq:MXgemm-acc}
         |\hat{C} - \bar{A}\bar{B}| \le (k+2)\epsilon_h|\bar{A}||\bar{B}|.
     \end{equation}
\jwdnote{This assume all operations done in
high precision. But above you make the point
that the preferred approach computes the
short dot product in low precision. Where do
the errors on line 6 of Alg 3 appear?} \svknote{Sorry, inputs are low precision in line 6, accumulation is in high precision. All the low precision error is captured by the truncation terms. Changed the comment in the algorithm to indicate accumulation is in high precision}
Finally, write $E_A = \bar{A} - A$ and $E_B = \bar{B} - B$, so $|E_A| \le \phi|A''|$ and $|E_B| \le \phi|B''|$ by \eqref{eq:MXmatrix}, and hence $|\bar{A}| \le |A| + \phi|A''|$ and $|\bar{B}| \le |B| + \phi|B''|$.
     \jwdnote{$E_A$ and $E_B$ are undefined.} \svknote{Fixed}
     Since $\bar{A}\bar{B} - AB = E_AB + AE_B + E_AE_B$, we get $|\bar{A}\bar{B} - AB| \le \phi(|A''||B| + |A||B''| + \phi|A''||B''|)$. Adding this to \eqref{eq:MXgemm-acc} gives \eqref{eq:MXGEMM}. If we have power-of-two scales, the $k+2$ factor can also be reduced to $k$ since the scaling is exact.
 \end{proof}

We also consider an alternative way to multiply microscaled matrices. This alternate method is depicted in Algorithm~\ref{alg:simulated MX GEMM}, which is clearly more expensive than Algorithm~\ref{alg:microscaling GEMM} since it does not use mixed-precision FMAs with low-precision inputs and high-precision accumulation. Instead, everything is done in high precision. Algorithm~\ref{alg:simulated MX GEMM} is useful in cases where low-precision hardware is not available, since it simulates Algorithm~\ref{alg:microscaling GEMM} while obeying the same error bound as Lemma~\ref{lem:MXGEMM-DS}. 

 \begin{algorithm}
\caption{Simulated GEMM of two microscaled matrices $\hat{A},\hat{B}$ with scale factors $S_A,S_B$, and microscaling block size $b_{\mathrm{mx}}$}
\label{alg:simulated MX GEMM}
\begin{algorithmic}[1]
\Require Matrices $\hat{A} \in \mathbb{R}^{m \times k}, \hat{B} \in \mathbb{R}^{k \times n}$ in precision $l$, Block size $b_{\mathrm{mx}} \leq k$, scale factors $S_A \in \mathbb{R}^{m \times \lceil \frac{k}{b_{\mathrm{mx}}} \rceil}, S_B \in \mathbb{R}^{\lceil \frac{k}{b_{\mathrm{mx}}} \rceil \times n}$, (high) computation precision $h$
\Ensure Output matrix $C \in \mathbb{R}^{m \times n}$ in precision $h$
\For{$i = 1, \ldots, m$}
    \For{$k' = 1, \ldots, k$}
        \State $\bar{a}_{i,k'} \gets (s_A)_{i,\lceil k'/b_{\mathrm{mx}} \rceil} \cdot \hat{a}_{i,k'}$
        \Comment{Dequantize $\hat{A}$; one rounding per entry, precision $h$}
    \EndFor
\EndFor
\For{$k' = 1, \ldots, k$}
    \For{$j = 1, \ldots, n$}
        \State $\bar{b}_{k',j} \gets (s_B)_{\lceil k'/b_{\mathrm{mx}} \rceil, j} \cdot \hat{b}_{k',j}$
        \Comment{Dequantize $\hat{B}$; one rounding per entry, precision $h$}
    \EndFor
\EndFor
\For{$i = 1, \ldots, m$}\label{line:sd-gemm-start}
    \For{$j = 1, \ldots, n$}
        \State $c_{ij} \gets 0$
        \For{$k' = 1, \ldots, k$}
            \State $c_{ij} \gets c_{ij} + \bar{a}_{i,k'} \cdot \bar{b}_{k',j}$
            \Comment{One length-$k$ dot product, entirely in precision $h$}
        \EndFor
    \EndFor
\EndFor\label{line:sd-gemm-end}
\State \Return $C$
\end{algorithmic}
\end{algorithm}

 \jwdnote{Ambiguous. Perhaps use the word
 "order" above to clarify which is "first order" and which is "second"?}
 \jwdnote{second what?}
 \jwdnote{I could only tell the difference between "tilde" and "bar" after viewing at 150\%. Also, meaning of fl() is ambiguous.} \svknote{Oh no. I think I can probably write the proof to not use $\tilde{A}, \tilde{B}$. It should probably also remove the issue with fl()}

 \begin{lemma} \label{lem:MXGEMM-SD}
     With the notation of Lemma \ref{lem:MXmatrix}, the output $\hat{C}$ of Algorithm \ref{alg:simulated MX GEMM} satisfies \eqref{eq:MXGEMM} as well.
 \end{lemma}
 \begin{proof}
     The two dequantization loops of Algorithm \ref{alg:simulated MX GEMM} put one rounding error on each entry of $\bar{A}$ and one on each entry of $\bar{B}$. We then do a regular GEMM from lines~\ref{line:sd-gemm-start} to~\ref{line:sd-gemm-end}, which adds $k$ roundoff terms \cite{Higham_main_book}. This gives $k + 2$ total roundoff terms, each of which is in high precision and is bounded by $\epsilon_h$. Thus, Equation \eqref{eq:MXentry} holds, and the rest of the proof is similar to the proof of Lemma \ref{lem:MXGEMM-DS}.
 \end{proof}

 \begin{theorem} \label{thm:MXGEMM}
     With the notation of Lemma \ref{lem:MXmatrix}, Algorithms \ref{alg:microscaling GEMM} and \ref{alg:simulated MX GEMM} satisfy the same bound \eqref{eq:MXGEMM}.
 \end{theorem}
 \begin{proof}
     Direct application of Lemmas \ref{lem:MXGEMM-DS} and \ref{lem:MXGEMM-SD}.
 \end{proof}

 For our experiments, $\epsilon_l$ will be $0.125$, and $\epsilon_h$ is $2^{-23}$ , so $\epsilon_h \ll \epsilon_l$, and we can simplify $\phi = (1+\epsilon_h)(1+\epsilon_l) - 1 \approx \epsilon_l$, and relax \eqref{eq:MXGEMM} to
 \begin{equation} \label{eq:MXGEMMsimple}
 \begin{aligned}
     |\hat{C} - AB| &\lesssim \epsilon_l(|A''||B| + |A||B''| + \epsilon_l|A''||B''|) + (k+2)\epsilon_h|A||B| .
 \end{aligned}
 \end{equation}
 Note that if we don't do any microscaling, then $N = 1$ with all scale factors equal to 1, and each $D_Z$ matrix is the identity. \jwdnote{each $D_Z$ matrix?} \svknote{fixed} If we use global scaling factors $\theta_A, \theta_B$, then $D_A = \theta_AI$ and $D_B = \theta_BI$. Additionally, if either of the matrices is scaled such that we fully avoid underflow when rounding to low precision, the corresponding $M_Z$ will be 0 and $|Z''| = |Z'| = |Z|$. \jwdnote{what is ``perfectly''?} \jwdnote{each $M_Z$}
 \jwdnote{Last sentence ambiguous.} \svknote{Removed the sentence. I will talk about the accumulation error from doing the + C only when needed (for example in Cholesky).}


\section{Mixed-Precision Cholesky}\label{sec:mp-chol}
To motivate the adaptive algorithm, we first consider the ``regular'' mixed-precision algorithm described in Algorithm \ref{A_chol}.\jwdnote{Now 3?} The algorithm can either use an emulated GEMM/SYRK, or a mixed-precision floating-point GEMM/SYRK (where inputs are assumed to be low-precision and accumulation and output are assumed to be high-precision) for the Schur update in line~\ref{line:chol-gemm}\jwdnote{in line 8.}. We assume that the algorithm receives a list of precisions/slices to use at each iteration, and the slice count/precision in iteration $I$ is given by $\pi_g(I)$ for GEMM and $\pi_t(I)$ for TRSM. The TRSM is implemented such that the off-diagonal GEMMs in the TRSM are done with emulation (as in Algorithm~\ref{alg:blocked-trsm}).
\newcommand{\Output}{\textbf{Output}}
\begin{algorithm}[t]
    \caption{Cholesky with statically scheduled emulated GEMM and TRSM}
    \label{A_chol}
    \begin{algorithmic}[1]
        \Statex \textbf{Input : } PSD Matrix \( A \) in precision \( h \), block size \( n_b\), and for \( 1 \le I \le n/n_b \):
        \Statex \qquad \( \pi_g(I) \): GEMM slice count (emulation) or low-precision format (floating point)
        \Statex \qquad \( \pi_t(I) \): TRSM slice count (emulation only)
        \Statex \textbf{Output : } Lower triangular matrix \( \hat{L} \) such that \( A \approx \hat{L}\hat{L}^T \)

        \For{$I = 1 \rightarrow n/n_b$}
            \State $\hat{L}_{II} \gets \textsc{Potrf}(A_{II})$ \Comment{$A_{II} = \hat{L}_{II}\hat{L}_{II}^T$, in precision $h$}
            \For{$J = I+1 \rightarrow n/n_b$}
                \State\label{line:chol-trsm} $\hat{L}_{JI} \gets \textsc{Trsm}(A_{JI}, \hat{L}_{II})$ \Comment{$\hat{L}_{JI}\hat{L}_{II}^T = A_{JI}$; $\pi_t(I)$ slices, or precision $h$}
            \EndFor
            \For{$J = I+1 \rightarrow n/n_b$}
                \For{$K = I+1 \rightarrow J$}
                    \State\label{line:chol-gemm} $A_{JK} \gets A_{JK} - \textsc{Gemm}(\hat{L}_{JI}, \hat{L}_{KI}^T)$ \Comment{format or slices $\pi_g(I)$}
                \EndFor
            \EndFor
        \EndFor
    \end{algorithmic}
\end{algorithm}

\subsection{Emulation Case}
\svknote{Jim, I have made some changes to section 4.2 according to your recommendations from an earlier email. This section seems like it is much easier to read now :). Please let me know what else I can edit}
\begin{theorem} \label{thm:chol-emu}
    If Algorithm \ref{A_chol} completes successfully with the $\pi_g(I)$ slices used for the Ozaki 1 GEMMs and the $\pi_t(I)$ slices used for the Ozaki 1 TRSM in iteration $I$, assuming that there is no overflow/underflow, provided $n_b\epsilon_{64} < 1$ and $\gamma_{\mathrm{emu}}(\psi, \Psi, n, n_b) < 1$, we obtain $A + E = \hat{L}\hat{L}^T$ with
    \begin{equation} \label{eq:chol-emu}
        |e_{st}| \le \frac{\gamma_{\mathrm{emu}}(\psi, \Psi, n, n_b)}{1 - \gamma_{\mathrm{emu}}(\psi, \Psi, n, n_b)}\sqrt{a_{ss}a_{tt}},
    \end{equation}
    where $\psi = \min_I \min (\pi_g (I), \pi_t (I))$ and $\Psi = \max_I \max (\pi_g (I), \pi_t (I))$ for $1 \le I \le n/n_b$, $\gamma_{\mathrm{emu}}(\psi, \Psi, n, n_b) = \Theta (\psi, n_b) + 4M\epsilon_{64}$ with $M = M(\Psi, n, n_b) = \frac{\Gamma (\Psi) + \max\left(n_b, \frac{n}{n_b}\right)}{1 - n_b\epsilon_{64}}$,  $\Theta (\psi,n_b) = 4 \left ( \psi + \frac{1}{2}\right )n_b2^{-8\psi +2}$, and $\Gamma (\Psi) = \frac{1}{2}\Psi(\Psi+1) - 1$.
\end{theorem}
\begin{proof}
    For simplicity, assume $n_b|n$.
    First consider diagonal blocks of $A$. The blocks $A_{II}$ are updated with Ozaki 1 GEMMs/SYRKs, and then we perform an FP64 Cholesky on the updated submatrix $\hat{A}^u_{II}$. From \cite{demmel1989cholesky}, we have
    \begin{equation} \label{diag cholesky}
        |\hat{A}^u_{II} - \hat{L}_{II}\hat{L}_{II}^T| \le n_b\epsilon_{64}\Omega_{II},
    \end{equation}
    where $\omega_{st} = \sqrt{|\hat{a}^u_{ss}||\hat{a}^u_{tt}|}$.
    Next, by triangle inequality,
    \begin{equation*}
    \begin{aligned}
    \Bigl|A_{II} - \sum_{K=1}^I \hat{L}_{IK}\hat{L}_{IK}^T\Bigr| &= \Bigl|\Bigl(A_{II} - \sum_{K=1}^{I-1} \hat{L}_{IK}\hat{L}_{IK}^T -\hat{A}^u_{II}\Bigr) + \bigl(\hat{A}^u_{II} -  \hat{L}_{II}\hat{L}_{II}^T\bigr)\Bigr| \le \Bigl|A_{II} - \sum_{K=1}^{I-1} \hat{L}_{IK}\hat{L}_{IK}^T -\hat{A}^u_{II}\Bigr|  + \bigl|\hat{A}^u_{II} -  \hat{L}_{II}\hat{L}_{II}^T\bigr|.
    \end{aligned}
    \end{equation*}
    \jwdnote{Already appears above, redundant.} \svknote{True. Removed.}
     The first term on the right hand side can be bounded by the Ozaki 1 error bound \eqref{eqn:GEMMerrbnd}, while accounting for $A_{II}$ in the accumulation roundoff (similar to the proof of Theorem~\ref{thm:trsm}).  Also note that for $p \in \mathbb{N}$, $\Theta(p,n)$ decreases as $p$ increases and $\Gamma(p)$ increases with $p$. Thus $\Theta(\psi, n_b) \ge \Theta (\pi_g(I), n_b)$ for every $1 \le I \le n/n_b$, and $\Gamma(\Psi) \ge \Gamma (\pi_g(I))$ for $1 \le I \le n/n_b$. Ditto for $\pi_t$. Using \eqref{eqn:GEMMerrbnd} and the same counting argument as in the proof of Theorem~\ref{thm:trsm}, we can write 
    \begin{equation}\label{Schur err emulate}
    \begin{aligned}
        \Bigl|A_{II} - \sum_{K=1}^{I-1}\hat{L}_{IK}\hat{L}_{IK}^T - \hat{A}^u_{II}\Bigr| & \le \sum_{K=1}^{I-1} \Theta (\pi_g(K),n_b)(\rho^{(I,K)})(\rho^{(I,K)})^T + \bigl(\max_{K \le I}\Gamma(\pi_g(K)) + I\bigr)\epsilon_{64}\Bigl(\sum_{K=1}^{I-1}|\hat{L}_{IK}||\hat{L}_{IK}^T| + |A_{II}|\Bigr) \\
        &\le \Theta (\psi, n_b)\sum_{K=1}^{I-1}(\rho^{(I,K)})(\rho^{(I,K)})^T + (\Gamma (\Psi) + I)\epsilon_{64}\Bigl(\sum_{K=1}^{I-1}|\hat{L}_{IK}||\hat{L}_{IK}^T| + |A_{II}|\Bigr),
    \end{aligned}
    \end{equation}
    where $\rho^{(I,K)}$ is the length $n_b$ vector defined as $\rho^{(I,K)}_s = \max_{(K-1)n_b < j \le Kn_b}|\hat{l}_{sj}|$ for $(I-1)n_b < s \le In_b$. \jwdnote{Nearly the same as above inequality, combine them?} \svknote{Thanks, done}
    Combining (\ref{Schur err emulate}) and (\ref{diag cholesky}) yields
    \begin{equation} \label{Chol diag block emu}
    \begin{aligned}
        \Bigl|A_{II} - \sum_{K=1}^I \hat{L}_{IK}\hat{L}_{IK}^T\Bigr| &\le \Theta (\psi,n_b)\sum_{K=1}^{I-1}(\rho^{(I,K)})(\rho^{(I,K)})^T + (\Gamma (\Psi) + I)\epsilon_{64}\Bigl(\sum_{K=1}^{I-1}|\hat{L}_{IK}||\hat{L}_{IK}^T| + |A_{II}|\Bigr) + n_b\epsilon_{64}\Omega_{II}.
    \end{aligned}
    \end{equation}
Assuming $s < t \le n$, since $\hat{L}\hat{L}^T = A + E$, we can write $\sum_{k=1}^s \hat{l}_{sk}\hat{l}_{tk} = a_{st} + e_{st}$, and  $\sum_{k=1}^t \hat{l}_{tk}^2 = a_{tt} + e_{tt}$. 
Ditto for $s$. Also, each element, $\rho^{(I,K)}_t$, of $\rho^{(I,K)}$ is bounded by $\sqrt{\sum_{k=(K-1)n_b+1}^{Kn_b}\hat{l}^2_{tk}}$, and so $\sum_{K=1}^{I-1}\bigl(\rho^{(I,K)}_t\bigr)^2 \le \sum_{k=1}^{(I-1)n_b}\hat{l}^2_{tk}$. Now zoom in to the diagonal elements of the diagonal blocks of $A$ in Equation (\ref{Chol diag block emu}). Looking at the diagonal element $a_{tt}$, and noting that $I \le n/n_b$, we can write
\jwdnote{Why is first inequality below
not an equality for $|e_{tt}|$?
Are subscripts of $\hat{l}_{kt}$ swapped again?} \svknote{Subscripts fixed.}
\begin{equation} \label{eq:diag_elem_emu}
\begin{aligned}
    |e_{tt}| = \Bigl|a_{tt} - \sum_{k=1}^t\hat{l}^2_{tk}\Bigr| &\le \Theta(\psi, n_b)(|a_{tt}| + |e_{tt}|) + \left(\Gamma (\Psi) + \frac{n}{n_b}\right)\epsilon_{64}\Bigl(\sum_{k=1}^{(I-1)n_b}\hat{l}^2_{tk} + |a_{tt}|\Bigr) + n_b\epsilon_{64}|\hat{a}^u_{tt}|.
\end{aligned}
\end{equation}
The entry $(t,t)$ of (\ref{diag cholesky}) is $|\hat{a}^u_{tt} - \sum_{k=(I-1)n_b+1}^{In_b}\hat{l}^2_{tk}| \le n_b\epsilon_{64}|\hat{a}^u_{tt}|$, so $|\hat{a}^u_{tt}| \le \frac{1}{1 - n_b\epsilon_{64}}\sum_{k=(I-1)n_b+1}^{In_b}\hat{l}^2_{tk}$. Plugging this into \eqref{eq:diag_elem_emu},
\jwdnote{Not clear what this refers to?} \svknote{Thank you. Wrote everything out explicitly}
\begin{equation*}
\begin{aligned}
    |e_{tt}| &\le \Theta(\psi, n_b)(|a_{tt}| + |e_{tt}|) + \left(\Gamma (\Psi) + \frac{n}{n_b}\right)\epsilon_{64}\Bigl(\sum_{k=1}^{(I-1)n_b}\hat{l}^2_{tk} \Bigr) + \left(\Gamma (\Psi) + \frac{n}{n_b}\right)\epsilon_{64}|a_{tt}| + \frac{n_b\epsilon_{64}}{1 - n_b\epsilon_{64}}\Bigl(\sum_{k=(I-1)n_b + 1}^{In_b}\hat{l}^2_{tk} \Bigr) \\
    \implies |e_{tt}| &\le \Theta(\psi, n_b)(|a_{tt}| + |e_{tt}|) + M\epsilon_{64}\Bigl(\sum_{k=1}^{(I-1)n_b}\hat{l}^2_{tk} \Bigr) + M\epsilon_{64}|a_{tt}| + M\epsilon_{64}\Bigl(\sum_{k=(I-1)n_b + 1}^{In_b}\hat{l}^2_{tk} \Bigr)
\end{aligned}
\end{equation*}
\begin{equation*}
\begin{aligned}
    \implies |e_{tt}| &\le \bigl(\Theta(\psi, n_b) + M\epsilon_{64}\bigr)(|a_{tt}| + |e_{tt}|) + M\epsilon_{64}|a_{tt}|.
\end{aligned}
\end{equation*}
 Thus, for the diagonal elements,
\begin{equation} \label{eq:ett-emu}
    |e_{tt}| \le \frac{\Theta (\psi,n_b) + 2M\epsilon_{64}}{1 - \Theta (\psi,n_b) - M\epsilon_{64}}|a_{tt}|.
\end{equation}
The same inequality holds for all indices from 1 to $n$.
For off-diagonal elements in the diagonal block, \eqref{Chol diag block emu} becomes
\jwdnote{Should $a_{st}$ be $|a_{st}|$?} \svknote{Changed below}
\jwdnote{Should the first inequality below be an equality? And I think this
is referred to later as "(4.2)", so it
needs an equation number.} \svknote{Yes on the equality, and changed from equation* to equation}
\begin{equation} 
\begin{aligned}
    |e_{st}| &= \Bigl|a_{st} - \sum_{k=1}^s\hat{l}_{sk}\hat{l}_{tk}\Bigr| \le \Theta (\psi, n_b)\sum_{K=1}^{I-1}\rho^{(I,K)}_s\rho^{(I,K)}_t + \left(\Gamma(\Psi) + \frac{n}{n_b}\right)\epsilon_{64}\Bigl(\sum_{k=1}^{(I-1)n_b}|\hat{l}_{sk}||\hat{l}_{tk}| + |a_{st}|\Bigr) + n_b\epsilon_{64}\,\omega_{st}
    \end{aligned}
    \end{equation}
First, using Cauchy-Schwarz, we can bound $\sum_{K=1}^{I-1}\rho^{(I,K)}_s\rho^{(I,K)}_t \le \sqrt{\sum_{k=1}^{t}\hat{l}_{tk}^2}\sqrt{\sum_{k=1}^{s}\hat{l}_{sk}^2} \le \sqrt{|a_{ss}| + |e_{ss}|}\sqrt{|a_{tt}| + |e_{tt}|}$. So, 
\begin{equation} \label{eq: e_st diag emu intermediate} 
    |e_{st}| \le \Theta (\psi, n_b)\sqrt{(|a_{tt}| + |e_{tt}|)(|a_{ss}| + |e_{ss}|)} + \left(\Gamma(\Psi) + \frac{n}{n_b}\right)\epsilon_{64}\Bigl(\sum_{k=1}^{(I-1)n_b}|\hat{l}_{sk}||\hat{l}_{tk}| + |a_{st}|\Bigr) + n_b\epsilon_{64}\,\omega_{st}.
\end{equation}
Applying Cauchy-Schwarz to $\sum_{k=1}^{(I-1)n_b}|\hat{l}_{sk}||\hat{l}_{tk}|$ and using $|\hat{a}^u_{tt}| \le \frac{1}{1 - n_b\epsilon_{64}}\sum_{k=(I-1)n_b+1}^{In_b}\hat{l}^2_{tk}$  we can simplify the second and third terms in the right hand side of \eqref{eq: e_st diag emu intermediate} as follows:\jwdnote{What
happened to the $a_{st}$ term?} \svknote{Added it in}
\begin{equation*}
\begin{aligned}
    &\left(\Gamma(\Psi) + \frac{n}{n_b}\right)\epsilon_{64} \left (\sum_{k=1}^{(I-1)n_b}|\hat{l}_{sk}||\hat{l}_{tk}| + |a_{st}|\right )\;+\; n_b\epsilon_{64}\,\omega_{st} \le\;
M\epsilon_{64}\sqrt{\sum_{k=1}^{(I-1)n_b}\hat{l}^2_{sk}\sum_{k=1}^{(I-1)n_b}\hat{l}^2_{tk}} + M\epsilon_{64}|a_{st}| \\
    &\qquad\qquad + M\epsilon_{64}\sqrt{\sum_{k=(I-1)n_b+1}^{In_b}\hat{l}^2_{sk}\sum_{k=(I-1)n_b+1}^{In_b}\hat{l}^2_{tk}}.
\end{aligned}
\end{equation*}
\jwdnote{The second square root, is this
supposed to be a bound on $|a_{st}|$? If so,
why does the sum start at $k=(I-1)n_b +1$ instead
of $k=1$?} \svknote{I switched the second and third terms, switched them back to the same order now.}
For non-negative reals $x_1,x_2,y_1,y_2$, we know that $\sqrt{x_1y_1} + \sqrt{x_2y_2} \le \sqrt{(x_1 + x_2)(y_1 + y_2)}$. Thus, we can merge $\hat{l}_{tk}$ and $\hat{l}_{sk}$ \jwdnote{${\hat{l}}_{tk}$ and ${\hat{l}}_{sk}$?} \svknote{fixed}
sums and upper bound the above expression by
\begin{equation*}
\begin{aligned}
M\epsilon_{64}\sqrt{\sum_{k=1}^{In_b}\hat{l}^2_{sk}\sum_{k=1}^{In_b}\hat{l}^2_{tk}} + M\epsilon_{64}|a_{st}| \;\le\; M\epsilon_{64}\sqrt{(|a_{ss}| + |e_{ss}|)(|a_{tt}| + |e_{tt}|)} + M\epsilon_{64}|a_{st}|.
\end{aligned}
\end{equation*}
Plugging this into (\ref{eq: e_st diag emu intermediate}), \jwdnote{(4.2) does not appear.} \svknote{Thanks, fixed} we get
\begin{equation} \label{diagblock offdiag emu-intermediate}
\begin{aligned}
    |e_{st}| &\le \left(\Theta (\psi,n_b) + M\epsilon_{64}\right)\sqrt{(|a_{ss}| + |e_{ss}|)(|a_{tt}| + |e_{tt}|)} + M\epsilon_{64}|a_{st}|, 
\end{aligned}
\end{equation}
Then, to bound $|a_{st}|$, note that $|a_{st}| \le \sum_{k=1}^s |\hat{l}_{tk}||\hat{l}_{sk}| + |e_{st}| \le \sqrt{(|a_{ss}| + |e_{ss}|)(|a_{tt}| + |e_{tt}|)} + |e_{st}|$. Substituting this into \eqref{diagblock offdiag emu-intermediate}, 
\begin{equation*}
\begin{aligned}
    |e_{st}| &\le  \left(\Theta (\psi,n_b) + 2M\epsilon_{64}\right)\sqrt{(|a_{ss}| + |e_{ss}|)(|a_{tt}| + |e_{tt}|)} + M\epsilon_{64}|e_{st}| \\
    \implies |e_{st}| &\le \frac{\Theta (\psi,n_b) + 2M\epsilon_{64}}{1 - M\epsilon_{64}}\sqrt{(|a_{ss}| + |e_{ss}|)(|a_{tt}| + |e_{tt}|)}.
\end{aligned}
\end{equation*}
Substituting $e_{ss}, e_{tt}$ from (\ref{eq:ett-emu}) then gives
\jwdnote{Again, should the first inequality below be an equality?} \svknote{Yes, fixed.}
\begin{equation} \label{diagblock offdiag emu}
\begin{aligned}
     |e_{st}| \le \frac{\Theta (\psi,n_b) + 2M\epsilon_{64}}{1 - M\epsilon_{64}} \cdot \frac{1 + M\epsilon_{64}}{1 - \Theta (\psi,n_b) - M\epsilon_{64}}\sqrt{|a_{ss}||a_{tt}|} \le \frac{\Theta (\psi,n_b) + 2M\epsilon_{64}}{1 - M\epsilon_{64}} \cdot \frac{1 - M\epsilon_{64}}{1 - \Theta (\psi,n_b) - 3M\epsilon_{64}}\sqrt{|a_{ss}||a_{tt}|} 
     \\
     = \frac{\Theta (\psi,n_b) + 2M\epsilon_{64}}{1 - \Theta (\psi,n_b) - 3M\epsilon_{64}}\sqrt{|a_{ss}||a_{tt}|}.
\end{aligned}
\end{equation}

To compute the error bound for the off-diagonal blocks, we must use the Ozaki 1 TRSM bound (\ref{b - lk final}) rather than the Cholesky error bound on the diagonal block, and we use \eqref{eqn:GEMMerrbnd} for the error bound on GEMM. Let $J < I \le n/n_b$. Bounding $|A_{IJ} - \sum_{K=1}^J\hat{L}_{IK}\hat{L}_{JK}^T|$ by the Schur-update error $|A_{IJ} - \sum_{K=1}^{J-1}\hat{L}_{IK}\hat{L}_{JK}^T - \hat{A}^u_{IJ}|$ plus the TRSM error $|\hat{A}^u_{IJ} - \hat{L}_{IJ}\hat{L}_{JJ}^T|$, we get 
\begin{equation}
\label{eq:off-diag-emu-block}
\begin{aligned}
    \Bigl|A_{IJ} - \sum_{K=1}^J \hat{L}_{IK}\hat{L}_{JK}^T\Bigr| &\le \sum_{K=1}^{J-1}\Theta (\pi_g(K), n_b) (\rho^{(I,K)})(\rho^{(J,K)})^T + \bigl(\max_{K \le J}\Gamma (\pi_g(K)) + J\bigr)\epsilon_{64}\Bigl(\sum_{K=1}^{J-1}|\hat{L}_{IK}||\hat{L}_{JK}^T| + |A_{IJ}|\Bigr) \\
    &\qquad + \frac{1}{1 - (\Gamma (\pi_t(J)) + n_b)\epsilon_{64}}\Bigl(\Theta (\pi_t(J),n_b)(\rho^{(I,J)})(\rho^{(J,J)})^T + 2(\Gamma (\pi_t(J)) + n_b)\epsilon_{64}|\hat{L}_{IJ}||\hat{L}_{JJ}^T|\Bigr) \\
    &\le \Theta (\psi, n_b)\sum_{K=1}^{J-1}(\rho^{(I,K)})(\rho^{(J,K)})^T + (\Gamma (\Psi) + J)\epsilon_{64}\Bigl(\sum_{K=1}^{J-1}|\hat{L}_{IK}||\hat{L}_{JK}^T| + |A_{IJ}|\Bigr) \\
    &\qquad + \frac{1}{1 - (\Gamma (\Psi) + n_b)\epsilon_{64}}\Bigl(\Theta (\psi,n_b)(\rho^{(I,J)})(\rho^{(J,J)})^T + 2(\Gamma (\Psi) + n_b)\epsilon_{64}|\hat{L}_{IJ}||\hat{L}_{JJ}^T|\Bigr)
\end{aligned}
\end{equation}
\jwdnote{Refer to equation number used in
rest of bound too.}
where the second line is (\ref{b - lk final}) applied to solving $\hat{L}_{IJ}\hat{L}_{JJ}^T = \hat{A}^u_{IJ}$ with some block size $n_b' \le n_b$. Note that although the theorem is proved for $LX = B$, it also holds for $XL^T = B$ (with appropriate reordering of terms) since all the matrices are $n_b \times n_b$ in this case. Now, we want to merge the error contributed by TRSM and the error contributed by GEMM/SYRK. Using $\Theta (\psi, n_b)\sum_{K=1}^{J-1}(\rho^{(I,K)})(\rho^{(J,K)})^T + \frac{\Theta (\psi, n_b)}{1 - M\epsilon_{64}}(\rho^{(I,J)})(\rho^{(J,J)})^T \le \frac{\Theta (\psi,n_b)}{1 - M\epsilon_{64}}\sum_{K=1}^{J}(\rho^{(I,K)})(\rho^{(J,K)})^T$, and noting that $\Gamma (\Psi) + J \le \Gamma (\Psi) + \frac{n}{n_b} \le M$, $\Gamma (\Psi) + n_b \le M$, and $M\epsilon_{64} \le \frac{2M\epsilon_{64}}{1 - M\epsilon_{64}}$, we can rewrite \eqref{eq:off-diag-emu-block} as 
\begin{equation}\label{offdiag block emu}
    \Bigl|A_{IJ} - \sum_{K=1}^J \hat{L}_{IK}\hat{L}_{JK}^T\Bigr| \le \frac{\Theta (\psi, n_b)}{1 - M\epsilon_{64}}\sum_{K=1}^{J}(\rho^{(I,K)})(\rho^{(J,K)})^T + \frac{2M\epsilon_{64}}{1 - M\epsilon_{64}}\sum_{K=1}^{J}|\hat{L}_{IK}||\hat{L}_{JK}^T| + M\epsilon_{64}|A_{IJ}|.
\end{equation}
Again, the element $(s,t)$ of $\sum_{K=1}^{J}|\hat{L}_{IK}||\hat{L}_{JK}^T|$, and $\sum_{K=1}^{J}(\rho^{(I,K)})(\rho^{(J,K)})^T$  can both be bounded by $\sqrt{(|a_{ss}| + |e_{ss}|)(|a_{tt}| + |e_{tt}|)}$ via Cauchy-Schwarz. Also, as done after \eqref{diagblock offdiag emu-intermediate}, we can bound $|a_{st}| \le \sqrt{(|a_{ss}| + |e_{ss}|)(|a_{tt}| + |e_{tt}|)} + |e_{st}|$  \jwdnote{Why not $|e_{ss}|$ and $|e_{tt}|$
as in other bound?}\svknote{Because I forgot the $||$ again :). Sorry}
. Hence, the element version of \eqref{offdiag block emu} is
\begin{equation}\label{eq:elemwise-offdiag-bound}
\begin{aligned} 
    |e_{st}| = \Bigl|a_{st} - \sum_{k=1}^t\hat{l}_{sk}\hat{l}_{tk}\Bigr| &\le \frac{\Theta (\psi,n_b)}{1 - M\epsilon_{64}}\sum_{K=1}^{J}\rho^{(I,K)}_s\rho^{(J,K)}_t + \frac{2M\epsilon_{64}}{1 - M\epsilon_{64}}\sum_{k=1}^{t}|\hat{l}_{sk}||\hat{l}_{tk}| + M\epsilon_{64}|a_{st}| \\
    &\le \frac{\Theta (\psi,n_b) + 2M\epsilon_{64}}{1 - M\epsilon_{64}}\sqrt{(|a_{ss}| + |e_{ss}|)(|a_{tt}| + |e_{tt}|)} + M\epsilon_{64}\Bigl(\sqrt{(|a_{ss}| + |e_{ss}|)(|a_{tt}| + |e_{tt}|)} + |e_{st}|\Bigr).
\end{aligned}
\end{equation}
Using (\ref{eq:ett-emu}) on rows $s$ and $t$,
\begin{equation*}
    (1 - M\epsilon_{64})|e_{st}| \le \frac{\Theta (\psi,n_b) + 3M\epsilon_{64} - M^2\epsilon_{64}^2}{1 - M\epsilon_{64}}\sqrt{(|a_{ss}| + |e_{ss}|)(|a_{tt}| + |e_{tt}|)} \le \frac{\Theta (\psi,n_b) + 3M\epsilon_{64} - M^2\epsilon_{64}^2}{1 - \Theta (\psi,n_b) - 3M\epsilon_{64}}\sqrt{|a_{ss}||a_{tt}|}
\end{equation*}
\begin{equation} \label{eq:final_chol_emu}
\begin{aligned}
    \implies |e_{st}| \le \frac{\Theta (\psi,n_b) + 3M\epsilon_{64} - M^2\epsilon_{64}^2}{1 - M\epsilon_{64}} \cdot \frac{1}{1 - \Theta (\psi,n_b) - 3M\epsilon_{64}}\sqrt{|a_{ss}||a_{tt}|} \\
    \le \frac{\Theta (\psi,n_b) + 4M\epsilon_{64}}{1 - \Theta (\psi,n_b) - 3M\epsilon_{64}}\sqrt{|a_{ss}||a_{tt}|} \le \frac{\gamma_{\mathrm{emu}}(\psi, \Psi, n, n_b)}{1 - \gamma_{\mathrm{emu}}(\psi, \Psi, n, n_b)}\sqrt{|a_{ss}||a_{tt}|},
\end{aligned}
\end{equation}
where the last step follows from  $(\Theta (\psi,n_b) + 4M\epsilon_{64})(1 - M\epsilon_{64}) - (\Theta (\psi,n_b) + 3M\epsilon_{64} - M^2\epsilon_{64}^2) = M\epsilon_{64}(1 - \Theta (\psi,n_b) - 3M\epsilon_{64}) \ge 0$ since $\Theta (\psi,n_b) + 3M\epsilon_{64} \le \Theta (\psi,n_b) + 4M\epsilon_{64} \le 1$ by assumption.

Since the bound in \eqref{eq:final_chol_emu} is looser than the bounds we derived for the elements of the diagonal block in \eqref{eq:ett-emu} and \eqref{diagblock offdiag emu}, it applies there too.
\end{proof}

\jwdnote{Do you mean $\Theta > M$? If $p=7$ than $\Theta = 30 \cdot 2^{-54}n$
and $M \approx 28 n$ for small $n_b$,
which is much larger.} \svknote{Thanks, I think it should be p < 7}

\subsection{Mixed-precision Floating-point case}

\begin{theorem} \label{thm:chol-fp}
Suppose that the factorization in Algorithm~\ref{A_chol} uses (potentially diagonally scaled or MX) GEMM with inputs in precision $\pi_g(I)$ in iteration $I$, with machine epsilon $\epsilon_{\pi_g(I)}$ and underflow threshold $\mu_{\pi_g(I)}$, accumulation in precision $h$, as well as block diagonal factorization and TRSM in precision $h$. Here $\pi_g(I)$ is a lower precision, such as FP8 or FP16, and $h$ is a higher precision, such as FP64 or FP32. Let $\psi$ be the format that has the largest machine epsilon in all iterations, and let $\Psi$ be the format with the largest underflow threshold in all iterations. Assume that the algorithm completes successfully and that there is no overflow in both low and high precision. If rounding blocks of $\hat{L}$ (before GEMM) results in underflow on at most $\nu<n$ entries of
$\hat{L}$ per row, and provided $n_b\eh<1$ and $\gfp<1$, then $A+E=\hat L\hat L^T$ with
\begin{equation} \label{eq:chol-fp}
  |e_{st}| \;\le\; \frac{\gfp}{1-\gfp}\sqrt{a_{ss} + \nu \mu_\Psi^2 d_{ss}^2}\sqrt{a_{tt} + \nu \mu_\Psi^2d_{tt}^2},
\end{equation}
where
\begin{equation} \label{eq:gamma-fp}
  \gfp = \gfp(\epsilon_\psi, n, n_b) = (2\epsilon_\psi + \epsilon_\psi^2) + 3C , \qquad C = \eh\max\left(\frac{n_b}{1 - n_b\eh},\, \left(n_b + 2 + \frac{n}{n_b}\right)(1+\epsilon_\psi)^2\right),
\end{equation}
$n_b$ is the block size used for factorization, $\epsilon_\psi$ is the machine epsilon of the format $\psi$, $\mu_\Psi$ is the smallest normal number in the low precision format $\Psi$, and $d_{rr}$ is the largest scale factor we encounter in row $r$ of $\hat{L}$ when scaling $\hat{L}$ before rounding to low precision. Concretely, if $\sigma_{ij}$ is the scale factor used when scaling $\hat{l}_{ij}$ to $\hat{l}_{ij}/\sigma_{ij}$ before rounding, followed by GEMM, $d_{rr} = 0$ for $r \le n_b$, and  $d_{rr} = \max_{k=1}^{n_b (\lceil \frac{r}{n_b} \rceil - 1)} \sigma_{rk}$ for $n_b < r \le n$ . Additionally, if we set the microscaling target of Algorithm~\ref{alg:microscaling} to $\beta$, and store the scale factor with some precision $\mathrm{scal}$, then, provided $\gfp(1 + \lambda_\Psi) < 1$, we can also say
\begin{equation} \label{eq:chol-fp-mx}
    |e_{st}| \le \frac{\gfp(1 + \lambda_\Psi)}{1 - \gfp(1 + \lambda_\Psi)} \sqrt{a_{ss}a_{tt}},
\end{equation}
where $\lambda_\Psi = \frac{\nu \mu_\Psi^2(1 + \epsilon_{\mathrm{scal}})^2}{\beta^2}$, and $\epsilon_{\mathrm{scal}}$ is the machine epsilon of the precision used to store the scale factor.
\end{theorem}

\begin{proof}
For simplicity, assume $n_b|n$.
Following the same structure as the proof of Theorem~\ref{thm:chol-emu}
\jwdnote{Say "proof of Thm. 4.1", again avoiding
indirect references.}, the roundoff error on diagonal block is the sum of the error from the block diagonal Cholesky $n_b\epsilon_h\Omega_{II}$ and the error from the MX GEMM given by Equation (\ref{eq:MXGEMM}):
\begin{equation}
\begin{aligned}
  \Bigl|A_{II}-\sum_{K=1}^{I}\hat L_{IK}\hat L_{IK}^T\Bigr|
  &\;\le\; n_b\epsilon_h\,\Omega_{II} + \left(n_b + 2 + \frac{n}{n_b}\right)\epsilon_h|A_{II}| + \sum_{K=1}^{I-1}\Bigl[\epsilon_{\pi_g(K)}\bigl(|\hat L''_{IK}||\hat L^T_{IK}| + |\hat L_{IK}||\hat L''^T_{IK}|
     + \epsilon_{\pi_g(K)}|\hat L''_{IK}||\hat L''^T_{IK}|\bigr) \\
     &\qquad + \left(n_b + 2 + \frac{n}{n_b}\right)\epsilon_h \bigl(|\hat L_{IK}|+\epsilon_{\pi_g(K)}|\hat L''_{IK}|\bigr)\bigl(|\hat L_{IK}|+\epsilon_{\pi_g(K)}|\hat L''_{IK}|\bigr)^T\Bigr] \\
  &\;\le\; n_b\epsilon_h\,\Omega_{II} + \left(n_b + 2 + \frac{n}{n_b}\right)\epsilon_h|A_{II}| + \sum_{K=1}^{I-1}\Bigl[\epsilon_\psi\bigl(|\hat L''_{IK}||\hat L^T_{IK}| + |\hat L_{IK}||\hat L''^T_{IK}|
     + \epsilon_\psi|\hat L''_{IK}||\hat L''^T_{IK}|\bigr) \\
     &\qquad + \left(n_b + 2 + \frac{n}{n_b}\right)\epsilon_h \bigl(|\hat L_{IK}|+\epsilon_\psi|\hat L''_{IK}|\bigr)\bigl(|\hat L_{IK}|+\epsilon_\psi|\hat L''_{IK}|\bigr)^T\Bigr].
\end{aligned}
  \label{eq:diagblock_FPchol}
\end{equation}
We have applied Theorem \ref{thm:MXGEMM}  with $A = \hat L_{IK}$, $B = \hat L_{IK}^T$ and $\phi \approx \epsilon_{\pi_g(K)}$ as in \eqref{eq:MXGEMMsimple}. Also, $|\hat L''_{IK}| = D_{IK}M_{IK} + |\hat{L}'_{IK}|$ with $D_{IK}$ a diagonal matrix holding the largest scale factor of each row, $\hat{L}'_{IK} = \hat{L}_{IK}$ in positions where $\hat{L}_{IK}$ doesn't underflow when rounding to low precision, and zero everywhere else
\jwdnote{and zero otherwise}, and $M_{IK}$ a matrix with $\mu_{\pi_g(K)}$ in positions where there is underflow, and zero everywhere else (same as Lemma \ref{lem:MXmatrix}). Note that \eqref{eq:diagblock_FPchol} has accumulation error $\left(n_b + 2 + \frac{n}{n_b}\right)\epsilon_h (|A_{II}| + \sum_{K=1}^{I-1}(|\hat{L}_{IK}| + \epsilon_\psi|\hat{L}''_{IK}|)(|\hat{L}_{IK}| + \epsilon_\psi|\hat{L}''_{IK}|)^T)$, which is slightly different from the accumulation error in $\eqref{eq:MXGEMM}$. 
\jwdnote{Ambiguous, there are three $\epsilon_h$ terms.}
To understand why the accumulation error has the coefficient $\left(n_b + 2 + \frac{n}{n_b}\right)\epsilon_h$ and the extra $|A_{II}|$ term, 
\jwdnote{Again, there are 2 such terms.}
we can count the roundoff terms the same way we did for Lemmas~\ref{lem:MXGEMM-DS} and~\ref{lem:MXGEMM-SD}. Consider an element of $\hat{L}_{IK}$. When we round $\hat{L}_{IK}$ to low precision and perform MXGEMM $\hat{L}_{IK}\hat{L}^T_{IK}$, each element in that product receives at most $n_b + 2$ roundoff terms according to Theorem \ref{thm:MXGEMM}. Two terms come from the application of scale factors, one from the multiplication, and the remaining $n_b - 1$ terms from the sum reduction. Then, as in the proof of Theorem~\ref{thm:chol-emu}, Algorithm~\ref{A_chol} subtracts the $I - 1$ such products $\hat{L}_{IK}\hat{L}^T_{IK}$ from $A_{II}$, which adds at most $I - 1 \le \frac{n}{n_b}$ roundoff terms.
\jwdnote{Do you need $n - n_b -1 +1$?} \svknote{Not really} This gives a total of at most $n_b + 2 + \frac{n}{n_b}$ roundoff terms on all elements of $\hat{L}_{IK}$ and $A_{II}$. Moving forward, assuming $(s,t)$ falls in block $(I,J)$ with $J \le I$, we will refer to the entry of $\epsilon_\psi \sum_{K=1}^{J-1} \left ( |\hat L''_{IK}||\hat L^T_{JK}| + |\hat L_{IK}||\hat L''^T_{JK}|
     + \epsilon_\psi|\hat L''_{IK}||\hat L''^T_{JK}|\right )$ in position $(s,t)$ as $\zeta^{\mathrm{trunc}}_{st}$, and the corresponding entry of $\left(n_b + 2 + \frac{n}{n_b}\right)\epsilon_h \sum_{K=1}^{J-1} \bigl(|\hat L_{IK}|+\epsilon_\psi|\hat L''_{IK}|\bigr)\bigl(|\hat L_{JK}|+\epsilon_\psi|\hat L''_{JK}|\bigr)^T$ as $\zeta^{\mathrm{L\_accum}}_{st}$. Note that we do not include $\left(n_b + 2 + \frac{n}{n_b}\right)\epsilon_h|A_{II}|$ in $\zeta^{\mathrm{L\_accum}}_{st}$.  Now, for the proof of the Cholesky error bound, we will use the same strategy as the proof of Theorem~\ref{thm:chol-emu}: we first find a bound on the diagonal elements of the diagonal blocks, then the off-diagonal elements of the diagonal blocks, then all elements in off-diagonal blocks.
Defining $\hat{L}''_{IK} = 0$ for $K \ge I$ (since we don't use the corresponding blocks of $L$ for the MX GEMMs) , we can write $\hat L''$ for the full-matrix assembly of all $\hat L''_{IK}$, so that
$\hat l''_{tk}$ is defined for every $1 \le k,t \le n$. \jwdnote{Isn't $\hat l''_{tk}=0$ for $k>t$, which
the above formulas assume?} \svknote{I think I technically haven't defined $\hat{L}''_{IK}$ for $K \ge I$, since I don't do an MX GEMM with them. But I suppose I can just define those to be 0 as well}
Also, as in the statement of the Theorem, let $d_{tt}$ be the maximum scale factor in all blocks that contain the row $t$. If row $t$ falls into block row $I$, $d_{tt} = \max_{K=1}^{I-1} \left(D_{IK}\right)_{t - (I-1)n_b,\, t - (I-1)n_b}$.  As in the proof of Theorem \ref{thm:chol-emu}, $\Omega_{II}$ is the block diagonal Cholesky bound of \cite{demmel1989cholesky}, with entries $\omega_{st} = \sqrt{|\hat{a}^u_{ss}||\hat{a}^u_{tt}|}$, and $\hat{A}^u_{II}$ \jwdnote{Missing subscripts?} \svknote{Added it in}
denotes the computed Schur-updated block.

Now zoom in on the diagonal entry $(t,t)$. The diagonal element version of \eqref{eq:diagblock_FPchol} is 
\begin{equation} \label{eq:diagelem_FPchol}
    |e_{tt}| \le n_b\epsilon_h \omega_{tt} + \left(n_b + 2 + \frac{n}{n_b}\right)\epsilon_h |a_{tt}| + \zeta^{\mathrm{trunc}}_{tt} + \zeta^{\mathrm{L\_accum}}_{tt}.
\end{equation}
We already saw in the proof of Theorem~\ref{thm:chol-emu} that
\begin{equation} \label{eq:diag-bound-au}
   \epsilon_h n_b \omega_{tt} = \epsilon_h n_b |\hat{a}^u_{tt}| \le \frac{\epsilon_h n_b}{1 - n_b \epsilon_h}\sum_{k=(I - 1)n_b + 1}^{In_b} \hat{l}^2_{tk}.
\end{equation}
\jwdnote{Actually, I think in Thm 4.1 the limit of the 
sum was $In_b$ instead of $t$, but of course the 
difference is just a sum of zeros.} \svknote{Thanks, changed to $In_b$ for consistency, but yes, it shouldn't matter}
We will now focus on $\zeta^{\mathrm{trunc}}_{tt}$. For $\epsilon_\psi \sum_{K=1}^{I - 1} |\hat{L}''_{IK}||\hat{L}^T_{IK}|$, first define $U_t \subseteq \{1,\ldots,(I-1)n_b\}$ as the set of $k$ for which $\hat l_{tk}$ underflows when rounded to low precision. Now, for $k \notin U_t$, we have $|\hat{l}''_{tk}| = |\hat{l}_{tk}|$, and for $k \in U_t$, $|\hat{l}_{tk}| \le (1 + \epsilon_h)d_{tt}\mu_\Psi$ and $|\hat{l}''_{tk}| \le d_{tt}\mu_\Psi$. \jwdnote{Remind reader other entries are zero.} The rest of the $\hat{l}''_{tk}$ with $k > (I - 1)n_b$ are all 0.
 Element $(t,t)$ of $\epsilon_\psi \sum_{K=1}^{I-1}|\hat{L}_{IK}''||\hat{L}_{IK}^T|$ is $\epsilon_\psi\sum_{k=1}^{(I-1)n_b} |\hat{l}''_{tk}||\hat{l}_{tk}|$. We will try to bound this by applying Cauchy-Schwarz, for which we need a bound on $\sum_{k=1}^{(I-1)n_b} |\hat{l}''_{tk}|^2$
\begin{equation} \label{uuc-split-diagonal-sum}
    \sum_{k=1}^{(I-1)n_b} |\hat{l}''_{tk}|^2
    = \sum_{k \in U_t}|\hat{l}''_{tk}|^2 + \sum_{k \notin U_t} |\hat{l}''_{tk}|^2 \le \nu \mu_\Psi^2 d_{tt}^2 +  \sum_{k=1}^{(I - 1)n_b}|\hat{l}_{tk}|^2.
\end{equation}
  Then
\begin{equation} \label{diag-elem2-bound}
    \epsilon_\psi \sum_{k=1}^{(I-1)n_b}|\hat{l}''_{tk}||\hat{l}_{tk}| \le \epsilon_\psi \sqrt{\sum_{k=1}^{(I-1)n_b}|\hat{l}''_{tk}|^2 \sum_{k=1}^{(I-1)n_b}\hat{l}^2_{tk}} \le \epsilon_\psi (\sum_{k=1}^{(I-1)n_b}\hat{l}^2_{tk} + \nu \mu_\Psi^2 d_{tt}^2)
\end{equation}
The next matrix term in $\zeta^{\mathrm{trunc}}_{tt}$ is $\epsilon_\psi \sum_{K=1}^{I - 1} |\hat{L}_{IK}||\hat{L}''^T_{IK}|$, which is the transpose of $\epsilon_\psi \sum_{K=1}^{I - 1} |\hat{L}''_{IK}||\hat{L}^T_{IK}|$. So, \eqref{diag-elem2-bound} works as an upper bound for the diagonal element of $\epsilon_\psi \sum_{K=1}^{I - 1} |\hat{L}_{IK}||\hat{L}''^T_{IK}|$ as well.  For $\epsilon_\psi^2\sum_{K=1}^{I - 1} |\hat{L}''_{IK}||\hat{L}''^T_{IK}| $, the diagonal term is $\epsilon_\psi^2 \sum_{k=1}^{(I-1)n_b} |\hat{l}''_{tk}|^2$, which we can get an upper bound on by multiplying \eqref{uuc-split-diagonal-sum} by $\epsilon_\psi^2$. 
So 
\begin{equation} \label{eq:trunc_diag_elem}
    \zeta^{\mathrm{trunc}}_{tt} \le (2\epsilon_\psi + \epsilon_\psi^2)\left ( \sum_{k=1}^{(I-1)n_b}\hat{l}^2_{tk} + \nu \mu_\Psi^2 d_{tt}^2 \right )
\end{equation}
To get a bound on $\zeta^{\mathrm{L\_accum}}_{tt}$, observe the accumulation error matrices in \eqref{eq:diagblock_FPchol}
\begin{equation} \label{eq:expanded_L_acc_err}
\begin{aligned}
    &\left(n_b + 2 + \frac{n}{n_b}\right)\epsilon_h \bigl(|\hat{L}_{IK}| + \epsilon_\psi |\hat{L}_{IK}''|\bigr)\bigl(|\hat{L}_{IK}| + \epsilon_\psi |\hat{L}_{IK}''|\bigr)^T \\
    &\qquad = \left(n_b + 2 + \frac{n}{n_b}\right)\epsilon_h \bigl(|\hat{L}_{IK}||\hat{L}^T_{IK}| + \epsilon_\psi |\hat{L}_{IK}||\hat{L}''^T_{IK}| + \epsilon_\psi |\hat{L}''_{IK}||\hat{L}^T_{IK}| + \epsilon_\psi^2 |\hat{L}''_{IK}||\hat{L}''^T_{IK}|\bigr) .
\end{aligned}
\end{equation}
On the right hand side of \eqref{eq:expanded_L_acc_err}, \eqref{diag-elem2-bound} and \eqref{uuc-split-diagonal-sum} already bound the diagonal elements of every term except $|\hat{L}_{IK}||\hat{L}^T_{IK}|$, which has diagonal element $\sum_{k=1}^{(I-1)n_b}\hat{l}^2_{tk}$. So, the diagonal element of the accumulation error matrices can be bounded by 
\begin{equation} \label{eq:accum_diag_elem}
\begin{aligned}
   \zeta^{\mathrm{L\_accum}}_{tt} &\le \left(n_b + 2 + \frac{n}{n_b}\right)\epsilon_h \Bigl(\sum_{k = 1}^{(I-1)n_b} \hat{l}^2_{tk} + (2\epsilon_\psi + \epsilon_\psi^2)\Bigl(\sum_{k = 1}^{(I-1)n_b} \hat{l}^2_{tk} + \nu \mu_\Psi^2 d_{tt}^2\Bigr)\Bigr) \\
   &\le \left(n_b + 2 + \frac{n}{n_b}\right)\epsilon_h(1 + \epsilon_\psi)^2\Bigl(\sum_{k = 1}^{(I-1)n_b} \hat{l}^2_{tk} + \nu \mu_\Psi^2 d_{tt}^2\Bigr).
\end{aligned}
\end{equation}
Now, we can rewrite \eqref{eq:diagelem_FPchol} using \eqref{eq:accum_diag_elem},  \eqref{eq:trunc_diag_elem} and \eqref{eq:diag-bound-au} as
\begin{equation} \label{eq:fp-ett}
\begin{aligned}
    |e_{tt}| \le \Bigl[(2\epsilon_\psi + \epsilon_\psi^2) + \left(n_b + 2 + \frac{n}{n_b}\right)\epsilon_h(1+\epsilon_\psi)^2\Bigr]
     \Bigl(\sum_{k=1}^{(I-1)n_b}\hat{l}^2_{tk} + \nu \mu_\Psi^2 d_{tt}^2\Bigr) + \frac{n_b\epsilon_h}{1 - n_b\epsilon_h}\sum_{k=(I-1)n_b+1}^{t}\hat{l}^2_{tk}
      + \left(n_b + 2 + \frac{n}{n_b}\right)\epsilon_h\,|a_{tt}|. 
\end{aligned}
\end{equation}
\begin{equation*}
\begin{aligned}
    \implies |e_{tt}| &\le \Bigl[(2\epsilon_\psi + \epsilon_\psi^2) + \max\Bigl\{\left(n_b + 2 + \frac{n}{n_b}\right)\eh(1+\epsilon_\psi)^2 ,\; \frac{n_b\eh}{1 - n_b\eh}\Bigr\}\Bigr](|a_{tt}| + |e_{tt}| + \nu \mu_\Psi^2 d_{tt}^2) + \left(n_b + 2 + \frac{n}{n_b}\right)\eh\,|a_{tt}| \\
    &\le \bigl[(2\epsilon_\psi + \epsilon_\psi^2) + C\bigr](|a_{tt}| + |e_{tt}| + \nu \mu_\Psi^2 d_{tt}^2) + C\,|a_{tt}| \;\le\; \gfp\,(|a_{tt}| + |e_{tt}| + \nu \mu_\Psi^2 d_{tt}^2) ,
\end{aligned}
\end{equation*}
\begin{equation} \label{eq:ett-fp}
   \implies |e_{tt}| \le \frac{\gfp}{1 - \gfp}\,(|a_{tt}| + \nu \mu_\Psi^2 d_{tt}^2).
\end{equation}
The same inequalities hold for all the indices from $1$ to $n$.

For the off-diagonal elements of the diagonal block, pick some $1 \le s < t$. The inequality we need to analyze is 
\begin{equation} \label{est-diagblock-fp}
    |e_{st}| \le n_b\epsilon_h \omega_{st} + \left(n_b + 2 + \frac{n}{n_b}\right)\epsilon_h|a_{st}| + \zeta^{\mathrm{trunc}}_{st} + \zeta^{\mathrm{L\_accum}}_{st}.
\end{equation}
Starting with $\omega_{st}$, 
\begin{equation} \label{eq:offdiag-bound-au}
    n_b\eh\,\omega_{st} = n_b\eh\sqrt{|\hat{a}^u_{ss}||\hat{a}^u_{tt}|} \le \frac{n_b\eh}{1 - n_b\eh}\sqrt{\sum_{k=(I-1)n_b+1}^{s}\hat{l}^2_{sk}\sum_{k=(I-1)n_b+1}^{t}\hat{l}^2_{tk}}. 
\end{equation}
Next, $|a_{st}| = |\sum_{k=1}^{s}\hat{l}_{sk}\hat{l}_{tk} - e_{st}| \le \sqrt{(|a_{ss}| + |e_{ss}|)(|a_{tt}| + |e_{tt}|)} + |e_{st}|$. Then, for $\zeta^{\mathrm{trunc}}_{st}$, 
consider the element $(s,t)$ of $\epsilon_\psi \sum_{K=1}^{I - 1}|\hat{L}''_{IK}||\hat{L}^T_{IK}|$. This is $\epsilon_\psi \sum_{k=1}^{(I-1)n_b} |\hat{l}''_{sk}||\hat{l}_{tk}|$. This can be bounded above by Cauchy-Schwarz. $\sum_{k=1}^{(I-1)n_b} |\hat{l}''_{sk}|^2 \le \sum_{k=1}^{(I-1)n_b}\hat{l}^2_{sk} + \nu \mu_\Psi^2 d_{ss}^2$ (follows from \eqref{uuc-split-diagonal-sum}), and ditto for $t$, giving us
\begin{equation} \label{eq:fp-offdiag}
\begin{aligned}
    \epsilon_\psi \sum_{k=1}^{(I-1)n_b} |\hat{l}''_{sk}||\hat{l}_{tk}|
    &\le \epsilon_\psi\sqrt{\sum_{k=1}^{(I-1)n_b}|\hat{l}''_{sk}|^2\sum_{k=1}^{(I-1)n_b}\hat{l}^2_{tk}} \le \epsilon_\psi
      \sqrt{\Bigl(\sum_{k=1}^{(I-1)n_b}\hat{l}^2_{sk} + \nu \mu_\Psi^2 d_{ss}^2\Bigr)\Bigl(\sum_{k=1}^{(I-1)n_b}\hat{l}^2_{tk} + \nu \mu_\Psi^2 d_{tt}^2\Bigr)}.
\end{aligned}
\end{equation}

Clearly, swapping the $s$ and $t$ doesn't change \eqref{eq:fp-offdiag}, so the same inequality works for the $(s,t)$ element of
$\epsilon_\psi(\sum_{K=1}^{I - 1}|\hat{L}''_{IK}||\hat{L}^T_{IK}|)^T =  \epsilon_\psi(\sum_{K=1}^{I - 1}|\hat{L}_{IK}||\hat{L}''^T_{IK}|)$. Element $(s,t)$ of $\epsilon_\psi^2(\sum_{K=1}^{I-1}|\hat{L}''_{IK}||\hat{L}''^T_{IK}|)$ is $\epsilon_\psi^2 \sum_{k=1}^{(I-1)n_b} |\hat{l}''_{sk}||\hat{l}''_{tk}|$, so by the same Cauchy-Schwarz argument, \eqref{eq:fp-offdiag} holds for it too after we multiply an extra $\epsilon_\psi$, which means 
\begin{equation} \label{eq:zeta-trunc-offdiag}
    \zeta^{\mathrm{trunc}}_{st} \le (2\epsilon_\psi + \epsilon_\psi^2) \sqrt{\Bigl(\sum_{k=1}^{(I-1)n_b}\hat{l}^2_{sk} + \nu \mu_\Psi^2 d_{ss}^2\Bigr)\Bigl(\sum_{k=1}^{(I-1)n_b}\hat{l}^2_{tk} + \nu \mu_\Psi^2 d_{tt}^2\Bigr)}.
\end{equation}
Using \eqref{eq:fp-offdiag}, we can construct an upper bound on $\zeta^{\mathrm{L\_accum}}_{st}$ in the same way as we constructed \eqref{eq:accum_diag_elem}, 
\begin{equation}\label{eq:zeta-acc-offdiag}
    \zeta^{\mathrm{L\_accum}}_{st} \le \left(n_b + 2 + \frac{n}{n_b}\right)\epsilon_h (1 + \epsilon_\psi)^2 \sqrt{\Bigl(\sum_{k=1}^{(I-1)n_b}\hat{l}^2_{sk} + \nu \mu_\Psi^2 d_{ss}^2\Bigr)\Bigl(\sum_{k=1}^{(I-1)n_b}\hat{l}^2_{tk} + \nu \mu_\Psi^2 d_{tt}^2\Bigr)}.
\end{equation}
Substituting \eqref{eq:zeta-acc-offdiag}, \eqref{eq:zeta-trunc-offdiag}, \eqref{eq:offdiag-bound-au} in \eqref{est-diagblock-fp} we get 
\begin{equation*}
\begin{aligned}
    |e_{st}| &\le \Bigl[(2\epsilon_\psi + \epsilon_\psi^2) + \left(n_b + 2 + \frac{n}{n_b}\right)\eh(1+\epsilon_\psi)^2\Bigr]
      \sqrt{\textstyle\Bigl(\sum_{k=1}^{(I-1)n_b}\hat{l}^2_{sk} + \nu \mu_\Psi^2 d_{ss}^2\Bigr)\Bigl(\sum_{k=1}^{(I-1)n_b}\hat{l}^2_{tk} + \nu \mu_\Psi^2 d_{tt}^2\Bigr)} \\
    &\qquad + \frac{n_b\eh}{1 - n_b\eh}
      \sqrt{(\textstyle\sum_{k=(I-1)n_b+1}^{s}\hat{l}^2_{sk})(\sum_{k=(I-1)n_b+1}^{t}\hat{l}^2_{tk})}
      + \left(n_b + 2 + \frac{n}{n_b}\right)\eh\Bigl(\sqrt{(|a_{ss}| + |e_{ss}|)(|a_{tt}| + |e_{tt}|)} + |e_{st}|\Bigr) \\
    &\le \bigl[(2\epsilon_\psi + \epsilon_\psi^2) + 2C\bigr]\sqrt{(|a_{ss}| + |e_{ss}| + \nu \mu_\Psi^2 d_{ss}^2)(|a_{tt}| + |e_{tt}| + \nu \mu_\Psi^2 d_{tt}^2)} + C|e_{st}| .
\end{aligned}
\end{equation*}

\begin{equation} \label{eq:est-fp-diagblock}
   \implies  |e_{st}| \le \frac{\gfp - C}{1 - C}\sqrt{(|a_{ss}| + |e_{ss}| + \nu \mu_\Psi^2 d_{ss}^2)(|a_{tt}| + |e_{tt}| + \nu \mu_\Psi^2 d_{tt}^2)}
\end{equation}
\begin{equation} \label{eq:est-fp-diagblock-final}
   \implies  |e_{st}| \le \frac{\gfp}{1 - \gfp}\sqrt{(|a_{ss}| + \nu \mu_\Psi^2 d_{ss}^2)(|a_{tt}| + \nu \mu_\Psi^2 d_{tt}^2)}.
\end{equation}
where we have substituted \eqref{eq:ett-fp} for rows $s$ and $t$ and used the fact that $\gfp - C\gfp - \gfp + C > 0$
since $\gfp < 1$.
 
For off-diagonal block $(I,J)$ with $J < I \le n/n_b$, instead of a block diagonal Cholesky, we have a TRSM. Unlike the emulation case we proved in Theorem~\ref{thm:chol-emu}, we do TRSM entirely in precision $h$. So, the error bound is $n_b \epsilon_h |\hat{L}_{IJ}||\hat{L}^T_{JJ}| \le  \frac{n_b\eh}{1 - n_b\eh}|\hat{L}_{IJ}||\hat{L}^T_{JJ}|$. The matrix inequality for the error is then
\begin{equation} \label{eq:offdiag-block-cholFP}
\begin{aligned}
    \Bigl|A_{IJ} - \sum_{K=1}^{J}\hat{L}_{IK}\hat{L}_{JK}^T\Bigr| &\le
      \frac{n_b\eh}{1 - n_b\eh}|\hat{L}_{IJ}||\hat{L}^T_{JJ}|
      + \sum_{K=1}^{J-1}\epsilon_{\pi_g(K)}\bigl(|\hat L''_{IK}||\hat L^T_{JK}| + |\hat L_{IK}||\hat L''^T_{JK}|
      + \epsilon_{\pi_g(K)}|\hat L''_{IK}||\hat L''^T_{JK}|\bigr) \\
    &\qquad + \left(n_b + 2 + \frac{n}{n_b}\right)\epsilon_h\Bigl(\sum_{K=1}^{J-1}\bigl(|\hat L_{IK}|+\epsilon_{\pi_g(K)}|\hat L''_{IK}|\bigr)\bigl(|\hat L_{JK}|+\epsilon_{\pi_g(K)}|\hat L''_{JK}|\bigr)^T + |A_{IJ}|\Bigr) \\
    &\le
      \frac{n_b\eh}{1 - n_b\eh}|\hat{L}_{IJ}||\hat{L}^T_{JJ}|
      + \sum_{K=1}^{J-1}\epsilon_\psi\bigl(|\hat L''_{IK}||\hat L^T_{JK}| + |\hat L_{IK}||\hat L''^T_{JK}|
      + \epsilon_\psi|\hat L''_{IK}||\hat L''^T_{JK}|\bigr) \\
    &\qquad + \left(n_b + 2 + \frac{n}{n_b}\right)\epsilon_h\Bigl(\sum_{K=1}^{J-1}\bigl(|\hat L_{IK}|+\epsilon_\psi|\hat L''_{IK}|\bigr)\bigl(|\hat L_{JK}|+\epsilon_\psi|\hat L''_{JK}|\bigr)^T + |A_{IJ}|\Bigr) .
\end{aligned}
\end{equation}
Using the same $\zeta_{st}$ notation, we write
\begin{equation} \label{eq:est-offdiag-block}
    |e_{st}| \le \text{TRSM error} + \left(n_b + 2 + \frac{n}{n_b}\right)\epsilon_h|a_{st}| + \zeta^{\mathrm{trunc}}_{st} + \zeta^{\mathrm{L\_accum}}_{st}.
\end{equation}
As we just saw, the $(s,t)$ element of the TRSM error is bounded by $\frac{n_b\epsilon_h}{1 - n_b\epsilon_h}\sum_{k = (J-1)n_b+1}^{t}|\hat l_{sk}||\hat l_{tk}|$, which we can further bound by Cauchy-Schwarz
\begin{equation} \label{eq:cholFP-TRSM-err}
    \frac{n_b\epsilon_h}{1 - n_b\epsilon_h}\sum_{k = (J-1)n_b+1}^{t}|\hat l_{sk}||\hat l_{tk}|
    \le \frac{n_b\epsilon_h}{1 - n_b\epsilon_h}\sqrt{\sum_{k = (J-1)n_b+1}^{t}\hat l^2_{sk}\sum_{k = (J-1)n_b+1}^{t}\hat l^2_{tk}}.
\end{equation}
Next, we need an upper bound on $\zeta^{\mathrm{trunc}}_{st}$, for which we can slightly modify \eqref{eq:zeta-trunc-offdiag}. Looking at \eqref{eq:offdiag-block-cholFP}, we can see that the matrices we include in $\zeta^{\mathrm{trunc}}_{st}$ sum from $k=1$ to $k=(J-1)n_b < (I-1)n_b$ rather than $k=1$ to $k=(I - 1)n_b$ (as seen in \eqref{eq:diagblock_FPchol}), so
\begin{equation} \label{eq:zeta_trunc_offdiag_mx}
\begin{aligned}
    \zeta^{\mathrm{trunc}}_{st} &\le (2\epsilon_\psi + \epsilon_\psi^2)\sqrt{\Bigl(\sum_{k=1}^{(J-1)n_b}\hat l^2_{sk} + \nu \mu_\Psi^2 d_{ss}^2\Bigr)\Bigl(\sum_{k=1}^{(J-1)n_b}\hat l^2_{tk} + \nu \mu_\Psi^2 d_{tt}^2\Bigr)} \\
    &\le (2\epsilon_\psi + \epsilon_\psi^2)\sqrt{(|a_{ss}| + |e_{ss}| + \nu \mu_\Psi^2 d_{ss}^2)(|a_{tt}| + |e_{tt}| + \nu \mu_\Psi^2 d_{tt}^2)}. 
\end{aligned}
\end{equation}
Similarly, for $\zeta^{\mathrm{L\_accum}}_{st}$, we can write the bound
\begin{equation}\label{eq:zeta_accum_offdiagblk_mx}
    \zeta^{\mathrm{L\_accum}}_{st} \le \left(n_b + 2 + \frac{n}{n_b}\right)\epsilon_h(1 + \epsilon_\psi)^2 \sqrt{\Bigl(\sum_{k=1}^{(J-1)n_b} \hat{l}_{sk}^2 + \nu \mu_\Psi^2 d_{ss}^2\Bigr)\Bigl(\sum_{k=1}^{(J-1)n_b} \hat{l}_{tk}^2 + \nu \mu_\Psi^2 d_{tt}^2\Bigr)}.
\end{equation}
Substituting \eqref{eq:zeta_accum_offdiagblk_mx}, \eqref{eq:zeta_trunc_offdiag_mx}, \eqref{eq:cholFP-TRSM-err} in \eqref{eq:est-offdiag-block}, and using $|a_{st}| = |\sum_{k=1}^{t}\hat{l}_{sk}\hat{l}_{tk} - e_{st}| \le \sqrt{(|a_{ss}| + |e_{ss}|)(|a_{tt}| + |e_{tt}|)} + |e_{st}|$, we get the same bounds as \eqref{eq:est-fp-diagblock} and \eqref{eq:est-fp-diagblock-final}.
Finally, to get \eqref{eq:chol-fp-mx} we need to bound $d_{rr}$ for $1 \le r \le n$. Since $d_{rr}$ is a scale factor in row $r$, $d_{rr} \le \frac{(1 + \epsilon_{\mathrm{scal}})}{\beta} \max_{k=1}^r |\hat{l}_{rk}| \le \frac{(1 + \epsilon_{\mathrm{scal}})}{\beta} \sqrt{|a_{rr}| + |e_{rr}|}$. Substituting this in Equation \eqref{eq:ett-fp}, we get 
\begin{equation} \label{eq:ett-mx}
\begin{aligned}
    (1 - \gfp)|e_{tt}| &\le \gfp\Bigl(|a_{tt}| + \frac{\nu\mu_\Psi^2(1 + \epsilon_{\mathrm{scal}})^2}{\beta^2}(|a_{tt}| + |e_{tt}|)\Bigr) \\
    \implies |e_{tt}| &\le \frac{\gfp (1 + \lambda_\Psi)}{1 - \gfp (1 + \lambda_\Psi)}|a_{tt}|.
\end{aligned}
\end{equation}
For $|e_{st}|$, using the bound on $d_{rr}$ in Equation \eqref{eq:est-fp-diagblock}, we obtain
\begin{equation}
    |e_{st}| \le \frac{\gfp - C}{1 - C} \sqrt{(1 + \lambda_\Psi)(|a_{ss}| + |e_{ss}|)} \sqrt{(1 + \lambda_\Psi)(|a_{tt}| + |e_{tt}|)} \le \gfp \sqrt{(1 + \lambda_\Psi)(|a_{ss}| + |e_{ss}|)} \sqrt{(1 + \lambda_\Psi)(|a_{tt}| + |e_{tt}|)}.
\end{equation}
Substituting \eqref{eq:ett-mx} above gives us the desired inequality.
\end{proof}

\begin{theorem} \label{thm:chol-fp-ns}
    Suppose we do a Cholesky with the same assumptions as Theorem~\ref{thm:chol-fp}, but we don't scale (for example, FP16 inputs accumulating into FP32, with FP32 output), and if $\gfp < 1$ we have $A + E = \hat L\hat L^T$ with
\begin{equation} \label{eq:chol-fp-ns}
  |e_{st}| \;\le\; \frac{\gfp}{1 - \gfp}\sqrt{a_{ss} + \nu \mu_\Psi^2} \sqrt{a_{tt} + \nu \mu_\Psi^2}.
\end{equation}
\end{theorem}
\begin{proof}
    As discussed in Section~\ref{sec:mx}, if we do not microscale then $d_{tt} = 1$. Plugging this into Equation \eqref{eq:chol-fp} gives us the desired inequality. 
\end{proof}

\subsection{Factorization Failures}\label{sec:failures}

It is shown in \cite{demmel1989cholesky} and \cite[Theorem 10.7]{Higham_main_book}  that we can have factorization failures when our machine epsilon gets too large. If Algorithm~\ref{A_chol} tries to give us \jwdnote{tries to give us} $\hat{L}\hat{L}^T = A + E$ such that $A + E$ ends up being indefinite, our factorization will fail. The failure/success guarantees from \cite[Theorem 10.7]{Higham_main_book} can be adapted to our error bounds in Section~\ref{sec:mp-chol}. Since all of our error bounds take the form $\frac{\gamma}{1 - \gamma} \sqrt{a_{ii}a_{jj}}$, we can state the following theorem:
\begin{theorem}
    Suppose we have a Cholesky factorization algorithm which on a successful run with input $A_0$, returns $\hat{L}$ \jwdnote{attempt to get} such that $\hat{L}\hat{L}^T = A_0 + E$ with $|e_{ij}| \le \frac{\gamma}{1 - \gamma}\sqrt{a_{ii}a_{jj}}$. Here $\gamma$ can be $\gamma_{\text{emu}}, \gfp$ or $\gfp (1 + \lambda_\Psi)$ from Theorems~\ref{thm:chol-emu}, \ref{thm:chol-fp} or~\ref{thm:chol-fp-ns}.   
    \jwdnote{Where is $\gamma$ defined?} Let $H = D^{-1}AD^{-1}$ with $D = \operatorname{diag}_i \sqrt{a_{ii}}$. Then the factorization is guaranteed to succeed if $\lambda_{\min}(H) >  \frac{n\gamma}{1 - \gamma}$, and the factorization is guaranteed to fail if  $\lambda_{\min}(H) < - \frac{n\gamma}{1 - \gamma}$.
\end{theorem}
\jwdnote{The proof in \cite[Theorem 10.7]{Higham_main_book} is more elaborate, just
cite it instead of reproducing it.}
\begin{proof}
See \cite[Theorem 10.7]{Higham_main_book}.
\end{proof}
We always have $\lambda_{\min} (H) \le 1$ since the sum of all eigenvalues of $H$ is $n$. 
\jwdnote{No, what about $H = [1, .5; .5, 1]$.
Can instead say $\lambda_{min}(H) =\min_{x \neq 0} x^THx/x^Tx \leq e_i^T H e_i = H_{ii}=1.$} \svknote{Sorry, meant to say the sum of all eigenvalues of H is n (trace of H), which would mean that the smallest eigenvalue is less than or equal to 1.}
Say we do a Cholesky with FP8 (E4M3) MX GEMMs, if we assume no underflow, and ignore higher precision rounding errors, by Theorem~\ref{thm:chol-fp}, we have $\gamma \approx 2\epsilon_l + \epsilon_l^2 = 0.265625$. So we cannot guarantee success for $n \ge 3$. Similarly, for FP16 inputs and FP32 accumulation, if there is no underflow, $\gamma \approx 0.00195$, so we cannot guarantee success for $n \ge 511$. \jwdnote{$n=512$?} To deal with this, we propose using a modified Cholesky factorization. An extensive overview and comparison of the various modified Cholesky factorizations is provided in \cite{McSweeney2020}. We consider two modified Cholesky methods which perturb the diagonal of $A$ during the factorization to ensure the diagonal stays positive. The first is a simple clamping algorithm,
which clamps each updated diagonal entry so that it never falls below a user-specified
$\delta^{\mathrm{clamp}} > 0$. The second is the algorithm of Gill, Murray, and Wright
\cite{GillMurrayWright1981}, which we use without pivoting. We will refer to their algorithm as GMW81. \jwdnote{Omit "(GMW11)"?} \svknote{I refer to their algorithm as GMW81 later. Perhaps I can refer to the algorithm as something else? Or I can say "we will refer to their algorithm as GMW81"?}
Both algorithms are equivalent to a Cholesky factorization of $A + F$ for some
diagonal $F \ge 0$, so Theorems~\ref{thm:chol-fp} and~\ref{thm:chol-fp-ns} still apply with $\sqrt{a_{ss}a_{tt}}$
replaced by $\sqrt{(a_{ss} + f^*)(a_{tt} + f^*)}$, where $f^* = \max_i f_{ii}$ (plus the underflow term of
Theorem~\ref{thm:chol-fp-ns} if underflow is expected to be significant). It therefore suffices to
bound $f^*$.

We will use a similar notation to Section~\ref{sec:mp-chol}, for $j \le i$ let $a^u_{ij} = a_{ij} - \sum_{k=1}^{j-1} l_{ik} l_{jk}$. Note that $a_{ii}^u$ may be negative. Both diagonal perturbation schemes, along with their respective bounds on $f^*$ are described in the following.
\subsubsection{Diagonal Clamping}
In this scheme, we set $l_{jj}^2 = \max( a^u_{jj}, \delta^{\mathrm{clamp}})$, so $a_{jj}^u + f_{jj} = \max(a_{jj}^u, \delta^{\mathrm{clamp}}) \implies f_{jj} = \max(0,\delta^{\mathrm{clamp}} - a^u_{jj})$. Also, $l_{jj} \ge \sqrt{\delta^{\mathrm{clamp}}}$, and for $i > j$,
$|l_{ij}| \le \frac{1}{\sqrt{\delta^{\mathrm{clamp}}}}\Bigl(|a_{ij}| + \sum_{k=1}^{j-1}|l_{ik}||l_{jk}|\Bigr)$
. If $|l_{ij}| \le \lambda$ for all $i > j$, then
\[
f_{ii} = \max\Bigl(0,\ \delta^{\mathrm{clamp}} - a_{ii} + \sum_{k=1}^{i-1} l_{ik}^2\Bigr) \implies f^*
\le \delta^{\mathrm{clamp}} + \max_t |a_{tt}| + (n-1)\lambda^2 .
\]
We cannot know $\lambda$ ahead of time; all we can say is that $\lambda$ is at least $O(1/\sqrt{\delta^{\mathrm{clamp}}})$, so a very small $\delta^{\mathrm{clamp}}$, or a very large $\delta^{\mathrm{clamp}}$ will both lead to large errors. Since the perturbation doesn't use data from anywhere else in the matrix, this perturbation scheme can be incorporated easily in the blocking used in Algorithm~\ref{A_chol}, which uses block diagonal Cholesky + TRSM + GEMM/SYRK.

\subsubsection{GMW81}
The perturbation rule here is
$l_{jj}^2 = \max\!\left(\delta^{\mathrm{GMW}},\ |a^u_{jj}|,\ \frac{\bigl(\theta^{\mathrm{GMW}}_j\bigr)^2}{\bigl(\beta^{\mathrm{GMW}}\bigr)^2}\right)$, with
$\theta^{\mathrm{GMW}}_j = \max_{i > j} |a^u_{ij}|$,
where $\beta^{\mathrm{GMW}}$ and $\delta^{\mathrm{GMW}}$ are user inputs. Since $|l_{jj}| \ge \theta^{\mathrm{GMW}}_j/\beta^{\mathrm{GMW}}$, for every $i > j$, we have
$|l_{ij}| = |a^u_{ij}|/|l_{jj}| \le \beta^{\mathrm{GMW}} |a^u_{ij}|/\theta^{\mathrm{GMW}}_j \le \beta^{\mathrm{GMW}}$, so the GMW algorithm gives us an a priori upper bound on the elements of $L$. Now, we want a bound on $|a_{jj}^u|$.
$|a^u_{jj}| \le |a_{jj}| + \sum_{k=1}^{j-1} l_{jk}^2 \le \max_t |a_{tt}| + (j-1)\bigl(\beta^{\mathrm{GMW}}\bigr)^2$. Next, for $\theta^{\mathrm{GMW}}_j$, we need a bound on $|a^u_{ij}|$ for $i > j$:
$|a^u_{ij}| \le |a_{ij}| + \sum_{k=1}^{j-1}|l_{ik}||l_{jk}| \le \max_{s \ne t}|a_{st}| + (j-1)\bigl(\beta^{\mathrm{GMW}}\bigr)^2$. Since $\theta^{\mathrm{GMW}}_j$ is the maximum of $|a^u_{ij}|$ over $i > j$,
$\theta^{\mathrm{GMW}}_j \le \max_{s \ne t} |a_{st}| + (j-1)\bigl(\beta^{\mathrm{GMW}}\bigr)^2$. Finally, since $f_{jj} = l_{jj}^2 - a^u_{jj}$ ,
\[
\begin{aligned}
f_{jj} &\le \delta^{\mathrm{GMW}} + \frac{\bigl(\theta^{\mathrm{GMW}}_j\bigr)^2}{\bigl(\beta^{\mathrm{GMW}}\bigr)^2} + 2|a^u_{jj}| \\
\implies f^* &\le \delta^{\mathrm{GMW}} + 2\Bigl(\max_t |a_{tt}| + (n-1)\bigl(\beta^{\mathrm{GMW}}\bigr)^2\Bigr)
+ \left(\frac{\max_{s \ne t} |a_{st}|}{\beta^{\mathrm{GMW}}} + (n-1)\beta^{\mathrm{GMW}}\right)^{2}.
\end{aligned}
\]
Gill, Murray, and Wright \cite{GillMurrayWright1981} prescribe
$\bigl(\beta^{\mathrm{GMW}}\bigr)^2 = \max\bigl(\max_t |a_{tt}|,\ \max_{s \ne t}|a_{st}|/\sqrt{n^2-1},\ \epsilon\bigr)$,
where $\epsilon$ is the machine precision.
 It is important to note that the perturbation rule uses all subdiagonal elements in column $j$ to compute $\theta^{\mathrm{GMW}}_j$. As a consequence of this, we cannot block this algorithm the same way as Algorithm~\ref{A_chol}. Instead, the TRSM and diagonal block solve need to be merged into one "panel solve" that is applied to a block column of width $n_b$.

\section{Adaptive precision Algorithms} \label{sec:adaptive_alg}

In this section, we introduce a framework for adaptive precision Cholesky algorithms which dynamically choose the number of slices (or the precision) of each GEMM or TRSM operation. At each iteration, the algorithms choose the smallest number of slices (or the smallest precision) that ensure that the computed $n \times n$ $\hat{L}$ factor satisfies $|a_{ij} - \sum_{k=1}^{\min (i,j)} \hat{l}_{ik}\hat{l}_{jk}| \le \epsilon' \sqrt{a_{ii}a_{jj}}$ for $1 \le i,j \le n$, where $A$ is a positive definite $n\times n$ input matrix, and $\epsilon'$ is a user-defined threshold. When $\epsilon'$ is chosen to be equal to $n\epsilon_{64}$, the bound $\epsilon'\sqrt{a_{ii}a_{jj}}$ coincides with the classical component-wise error bound of the Cholesky factorization \cite{demmel1989cholesky, Higham_main_book}.
 
\begin{algorithm}
\caption{Blocked Cholesky Factorization with Adaptive Slice Ozaki 1 GEMM}
\label{alg:adaptive}
\begin{algorithmic}[1]
\Require Matrix $A \in \mathbb{R}^{n \times n}$ (SPD) in some high precision $h$, block size $n_b$, tolerance $\epsilon'$, max slice count $\Psi'$
\jwdnote{Say what types of all variables are.} \svknote{Thanks, just added all types/ slice counts. Let me know what you think}
\Ensure Lower Cholesky factor $L$ (in FP64) such that $A = LL^T$ (in-place)
\State $\mathrm{diag}_i \gets a_{ii}$ for $i = 1,\ldots,n$ \Comment{$\mathrm{diag} \in \mathbb{R}^n$: original diagonal of $A$ in $h$}
\State $\mathrm{diag}^u_i \gets a_{ii}$ for $i = 1, \ldots , n$
\Comment{$\mathrm{diag}^u \in \mathbb{R}^n$, in $h$}
\State $m_i \gets 0$ for $i = 1,\ldots,n$ \Comment{$m \in \mathbb{R}^n$ in $h$}
\State $s_i \gets 0$ for $i = 1,\ldots,n$
\Comment{$s \in \mathbb{R}^n$ in $h$}
\State $g_i \gets 0$ for $i = 1,\ldots,n$ \Comment{$g \in \mathbb{R}^n$ in $h$}
\State $c \gets 0$ \Comment{in $h$}
\State $\psi \gets \Psi'$ 
\For{$I = 1, \ldots, n/n_b$}
    \State $o \gets (I{-}1)n_b$ 
    \State $w \gets \min(n_b,\, n - o)$ 
    \State $A_{II} \gets A[o{+}1{:}o{+}w,\; o{+}1{:}o{+}w]$ \Comment{diagonal block}
    \State $A_{II} \gets \textsc{Potrf}(A_{II})$ \Comment{$A_{II} \leftarrow L_{II} L_{II}^T$ in $h$}
    \If{$o + w < n$}
        \State $A_{I{+}1,I} \gets A[o{+}w{+}1{:}n,\; o{+}1{:}o{+}w]$ \Comment{off-diagonal panel}
        \State $A_{I{+}1,I} \gets A_{I{+}1,I} \cdot A_{II}^{-T}$  \Comment{TRSM is emulated using $\psi$ slices}
        \For{$i = o{+}w{+}1, \ldots, n$} \Comment{one pass over row $i$ of the panel}
            \State $g_i \gets \max_{o<k\le o+w} |l_{ik}|$
            \State $s_i \gets s_i + \sum_{k=o+1}^{o+w} l_{ik}^2$
            \State $\mathrm{diag}^u_i \gets \mathrm{diag}^u_i - \sum_{k=o+1}^{o+w} l_{ik}^2$
        \EndFor
        \For{$p_1 = 1, \ldots , \Psi'$}\label{line:emu-select-start}
        \For{$p_2 = 1, \ldots , p_1$}
            \If{$\Theta (p_2, n_b)g_i^2 + (\max(c, \Gamma (p_2)) + I + 1)\epsilon_{64}(s_i + a_{ii}) + \frac{\gamma_{\mathrm{emu}}(p_1, \Psi', n - In_b, n_b)}{1 - \gamma_{\mathrm{emu}}(p_1, \Psi', n - In_b, n_b)}a_{ii}^u \le \epsilon' a_{ii} - m_i$ for $o+w < i \le n$ and $ \gamma_{\mathrm{emu}}(p_1, \Psi', n - I n_b, n_b ) < 1$}\label{line:emu-select-if}
            \State $\pi_g(I) \gets p_2, \psi \gets p_1$
            \State goto line~\ref{line:emu-after-select}
            \EndIf
        \EndFor
        \EndFor\label{line:emu-select-end}
        \State\label{line:emu-after-select} $c \gets \max(c, \Gamma(\pi_g(I)))$
        \For{$i = o{+}w{+}1, \ldots, n$}
            \State $m_i \gets m_i + \Theta(\pi_g(I), n_b)\,g_i^2$
        \EndFor
        \State $A_{I{+}1,I{+}1} \gets A[o{+}w{+}1{:}n,\; o{+}w{+}1{:}n]$ \Comment{trailing submatrix}
        \State $A_{I{+}1,I{+}1} \gets A_{I{+}1,I{+}1} - A_{I{+}1,I}\, A_{I{+}1,I}^T$ 
        \Comment{GEMM/SYRK is emulated using $\pi_g(I)$ slices}
    \EndIf
\EndFor
\end{algorithmic}
\end{algorithm}

\subsection{Adaptive precision with Ozaki 1}
For simplicity, assume $n_b|n$, and to be concrete, we focus on the description of the algorithm using emulation based on the Ozaki 1 scheme. The mixed-precision floating-point case, and the Ozaki 2 case are very similar, and we state the algorithm without proof.  Algorithm \ref{A_chol} serves as the basis for our adaptive precision algorithm. The key modification is to select the number of slices used in the TRSM and GEMM operations (lines~\ref{line:chol-trsm} and~\ref{line:chol-gemm} of Algorithm~\ref{A_chol}, respectively) at every iteration, instead of using a pre-determined list of slices throughout the factorization.  The following theorem defines the conditions under which the number of slices can be changed while ensuring that the output satisfies $|a_{ij} - \sum_{k=1}^{\min (i,j)} \hat{l}_{ik} \hat{l}_{jk}| \le \epsilon' \sqrt{a_{ii}a_{jj}}$.  Given an input matrix $A$, and considering the iteration $I$ of Algorithm \ref{A_chol}, the input matrix $A$ and the calculated factor $\hat{L}$ are partitioned into block matrices of dimension $3 \times 3$ such that $A - \hat{L}\hat{L}^T$ is written as 
\begin{equation}
\begin{pmatrix}
  A_{11}   & A_{21}^{T} & A_{31}^{T} \\
  A_{21}   & A_{22}     & A_{32}^{T} \\
  A_{31}   & A_{32}     & A_{33}
\end{pmatrix}
-
\begin{pmatrix}
  \hat{L}_{11} & 0            & 0 \\
  \hat{L}_{21} & \hat{L}_{22} & 0 \\
  \hat{L}_{31} & \hat{L}_{32} & \hat{L}_{33}
\end{pmatrix}
\begin{pmatrix}
  \hat{L}_{11} & 0            & 0 \\
  \hat{L}_{21} & \hat{L}_{22} & 0 \\
  \hat{L}_{31} & \hat{L}_{32} & \hat{L}_{33}
\end{pmatrix}^{T},
\end{equation}
where $A_{11}$ is $(I-1)n_b \times (I-1)n_b$, $A_{22}$ is $n_b \times n_b$, and $A_{33}$ is $(n - In_b) \times (n - In_b)$. The blocks of $\hat{L}$ are defined in a similar way. The only submatrix of $\hat{L}$ whose computation is not complete is $\hat{L}_{33}$. Once the Schur complement of $A_{33}$ is computed as $A_{33} - \hat{L}_{31}\hat{L}_{31}^T - \hat{L}_{32}\hat{L}_{32}^T$, according to Algorithm \ref{A_chol}, a mixed-precision Cholesky factorization of that Schur complement is performed. The backward error on $A_{33}$ is then $|A_{33} - \hat{L}_{31}\hat{L}_{31}^T - \hat{L}_{32}\hat{L}_{32}^T - \hat{L}_{33}\hat{L}_{33}^T|$. Denote $A_{33}^{(\mathrm{curr})}$ as the computed result of $ A_{33} - \hat{L}_{31}\hat{L}_{31}^T$ and $A_{33}^{(\mathrm{next})}$ as the computed result of $A_{33}^{(\mathrm{curr})} - \hat{L}_{32}\hat{L}_{32}^T$. Then, we can write the following telescoping error bound for $A_{33}$,
\begin{equation} \label{eq:telescope bound}
\begin{aligned}
    |A_{33} - \hat{L}_{31}\hat{L}_{31}^T - \hat{L}_{32}\hat{L}_{32}^T - \hat{L}_{33}\hat{L}_{33}^T| &\le |(A_{33} - \hat{L}_{31}\hat{L}_{31}^T) - A_{33}^{(\mathrm{curr})}| \\ &\quad + |(A_{33}^{(\mathrm{curr})} - \hat{L}_{32}\hat{L}_{32}^T) - A_{33}^{(\mathrm{next})}| + |(A_{33}^{(\mathrm{next})} - \hat{L}_{33}\hat{L}_{33}^T)| \\
    &\le E^{(1)} + E^{(2)} + E^{(3)},
\end{aligned}
\end{equation}
 where $E^{(1)}, E^{(2)}, E^{(3)}$ are slightly loosened worst case error bounds on $|(A_{33} - \hat{L}_{31}\hat{L}_{31}^T) - A_{33}^{(\mathrm{curr})}|$, $|(A_{33}^{(\mathrm{curr})} - \hat{L}_{32}\hat{L}_{32}^T) - A_{33}^{(\mathrm{next})}|$ and $|(A_{33}^{(\mathrm{next})} - \hat{L}_{33}\hat{L}_{33}^T)|$ obtained from \eqref{eqn:GEMMerrbnd}, and $\eqref{eq:chol-emu}$, respectively. In the proof of the following Theorem \ref{th:adaptPrecEisSPD} we explain why we need to  relax the bound and how we achieve this.
 
\begin{theorem}
\label{th:adaptPrecEisSPD}
    Consider Algorithm \ref{A_chol} that uses a $\pi_g(I)$-slice Ozaki 1 at iteration $I$ for the Schur update and a $\pi_t(I)$-slice Ozaki 1 for the TRSM. Then in Equation \eqref{eq:telescope bound},  $E' = E^{(1)} +  E^{(2)} + E^{(3)}$ is symmetric positive semidefinite. 
\end{theorem} 
\begin{proof}
For $1 \le K \le I$, define $\hat{L}_{3}^{(K)} = \hat{L}[I\,n_b+1{:}n,\; (K-1)\,n_b+1{:}K\,n_b]$. Then $\hat{L}_{32} = \hat{L}_3^{(I)} = \hat{L}[I\,n_b+1{:}n,\; (I-1)\,n_b+1{:}I\,n_b]$ and
\begin{equation*}
\hat{L}_{31}\hat{L}_{31}^T = \sum_{K=1}^{I-1}\hat{L}_{3}^{(K)}\hat{L}_{3}^{(K)T},
\end{equation*}
where the $K$-th GEMM is done with $\pi_g(K)$ slices. Our goal now is to apply \eqref{eqn:GEMMerrbnd}. For this, we define $g^{(K)}$ as the length $n - In_b$ vector \jwdnotenew{Fix grammar; remove first "vector"?} containing the ufp max norms of the rows of $\hat{L}_{3}^{(K)}$, or equivalently, the ufp max norms of the columns of $\hat{L}_{3}^{(K)T}$. Then
\begin{equation} \label{eq:err_bnd_e1+e2}
    E^{(1)} + E^{(2)} = \sum_{K=1}^I \Theta(\pi_g(K), n_b)g^{(K)}g^{(K)T} + \bigl(\max_{K \le I}\Gamma(\pi_g(K)) + I + 1\bigr)\epsilon_{64}\Bigl(\sum_{K=1}^I |\hat{L}^{(K)}_3||\hat{L}^{(K)}_3|^T + dd^T\Bigr),
\end{equation}
where $d$ is a length $n - In_b$ vector with $d_r = \sqrt{a_{r + In_b, r + In_b}}$. We use $dd^T$ instead of $|A_{33}|$ since although $A_{33}$ is positive semidefinite, $|A_{33}|$ need not be. But, since $A_{33}$ is a principal submatrix of $A$, $|A_{33}| \le dd^T$ by Cauchy-Schwarz. In equation \eqref{eq:err_bnd_e1+e2},  $g^{(K)}g^{(K)T}$,  $|\hat{L}^{(K)}_3||\hat{L}^{(K)}_3|^T$, and $dd^T$, are all positive semi-definite, and since \eqref{eq:err_bnd_e1+e2} only has positive coefficients, the sum of all matrices is also positive semidefinite. Next, for $E^{(3)}$, we apply Theorem~\ref{thm:chol-emu} to $A_{33}^{(\mathrm{next})}$ (under the assumption that Cholesky is successful). Denote its diagonal entries by $a_{rr}^u$ for $1 \le r \le n - In_b$. Then 
$e^{(3)}_{ij} = \frac{\gamma_{\mathrm{emu}}(\psi, \Psi, n - In_b, n_b)}{1 - \gamma_{\mathrm{emu}}(\psi, \Psi, n - In_b, n_b)}\sqrt{a^u_{ii}a^u_{jj}}$, where $\psi$ and $\Psi$ are $\min_J \min (\pi_g (J), \pi_t (J))$ and $\max_J \max (\pi_g (J), \pi_t (J))$ respectively for $I < J \le n/n_b$. This is also positive semidefinite, since it is another outer product of a vector with itself. 
\end{proof}
Since $E' = E^{(1)} + E^{(2)} + E^{(3)}$ is positive semidefinite, using Cauchy-Schwarz, we obtain $e'_{ij}\le \sqrt{e'_{ii}e'_{jj}}$.  Given a tolerance $\epsilon'$, we want our algorithm
\jwdnotenew{our algorithm} to ensure that the computed $\hat{L}$ satisfies $|a_{ij} - \sum_{k=1}^{\min(i,j)}l_{ik}l_{jk}| \le \epsilon' \sqrt{a_{ii}a_{jj}}$. This requires $\sqrt{e'_{ii}e'_{jj}} \le \epsilon'\sqrt{a_{ii}a_{jj}}$ at every iteration.  Equivalently, it is much cheaper to ensure that $e'_{ii} \le \epsilon' a_{ii}$ for every $1 \le i \le n$. Using \eqref{eq:err_bnd_e1+e2} and $e^{(3)}_{ii} = \frac{\gamma_{\mathrm{emu}}(\psi, \Psi, n - In_b, n_b)}{1 - \gamma_{\mathrm{emu}}(\psi, \Psi, n - In_b, n_b)}a^u_{ii}$, we can write 
\begin{equation} \label{formula:compute_eii}
\begin{aligned}
    e'_{ii} &= \sum_{K=1}^{I-1}\Theta (\pi_g (K), n_b)\bigl(g^{(K)}_i\bigr)^2 + \Theta (\pi_g (I), n_b)\bigl(g^{(I)}_i\bigr)^2 \\ &\qquad + \bigl(\max_{K \le I}\Gamma(\pi_g(K)) + I + 1\bigr)\epsilon_{64}\Bigl(\sum_{k=1}^{In_b}\hat{l}_{ik}^2 + a_{ii}\Bigr) + \frac{\gamma_{\mathrm{emu}}(\psi, \Psi, n - In_b, n_b)}{1 - \gamma_{\mathrm{emu}}(\psi, \Psi, n - In_b, n_b)} a_{ii}^u.
\end{aligned}
\end{equation}

Then, bounding $e'_{ii}$ by $\epsilon' a_{ii}$ is equivalent to checking the following inequality 
\begin{equation} \label{formula:e_ii-ineq}
\begin{aligned}
     &\Theta (\pi_g (I), n_b)\bigl(g^{(I)}_i\bigr)^2 + \bigl(\max_{K \le I}\Gamma(\pi_g(K)) + I + 1\bigr)\epsilon_{64}\Bigl(\sum_{k=1}^{In_b}\hat{l}_{ik}^2 + a_{ii}\Bigr) + \frac{\gamma_{\mathrm{emu}}(\psi, \Psi, n - In_b, n_b)}{1 - \gamma_{\mathrm{emu}}(\psi, \Psi, n - In_b, n_b)} a_{ii}^u \\ &\qquad \le \epsilon'a_{ii} - \sum_{K=1}^{I-1}\Theta (\pi_g (K), n_b)\bigl(g^{(K)}_i\bigr)^2.
\end{aligned}
\end{equation}

The  blocked Cholesky factorization with adaptive precision is presented in Algorithm~\ref{alg:adaptive}.  This algorithm incorporates an efficient adaptive slice selection and is based on the observation that every quantity required to evaluate $e'_{ii}$ in \eqref{formula:compute_eii} can be computed with $O(n\cdot n_b)$ work.
\jwdnotenew{Clarify: if this is the additional
cost per diagonal, state what the overall
additional cost is, and that it is a lower order
term, and so worth doing. Perhaps at the end of
the paragraph?} \svknotenew{Thanks, added to the end of the paragraph}
First, all elements $a_{ii}$ can be read from $A$ and stored in a separate array of length $n$ before starting the factorization. Each $\bigl(g_i^{(K)}\bigr)^2$ for $1 \le K \le \frac{n}{n_b}$ can be evaluated by computing $g_i^{(K)} = \mathrm{ufp}\bigl(\max_{k=(K-1)n_b+1}^{Kn_b} |\hat{l}_{ik}|\bigr)$ before each GEMM. Instead of the ufp, we use the plain max $\max_{(K-1)n_b < k \le Kn_b} |\hat{l}_{ik}|$ in Algorithm~\ref{alg:adaptive} (and in Algorithm~\ref{alg:adaptive-oz2}); since $\mathrm{ufp}(x) \le x$, this only makes the bound larger.  An array of length $n$ can be used to store and update at each iteration $m_i = \sum_{K=1}^{I-1} \Theta (\pi_g (K), n_b)(g_i^{(K)})^2$. We also need to compute $\sum_{k=(K-1)n_b+1}^{Kn_b}\hat{l}_{ik}^2$ in every iteration, and this computation can be done in the same pass used to find the maximum. The last quantity we need to compute is $a_{ii}^u$. This can also be approximated in the same pass, since we need to approximate $a_{ii} - \sum_{k=1}^{In_b} \hat{l}^2_{ik}$. To get an accurate approximation of this expression, we would essentially need an Ozaki style dot product function, but we find in our experiments that using regular double precision for this operation still gives us error guarantees, hence our algorithm assumes that $a_{ii} - \sum_{k=1}^{In_b} \hat{l}^2_{ik}$ is computed in double precision.  Since there are $n/n_b$ iterations, tracking  \eqref{formula:compute_eii} costs $O(n^2)$ work throughout the course of the $O(n^3)$ Cholesky algorithm. So, computing \eqref{formula:compute_eii} adds minimal overhead to Algorithm \ref{alg:adaptive}.  
\jwdnotenew{Perhaps again point out that
since only does $O(n^2)$ flops, it is cheap.}

Thus Algorithm~\ref{alg:adaptive} uses five additional arrays of length $n$ to monitor the roundoff error. Array $\mathrm{diag}$ with $\mathrm{diag}_i = a_{ii}$, $\mathrm{diag}^u$, which is updated at every iteration and stores $\mathrm{diag}^u_i = a^u_{ii} = a_{ii} - \sum_k l_{ik}^2$, array $g$ which stores the row-wise ufp max norms of the panel with which we are just about to perform SYRK/GEMM, array $m$ which stores the error from truncation as $m_i = \sum \Theta (\pi_g (K), n_b) (g_i^{(K)})^2$, and finally $s$ which stores the error from accumulation as $s_i = \bigl(\max_{K \le I}\Gamma(\pi_g(K)) + I + 1\bigr)\epsilon_{64}\Bigl(\sum_{k=1}^{In_b}\hat{l}_{ik}^2 + a_{ii}\Bigr)$.  Since computing $e'_{ii}$ is cheap, the goal of our adaptive precision algorithm is to solve $\eqref{formula:e_ii-ineq}$ for $\pi_g(I)$ and $\psi$ in every iteration. We do this in lines~\ref{line:emu-select-start}--\ref{line:emu-select-end} of Algorithm \ref{alg:adaptive}.
\jwdnotenew{Refer to this being done in lines 21-28 of Alg 6.}
When we want a double precision Cholesky factorization, we set $\epsilon' = n\epsilon_{64}$ and set $\Psi' = 8$, but $\Psi' = 7$ should also work. Also note that we default to 8 slices, not FP64, if we can't find any solutions to \eqref{formula:e_ii-ineq} in line~\ref{line:emu-select-if}. Different possible solutions exist for $\pi_g(I), \psi$. We consider two possibilities, $\pi_g(I) =  \psi$ and $\pi_g(I) \le \psi$, Algorithm \ref{alg:adaptive} displays the latter. We enforce $\pi_g(I) \le \psi$ in the nested loops in lines~\ref{line:emu-select-start}--\ref{line:emu-select-end} by using the loop index of the $\psi$ (outer) loop as the upper bound of the $\pi_g(I)$ (inner) loop. Solutions with $\pi_g(I) \ge \psi$ are not interesting since SYRK/GEMM is the computationally expensive part of Cholesky, and TRSM is more numerically sensitive than GEMM. So ideally, TRSM should use at least as many slices as GEMM.
\jwdnotenew{Is this because of way loops are
nested in lines 21-22, or something else?} \svknotenew{Yes, the loop bounds enforce $\pi_g(I) \le \psi$. Does this answer your question? Or was your question about why we don't consider $\pi_g(I) \ge \psi$?}
Experimentally, we found that choosing a tighter $\Psi'$ (for example, 7, which also gives us
\jwdnotenewnew{change ``less than'' to ``more accurate than''} double precision accuracy for sufficiently large $n$) tends to make Algorithm \ref{alg:adaptive} slower \jwdnotenew{``performance'' could
mean speed or accuracy, clarify.}
than choosing $\Psi' = 8$ since there is less flexibility for GEMMs to utilize lower precisions. Additionally, in our experiments, we find that the choice $\pi_g(I) = \psi$ gives us the same slice choice at every iteration both with and without equilibration by a diagonal matrix of powers of 2, but $\pi_g(I) \le \psi$ does not have that property (unless all the diagonal values are equal). Finally, although Algorithm \ref{alg:adaptive} provides the bound $\epsilon' \sqrt{a_{ii}a_{jj}}$, we can also change the algorithm to return a factorization that satisfies $\epsilon'\max_i a_{ii}$ or $\epsilon' ||A||$  by modifying \eqref{formula:e_ii-ineq} appropriately.
\subsection{Adaptive precision with mixed-precision floating-point}
The blocked floating-point Cholesky factorization with adaptive precision GEMM is presented in Algorithm~\ref{alg:adaptive-FP}. 
It relies on an adaptation of \eqref{formula:e_ii-ineq} to the floating-point case, which requires using Theorem~\ref{thm:chol-fp} in place of Theorem \ref{thm:chol-emu}, and Theorem~\ref{thm:MXGEMM} in place of \eqref{eqn:GEMMerrbnd}. 
The floating-point version of \eqref{formula:e_ii-ineq} then becomes
\begin{equation} \label{formula:e_ii-ineq-fp}
\begin{aligned}
    &\tau(\epsilon_{\pi_g(I)})\Bigl(\sum_{k=(I-1)n_b+1}^{In_b}\hat{l}_{ik}^2 + \nu\mu_{\pi_g(I)}^2\bigl(d_i^{(I)}\bigr)^2\Bigr) + \left(n_b + 2 + \frac{n}{n_b}\right)\eh\,a_{ii} + \frac{\gamma_{\mathrm{fp}}(\epsilon_\psi, n - In_b, n_b)(1 + \lambda_\Psi)}{1 - \gamma_{\mathrm{fp}}(\epsilon_\psi, n - In_b, n_b)(1 + \lambda_\Psi)}a_{ii}^u \\
    &\qquad \le \epsilon' a_{ii} - \sum_{K=1}^{I-1}\tau(\epsilon_{\pi_g(K)})\Bigl(\sum_{k=(K-1)n_b+1}^{Kn_b}\hat{l}_{ik}^2 + \nu\mu_{\pi_g(K)}^2\bigl(d_i^{(K)}\bigr)^2\Bigr),
\end{aligned}
\end{equation}
where
\[
    \tau(\epsilon_p) = \begin{cases}
        2\epsilon_p + \epsilon_p^2 + \left(n_b + 2 + \frac{n}{n_b}\right)\eh(1 + \epsilon_p)^2 \approx 2\epsilon_p + \epsilon_p^2, & \text{if } p \ne h, \\
        \left(n_b + 2 + \frac{n}{n_b}\right)\eh, & \text{if } p = h,
    \end{cases}
\]
since for $p = h$ the inputs of the GEMM are not rounded, and only the accumulation error and the rounding from applying the scale factors remain. Here $d_i^{(K)}$ is the scale factor of row $i$ in panel $K$. With microscaling, the scale factor of row $i$ in panel $K$ is at most $(1 + \epsilon_{\mathrm{scal}})\max_{k}|\hat{l}_{ik}|/\beta$, so $\nu\mu_{\pi_g(K)}^2 \bigl(d_i^{(K)}\bigr)^2 \le \lambda_{\pi_g(K)}\sum_{k=(K-1)n_b+1}^{Kn_b}\hat{l}_{ik}^2$ with $\lambda_p = \frac{\nu\mu_p^2(1 + \epsilon_{\mathrm{scal}})^2}{\beta_p^2}$, and so $\tau(\epsilon_{\pi_g(K)})\bigl(\sum_{k=(K-1)n_b+1}^{Kn_b}\hat{l}_{ik}^2 + \nu\mu_{\pi_g(K)}^2\bigl(d_i^{(K)}\bigr)^2\bigr)$ reduces to $\tau(\epsilon_{\pi_g(K)})(1 + \lambda_{\pi_g(K)})\sum_{k=(K-1)n_b+1}^{Kn_b}\hat{l}_{ik}^2$. Without scaling, $d_i^{(K)} = 1$ and it reduces to $\tau(\epsilon_{\pi_g(K)})\bigl(\sum_{k=(K-1)n_b+1}^{Kn_b}\hat{l}_{ik}^2 + \nu\mu_{\pi_g(K)}^2\bigr)$. The $a^u_{ii}$ term is kept in the MX form of \eqref{eq:chol-fp-mx} since our experiments use MXFP8 (E4M3) as the lowest precision, but if the lowest precision was FP16, we would use Theorem \ref{thm:chol-fp-ns} instead.
Note that in \eqref{formula:e_ii-ineq-fp} one can also include the diagonal perturbation errors from Section~\ref{sec:failures}, by turning $a_{ii}$ into $a_{ii} + f^*$, and by accounting for the perturbation rule in $a^u_{ii}$. But there are two problems with this. If we use GMW81, $f^*$ is $O(n^2)$, and assuming $\beta^{\mathrm{GMW}} = 1$, and $a_{ii} = 1$ for $1 \le i \le n$, for $n= 40000$, $f^*$ is $1.6 \times 10^9$ which is clearly much larger than 1, making the $a_{ii} + f^*$ a very impractical bound. If we used $a_{ii} + f^*$ the algorithm would essentially default to using low precision for every operation even for modest $\epsilon'$. The other problem, specifically with GMW81, is the fact that $a^u_{ii}$ cannot be accurately estimated since it depends on $\theta^{\mathrm{GMW}}$ which depends on the values $a^u_{ji}$ for $j > i$ which are not known ahead of time. For these reasons, when using diagonal perturbation methods, we treat \eqref{formula:e_ii-ineq-fp} as a heuristic rather than a strict error guarantee. If we use \jwdnotenew{use} clamping instead of GMW81, then we don't have a concrete $f^*$. For the floating-point case, since TRSM is not emulated, we don't have the choice of solving for both GEMM precision and TRSM precision. So, rather than having a nested loop like lines~\ref{line:emu-select-start} to~\ref{line:emu-select-end} of Algorithm~\ref{alg:adaptive}, we have only one loop from line~\ref{line:fp-select-start} to~\ref{line:fp-select-end} in Algorithm~\ref{alg:adaptive-FP}.

The inputs of Algorithm~\ref{alg:adaptive-FP} are as follows. The matrix $A$ is stored in the high precision $h$, which is also the precision used for POTRF, TRSM, and the accumulation in each GEMM. The list $\mathcal{P}$ contains the candidate input formats for the GEMMs, ordered from least to most accurate, and each format $p \in \mathcal{P}$ is described by its machine epsilon $\epsilon_p$ and its underflow threshold $\mu_p$. If no format in $\mathcal{P}$ satisfies \eqref{formula:e_ii-ineq-fp}, the GEMM falls back to precision $h$. The underflow count $\nu$ is an assumed upper bound on the number of entries per row of $\hat{L}$ that underflow when rounded, as in Theorem~\ref{thm:chol-fp}, and $\beta$ is the microscaling target of Section~\ref{sec:mx}. The last input selects whether a diagonal perturbation is used. Diagonal clamping (Section~\ref{sec:failures}) only changes the pivots, so it can be done inside POTRF, but GMW81 needs the whole column below each pivot, so POTRF and TRSM are replaced by a single GMW81 factorization of the tall panel. For GMW81, the parameters $\beta^{\mathrm{GMW}}_{p}$ and $\delta^{\mathrm{GMW}}_p$ are given for each format $p$, and at iteration $I$ we use those of $\pi_g(I-1)$, the format of the previous GEMM. This is because the entries of panel $I$ were last updated by that GEMM, so their rounding errors are on the scale of $\epsilon_{\pi_g(I-1)}$. The first panel is not updated by any GEMM, so we set $\pi_g(0) = h$. The values used in our experiments are given in Section~\ref{sec:experiments}.
\begin{algorithm}
\caption{Blocked Floating-Point Cholesky Factorization with Adaptive Precision GEMM}
\label{alg:adaptive-FP}
\begin{algorithmic}[1]
\Require Matrix $A \in \mathbb{R}^{n \times n}$ (SPD) in some high precision $h$, block size $n_b$, tolerance $\epsilon'$, list of precisions $\mathcal{P}$ ordered from least to most accurate, with machine epsilons $\epsilon_p$ and underflow thresholds $\mu_p$ for $p \in \mathcal{P}$, (expected) underflow count $\nu$, microscaling target $\beta$, perturbation $\in \{\text{none}, \text{diagonal clamping}, \text{GMW81}\}$ (with $\beta^{\mathrm{GMW}}_{p}, \delta^{\mathrm{GMW}}_p$ for each $p \in \mathcal{P}$ for GMW81)
\Ensure Lower Cholesky factor $L$ (in $h$) such that $A = LL^T$ (in-place), or $A + F = LL^T$ with diagonal $F \ge 0$ for diagonal clamping or GMW81
\State $\mathrm{diag}_i \gets a_{ii}$ for $i = 1,\ldots,n$ \Comment{$\mathrm{diag} \in \mathbb{R}^n$: original diagonal of $A$ in $h$}
\State $\mathrm{diag}^u_i \gets a_{ii}$ for $i = 1, \ldots , n$
\Comment{$\mathrm{diag}^u \in \mathbb{R}^n$, in $h$}
\State $s_i \gets 0$ for $i = 1,\ldots,n$
\Comment{$s \in \mathbb{R}^n$ in $h$}
\State $t_i \gets 0$ for $i = 1,\ldots,n$ \Comment{$t \in \mathbb{R}^n$ in $h$}
\State $\Psi \gets \arg\max_{p \in \mathcal{P}} \mu_p$
\State $\pi_g(0) \gets h$
\For{$I = 1, \ldots, n/n_b$}
    \State $o \gets (I{-}1)n_b$ 
    \State $w \gets \min(n_b,\, n - o)$ 
    \State $A_{II} \gets A[o{+}1{:}o{+}w,\; o{+}1{:}o{+}w]$ \Comment{diagonal block}
    \State $A_{I{+}1,I} \gets A[o{+}w{+}1{:}n,\; o{+}1{:}o{+}w]$ \Comment{off-diagonal panel (empty if $o + w = n$)}
    \If{perturbation $=$ GMW81}
        \State $A[o{+}1{:}n,\; o{+}1{:}o{+}w] \gets \textsc{PanelGMW81}(A[o{+}1{:}n,\; o{+}1{:}o{+}w], \beta^{\mathrm{GMW}}_{\pi_g(I-1)}, \delta^{\mathrm{GMW}}_{\pi_g(I-1)})$ \Comment{in $h$}
    \Else
        \State $A_{II} \gets \textsc{Potrf}(A_{II})$ \Comment{in $h$, optionally with diagonal clamping}
        \State $A_{I{+}1,I} \gets A_{I{+}1,I} \cdot A_{II}^{-T}$  \Comment{TRSM in precision $h$}
    \EndIf
    \If{$o + w < n$}
        \For{$i = o{+}w{+}1, \ldots, n$} \Comment{one pass over row $i$ of the panel}
            \State $t_i \gets \sum_{k=o+1}^{o+w} l_{ik}^2$
            \State $\mathrm{diag}^u_i \gets \mathrm{diag}^u_i - t_i$
        \EndFor
        \For{$p \in \mathcal{P}$}\label{line:fp-select-start}
            \If{$\tau(\epsilon_p)(1 + \lambda_p)t_i + \left(n_b + 2 + \frac{n}{n_b}\right)\epsilon_h a_{ii} + \frac{\gamma_{\mathrm{fp}}(\epsilon_p, n - In_b, n_b)(1 + \lambda_\Psi)}{1 - \gamma_{\mathrm{fp}}(\epsilon_p, n - In_b, n_b)(1 + \lambda_\Psi)}a_{ii}^u \le \epsilon' a_{ii} - s_i$ for $o+w < i \le n$ and $ \gamma_{\mathrm{fp}}(\epsilon_p, n - In_b, n_b)(1 + \lambda_\Psi) < 1$}
            \State $\pi_g(I) \gets p, \psi \gets p$
            \State goto line~\ref{line:fp-after-select}
            \EndIf
        \EndFor\label{line:fp-select-end}
        \For{$i = o{+}w{+}1, \ldots, n$}\label{line:fp-after-select}
            \State $s_i \gets s_i + \tau(\epsilon_{\pi_g(I)})(1 + \lambda_{\pi_g(I)})\,t_i$
        \EndFor
        \State $A_{I{+}1,I{+}1} \gets A[o{+}w{+}1{:}n,\; o{+}w{+}1{:}n]$ \Comment{trailing submatrix}
        \State microscale $A_{I{+}1,I}$ and round it to format $\pi_g(I)$
        \State $A_{I{+}1,I{+}1} \gets A_{I{+}1,I{+}1} - A_{I{+}1,I}\, A_{I{+}1,I}^T$ \Comment{MX GEMM/SYRK, accumulating in $h$}
    \EndIf
\EndFor
\end{algorithmic}
\end{algorithm}

\section{Experiments}\label{sec:experiments}

In this section we present numerical experiments for the Cholesky decompositions from Algorithms~\ref{alg:adaptive} and~\ref{alg:adaptive-FP}. We consider two types of matrices.
\begin{itemize}
    \item Matrices of the form $A = DU\Lambda U^TD$, where $U$ is a random orthonormal matrix of the Haar distribution, $D$ is a random diagonal matrix with entries $2^{\mathrm{exp}}$ where $\mathrm{exp}$ is sampled uniformly from $[-R, R]$, and $\Lambda$ is a diagonal matrix with values that are defined by  the distributions in Table~\ref{tab:eigen_dists}. For Ozaki 1 experiments, $R = 511$, for PCG and Ozaki 2, $R = 0$.
    \item Matrices from SuiteSparse \cite{SuiteSparse2011}
    \jwdnotenew{Add citation} that we consider dense, and thus use our dense Cholesky algorithm to compute their factorization (to test accuracy, not performance).  
    \jwdnotenewnew{Note that this is
    good for testing precision, not performance.} \svknotenewnew{I agree, but I wonder if we should then change the plot in Fig 4? It currently shows speed. Should we change it to show backward error instead?}
\end{itemize} 
In the case of the Ozaki 1, double precision is emulated directly using the cuBLAS 13.2 implementation \jwdnotenew{Say how: Using cuBLAS xxx, or something like that.}. Our goal is to ensure that $|A - \hat{L}\hat{L}^T|$ satisfies the standard $n\epsilon_{64}\sqrt{a_{ii}a_{jj}}$ bound of \cite{demmel1989cholesky}. In the case of the low-precision floating-point Cholesky decomposition, rather than directly emulating double precision, we sweep many different values of $\epsilon'$ and check the effectiveness of the resulting factorizations as preconditioners for a Preconditioned Conjugate Gradient (PCG) solver. The PCG solves $Ax = b$ until double precision accuracy is achieved, i.e., until $\frac{\|b - Ax\|_2}{\sqrt{n}\|A\|_2\|x\|_2} \le \epsilon_{64}$ is achieved. We check this residual at every iteration.

\begin{table}[h]
\centering
\renewcommand{\arraystretch}{1.5}
\begin{tabular}{c|l|l}
\textbf{Type} & \textbf{Distribution} & \textbf{Eigenvalue Formula} ($\lambda_k,\ k=1 \dots n$) \\ \hline
1 & Geometric & $\lambda_k = \kappa^{-\frac{k-1}{n-1}}$ \\
2 & Arithmetic & $\lambda_k = 1 - \frac{k-1}{n-1}\left(1 - \frac{1}{\kappa}\right)$ \\
3 & Random & $\lambda_1=1, \lambda_n=\frac{1}{\kappa}, \lambda_k \sim \text{Uniform}(\frac{1}{\kappa}, 1)$ \\
4 & Half-Half & $\lambda_k = 1$ for $k \le \frac{n}{2}$,  $\lambda_k = \frac{1}{\kappa}$ for $k > \frac{n}{2}$
\end{tabular}
\caption{Eigenvalue distributions used for our synthetic matrices of dimension $n \times n$ and condition number $\kappa$.}
\label{tab:eigen_dists}
\end{table}
\subsection{Ozaki 1 Emulation on Synthetic Matrices and SuiteSparse}\label{sec:oz1-exp}
\begin{figure*}[tp]
  \centering
  \includegraphics[width=\linewidth]{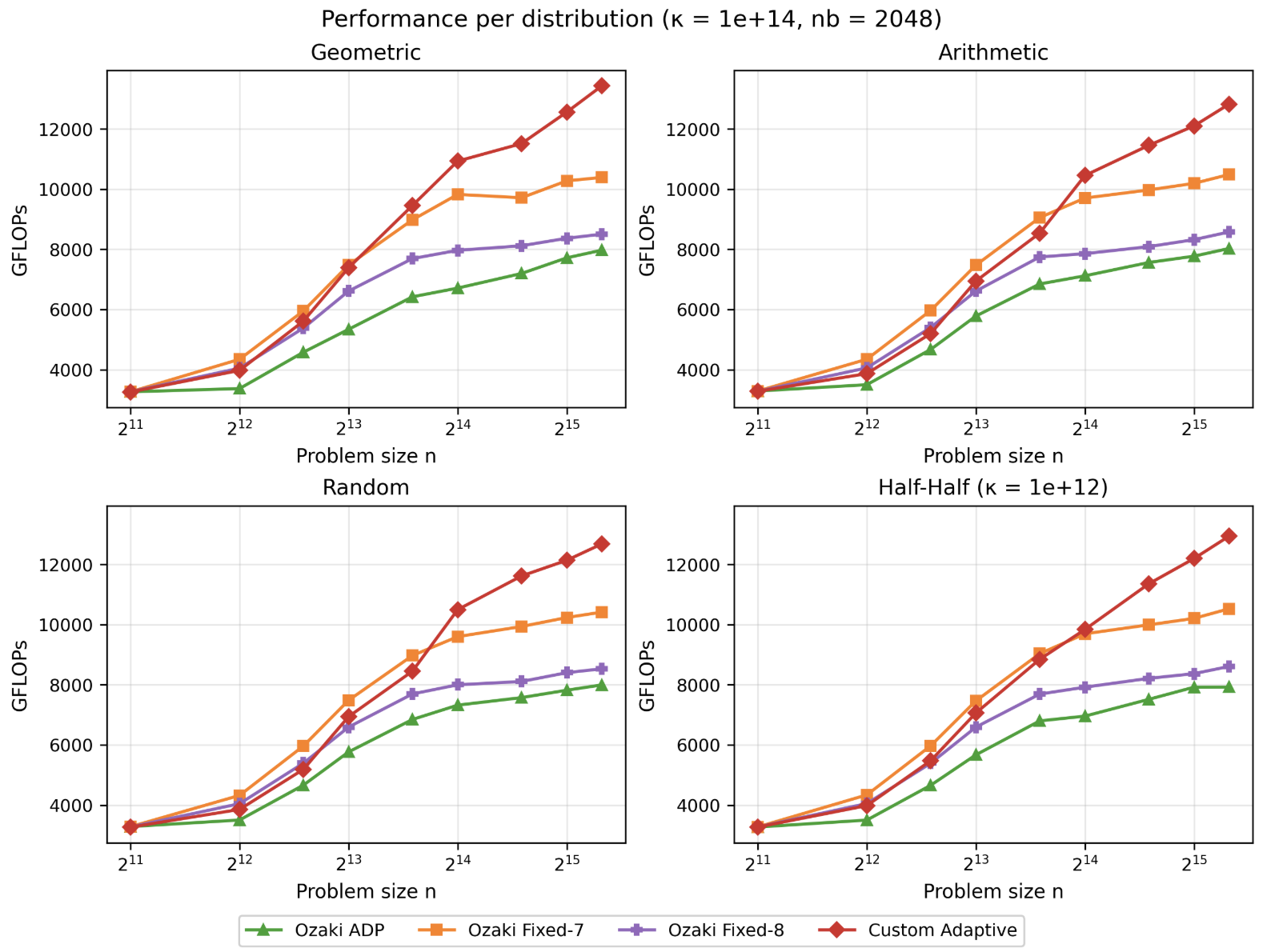}
  \caption{Performance in GFLOPs = $(n^3/3)/\text{factorization\_time}$, of the Ozaki 1 emulation schemes on the synthetic distributions.}
  \label{fig:oz1-perf}
\end{figure*}

\begin{figure*}[tp]
  \centering
\includegraphics[width=\linewidth]{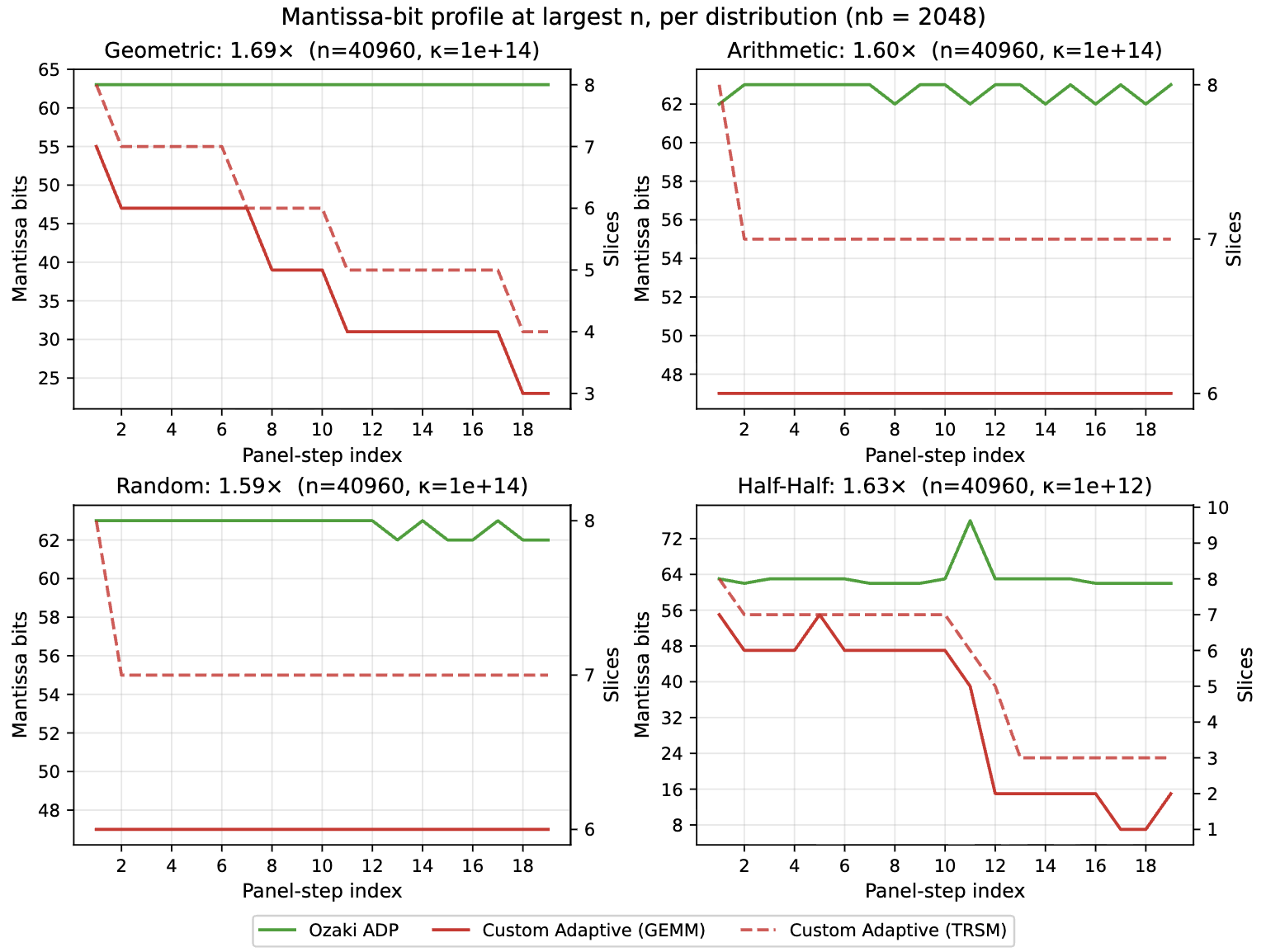}
  \caption{Slice (mantissa bit) profile per panel for Ozaki 1 on the synthetic distributions. 
  \jwdnotenew{Does the Ozaki ADP \#bits really
  just change by 1 as shown in green? Not sure what this
  means.} \svknotenew{I think the cuBLAS ADP algorithm returns a number of bits rather than a number of slices. They then calculate the number of slices needed to store those number of bits. I use this function to get the required number of bits (removed cuBLAS link becuase it caused an arXiv error).}
  \jwdnotenew{A little confusing that the Slices and
  Panel-step index axes are labeled with non-integer
  values, when they can only be integers.} \svknotenew{Thanks. Fixed.}}
  \label{fig:oz1-bits}
\end{figure*}

\begin{figure*}[tp]
  \centering
  \includegraphics[width=\linewidth]{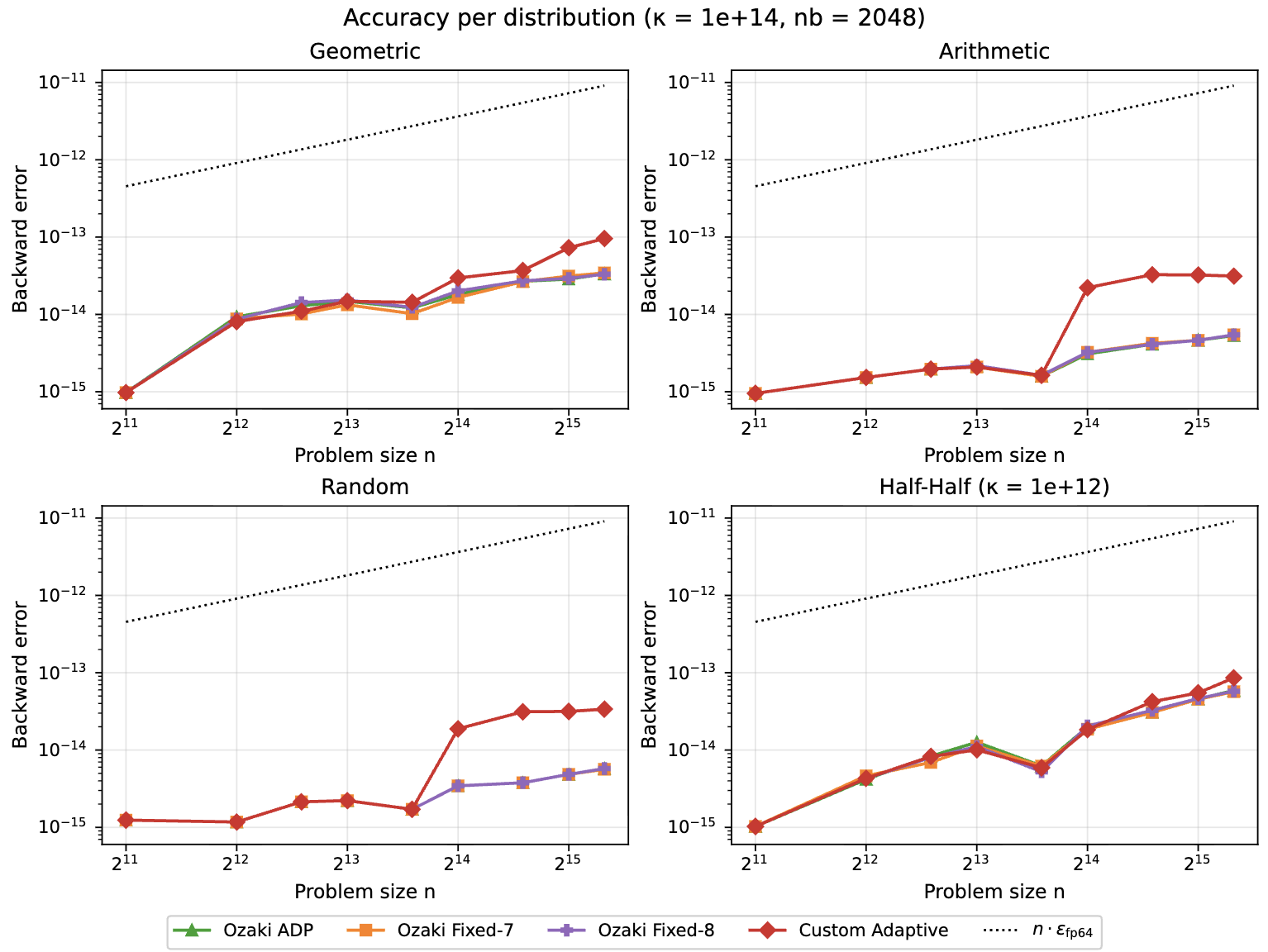}
  \caption{Backward error, $\max_{i,j}|a_{ij} - \sum_k l_{ik}l_{jk}|/\sqrt{a_{ii}a_{jj}}$, of the Ozaki 1 emulation schemes on the synthetic distributions.}
  \label{fig:oz1-acc}
\end{figure*}

\begin{figure*}[tp]
  \centering
  \includegraphics[width=\linewidth]{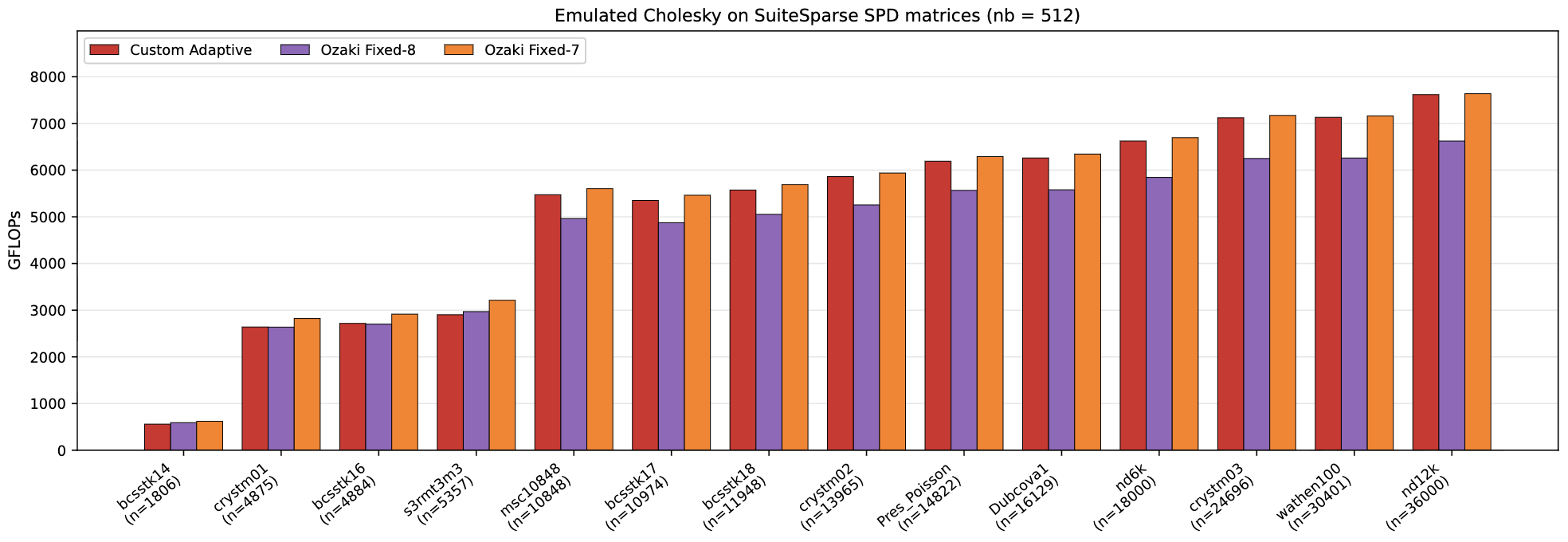}
  \caption{Performance in GFLOPs = $(n^3/3)/\text{factorization\_time}$ of fixed-7, fixed-8 and Adaptive on selected SuiteSparse matrices}
  \label{fig:oz1-suitesparse}
\end{figure*}

We first discuss the performance of the blocked Cholesky factorization with an adaptive number of slices based on the Ozaki 1 scheme, described in Algorithm~\ref{alg:adaptive}.
We use the cuBLAS Ozaki 1 implementation and our experiments are run on a GH200 with the CUDA \jwdnotenew{CUDA} toolkit version 13.2. CUDA 13.2 provides emulated DSYRK  only on Blackwell and not on Hopper. However, emulated DGEMM is available on Hopper. Thus, in our experiments, the Cholesky factorization is implemented using emulated DGEMM instead of DSYRK. Because our DGEMM-based implementation does not exploit symmetry, it performs twice as many floating-point operations as a DSYRK-based implementation. Although this roughly doubles the wall-clock times, the relative speedups should remain comparable if DGEMM is replaced by DSYRK. We compare four emulation variants: fixed 7 slice Ozaki 1, fixed 8 slice Ozaki 1, cuBLAS Automatic Dynamic Precision (\jwdnotenew{Say ``Ozaki ADP'' too to match the label in Fig. 1.} Ozaki ADP) and our adaptive factorization from Algorithm~\ref{alg:adaptive} (labeled Custom Adaptive in Figures \ref{fig:oz1-perf}, \ref{fig:oz1-bits}, \ref{fig:oz1-acc}, \ref{fig:oz1-suitesparse}). 
\jwdnotenew{Specify that this is labeled
Custom Adaptive in Fig. 1.}
As all the emulation variants are slower than the traditional floating-point DGEMM available on Hopper, we restrict our performance comparison to emulation based factorizations. CUBLAS\_EMULATION\_STRATEGY is set to eager so that DGEMMs and DTRSMs always use the Ozaki 1 scheme. For the synthetic matrices in Section~\ref{sec:oz1-exp}, Algorithm \ref{alg:adaptive} uses a block size of  $2048$. Figure~\ref{fig:oz1-perf} displays the performance of all four algorithms for the synthetic matrices in our test set as their dimension $n$ increases. We use $n = 2048 * i$ where $i \in [1,2,3,4,6,8,12,16,20]$.
\jwdnotenew{Are all values of $n$ multiples
of $n_b=2048$? Or does $n=40000$ appear too?} \svknotenew{40000 appears too, here are all the numbers I used - [2048, 4096, 6144, 8192, 12288, 16384, 24576, 32768, 40000]}\jwdnotenew{You should list these 
dimensions somewhere, perhaps in the figure caption.
Also write ``6144 = 3*2048'' etc. to clarify which
are multiples of $n_b$.}
The condition number is $10^{14}$ for the first three distributions of Table~\ref{tab:eigen_dists}, and $10^{12}$ for the half-half distribution (since $10^{14}$ failed in this case even for native FP64).  Figure~\ref{fig:oz1-bits} presents the number of slices (or bits) selected for each GEMM and TRSM by our adaptive algorithm and the number of slices (bits) selected by cuBLAS ADP for each GEMM. We report both the number of slices and the number of bits, since the cuBLAS API tracks the number of bits required for Ozaki 1. We refer to the sequence of slice counts used for GEMM and TRSM at successive iterations as the bit profile of the matrix. Figure \ref{fig:oz1-bits} displays the bit profile of matrices with $n = 40960$ and the largest condition number for which the Cholesky factorization did not fail in our experiments. \jwdnotenew{So 20*2048=40960 failed?} \svknotenew{No, I just didn't run it. I forgot that every other point was a multiple of 2048. I can try running 40960 too, shouldn't take too long.} The number on top of each \jwdnotenew{each} graph represents the speedup of our adaptive algorithm over cuBLAS ADP. As can be seen from these bit profiles, ADP chooses $8-10$ slices for all GEMMs in the factorization, while the adaptive algorithm chooses fewer slices when the matrix has  more rapidly decaying eigenvalues, as in the case of geometric and half-half eigenvalue distributions. Even for the more "uniform" cases of arithmetic and random distributions, the adaptive algorithm chooses a much smaller slice count for GEMMs, and uses mostly $7$ slices for TRSMs, and thus speedups are obtained for these cases as well. 
Figure~\ref{fig:oz1-acc} displays the backward error obtained by the different Ozaki 1 emulation schemes. Our algorithm leads to a backward error smaller than the limit $n\epsilon_{64}$. The loss in accuracy relative to native FP64 is higher for the arithmetic and random distributions compared to the geometric and half-half distributions.

\begin{figure*}[tp]
  \centering
  \includegraphics[width=\linewidth]{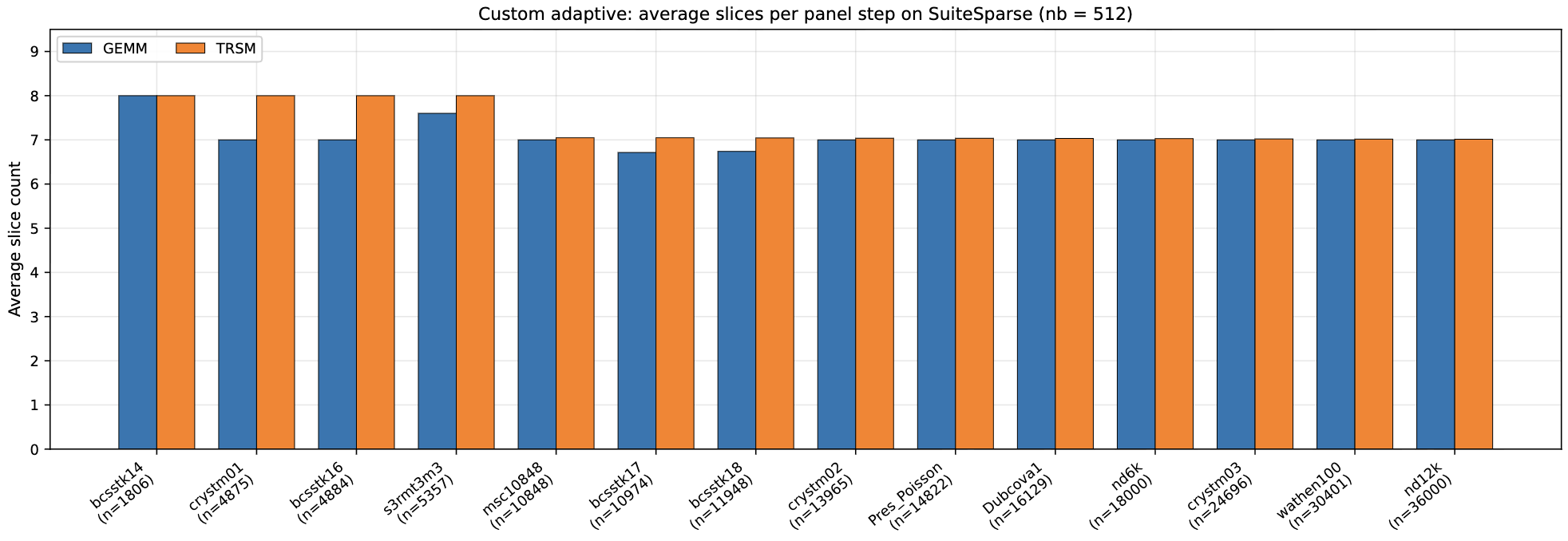}
  \caption{average number of slices for GEMMs and TRSMs on SuiteSparse matrices}
  \label{fig:oz1-suitesparse-bitprof}
\end{figure*}

Figure~\ref{fig:oz1-suitesparse} \jwdnotenew{You might
want to number the figures in the same order in
which the text discusses them.} \svknotenew{Agreed. Fixed, thanks.} displays the results obtained for matrices from the SuiteSparse collection. For the SuiteSparse experiments, we use Algorithm \ref{alg:adaptive} with a block size of 512. Ozaki Fixed-7 is the fastest, Ozaki Fixed-8 is the slowest, and for matrices with
  $n > 10^4$, Ozaki Fixed-7 is $1.00$ to $1.02$ times faster than Custom Adaptive, and
  $1.12$ to $1.15$ times faster than Ozaki Fixed-8. Additionally, all 3 variants produce nearly identical backward errors,
  in the range from $1.4 \times 10^{-4}$ to $1.4 \times 10^{-2}$ times the target upper bound of
  $n\epsilon_{64}$. We do not include a slice count for cuBLAS ADP since it returned very high slice counts for the SuiteSparse matrices, which caused the algorithm to default to regular FP64. This is likely because ADP is run for each GEMM and returns a larger number of slices when the dynamic range of each matrix ($\max_{ij}|a_{ij}|/\min_{ij : a_{ij} \neq 0}|a_{ij}|$) is large.  For most of these matrices, the number of GEMM slices used in our adaptive algorithm remains close to $7$. The number of TRSM slices stays at 8 for some of the SuiteSparse examples (bcsstk14, crystm01, bcsstk16 and s3rmt3m3), but falls to 7 for all other SuiteSparse matrices that we tested. In Figure \ref{fig:oz1-suitesparse-bitprof} we display the average number of slices used for GEMMs and TRSMs for the entire factorization. Note that the average is unweighted. So, the average slice count for the GEMMs is computed as $\frac{1}{\lceil n/n_b \rceil}\sum_{K=1}^{\lceil n/n_b \rceil} \pi_g(K)$. \jwdnotenew{Do you mean it 
occasionally chose more slices? How many?} \svknotenew{The first TRSM uses 8 slices because that's the default, but other than that, it is pretty much 7 slices throughout. I think there is one matrix the last GEMM uses 3 slices, and I think there are a few where we see 6 slices for the last few iterations. But only 7 slices for both GEMM and TRSM other than those minor exceptions.} \svknotenew{Correction : TRSM uses 8 slices for almost all iteratons. Added a new plot.}
Thus, the performance of Algorithm \ref{alg:adaptive} on SuiteSparse matrices corresponds to that of an algorithm that uses a fixed number of $7$ slices.

\subsection{Algorithm~\ref{alg:adaptive-FP} with PCG on Synthetic Matrices and SuiteSparse}

\begin{figure*}[tp]
  \centering
  \includegraphics[width=0.95\linewidth]{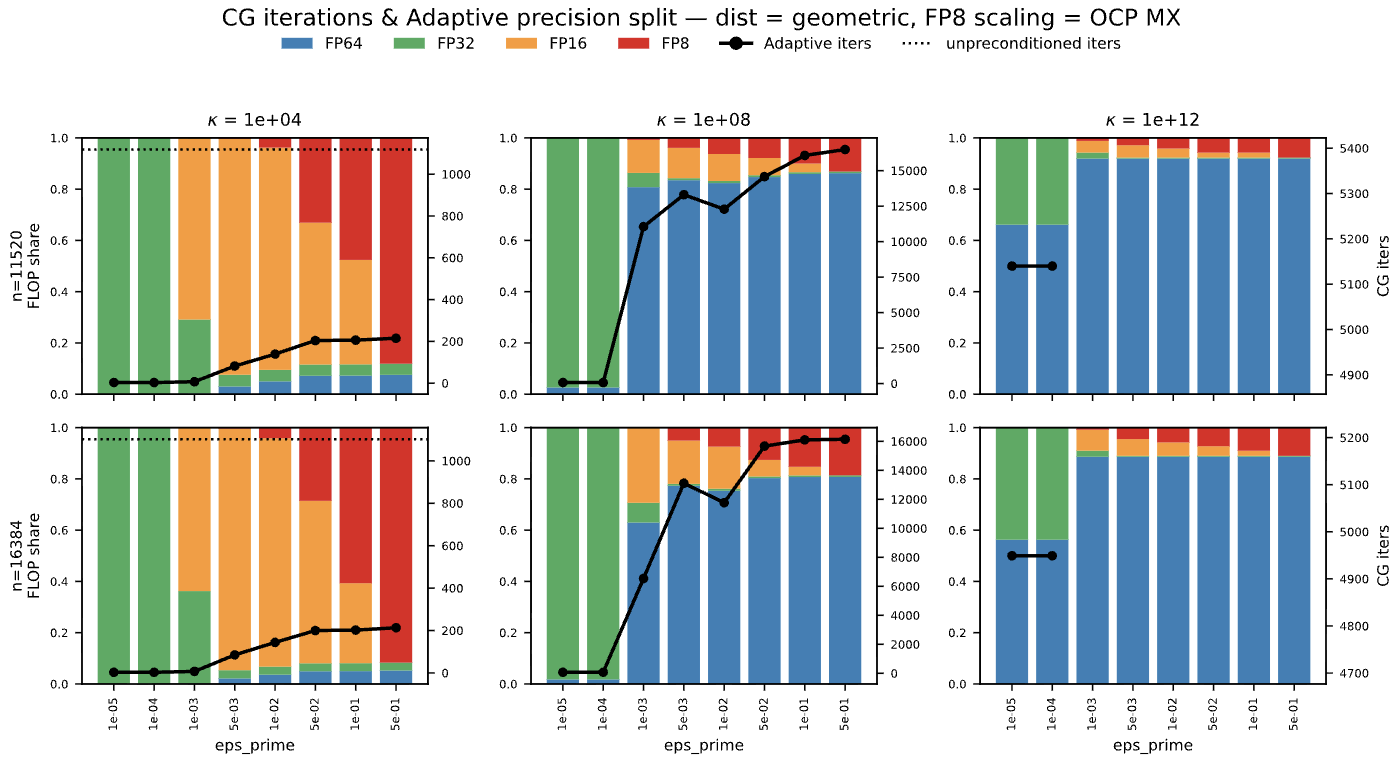}
  \includegraphics[width=0.95\linewidth]{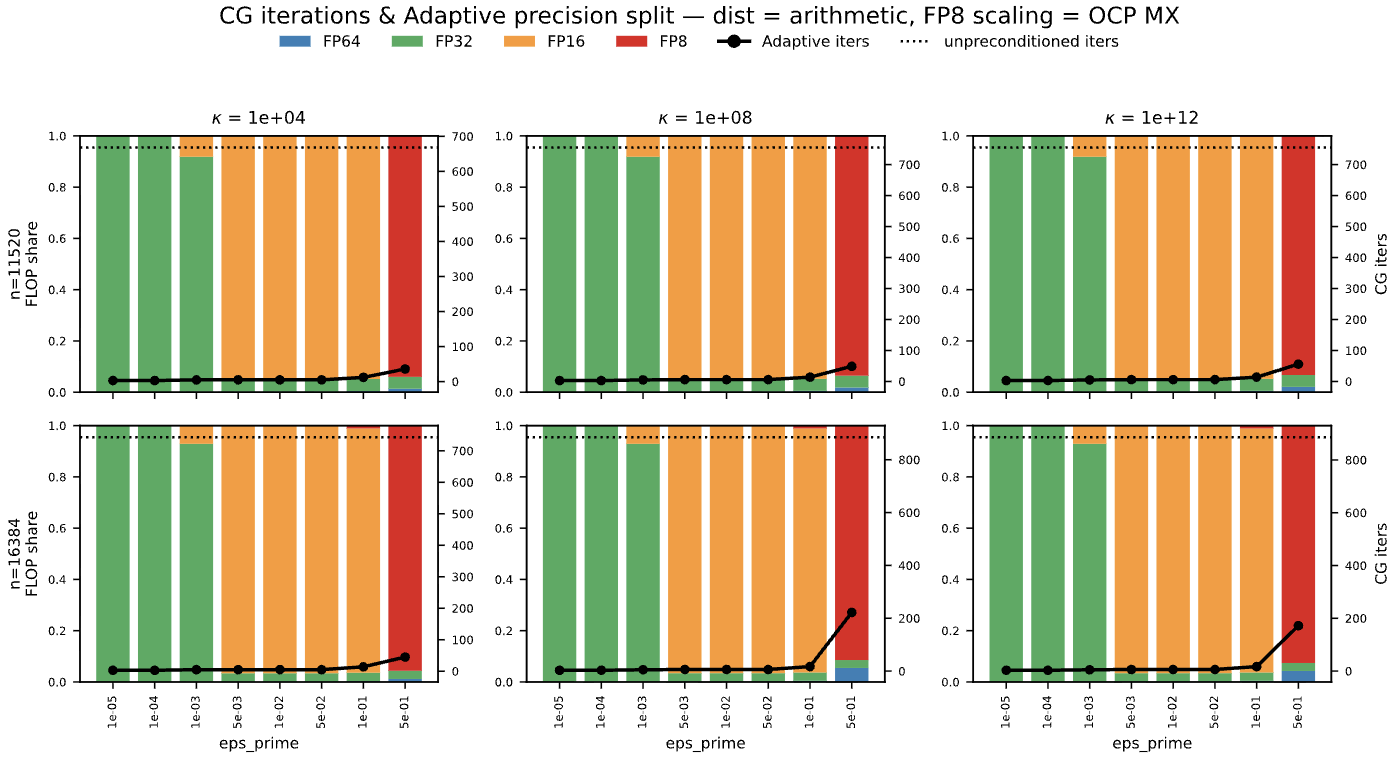}
  \caption{PCG iterations and FLOP share per precision for the Geometric (top) and Arithmetic (bottom) distributions. FLOP share is computed as (\#FLOPs in precision X)/(total number of FLOPs)}
  \label{fig:pcg-synth}
\end{figure*}

In the case of the low-precision floating-point Cholesky decomposition in Algorithm~\ref{alg:adaptive-FP}, we use the produced factorization as a preconditioner for the preconditioned Conjugate Gradient (PCG) to solve a linear system $Ax = b$. For the SuiteSparse matrices, before factorization, the matrix $A$ is equilibrated using a diagonal matrix $D$ whose diagonal elements $d_{ii}$ are chosen to be the power of two closest to $\sqrt{a_{ii}}$. The equilibrated matrix $H = D^{-1}AD^{-1}$ is then factored using Algorithm~\ref{alg:adaptive-FP} to obtain $\hat{L}$. We do not equilibrate the synthetic matrices. For the synthetic matrices, CG on equilibrated matrices does not yield much benefit over CG on the original matrix. We consider microscaling with two different block sizes and scaling formats for the case where FP8 (E4M3)
\jwdnotenew{Add "(E4M3)", or whatever is correct,
to remind the reader.}
is used in Algorithm~\ref{alg:adaptive-FP}. The first strategy, described in Section~\ref{sec:mx}, is OCP MXFP8-style scaling with a block size of $32$ along each row and E8M0 scale factors. Since cuBLAS does not support this on Hopper, we use the LoFloat simulator~\cite{LoFloat}, along with Algorithm \ref{alg:simulated MX GEMM} to simulate the MX GEMM. The second strategy uses a single FP32 scale factor for each row of the panel of $L$ (this scaled GEMM is available in cuBLASLt). For both scaling methods, $\beta$ is chosen as 448, the largest finite number in E4M3. The amount of underflow per row at $n = 16384$ for both strategies was smaller than $10$ for all synthetic test matrices, and we did not see any difference in the quality of the preconditioner. All our experiments implement Algorithm \ref{alg:adaptive-FP} + GMW81 with a block size of $256$. We let the PCG run for either 30000 iterations or until we reach $\frac{\|b - Ax\|_2}{\sqrt{n}\|A\|_2\|x\|_2} \le \epsilon_{64}$. We set FP32 to be the highest available precision and perform panel factorization in FP32. Finding the best $\epsilon'$ is a tuning problem that is beyond the scope of this paper. 
The iteration count and the distribution of FLOPs across Cholesky and PCG for synthetic matrices is shown in Figure~\ref{fig:pcg-synth}. We measure the total number of FLOPs in each precision plus the number of PCG iterations needed for convergence. If there is no point in a bar, that means that the algorithm failed to converge within 30000 iterations. We see that an FP8 preconditioner leads to quick convergence for the arithmetic cases, and it can somewhat speed up convergence for the geometric case if the condition number is small enough. However, for large condition numbers, FP8 fails to produce a good preconditioner when singular values have the geometric distribution. We don't present the random case since it behaves the same way as the arithmetic case, and we don't present the half-half case since it is very easy for CG to solve without any preconditioning.  \jwdnotenew{Confused: Fig 6 shows arithmetic and geometric, not random cases.
And FP8 preconditioner doesn't help for geometric
with $\kappa \geq 1e8$.}
These results suggest that a preconditioner produced using FP8 GEMMs may be useful for some matrices as long as the FP8 GEMMs are fast enough to make up for the cost of increased number of PCG iterations performed in double precision. 

\begin{figure*}[tp]
  \centering
  \includegraphics[width=1.1\linewidth]{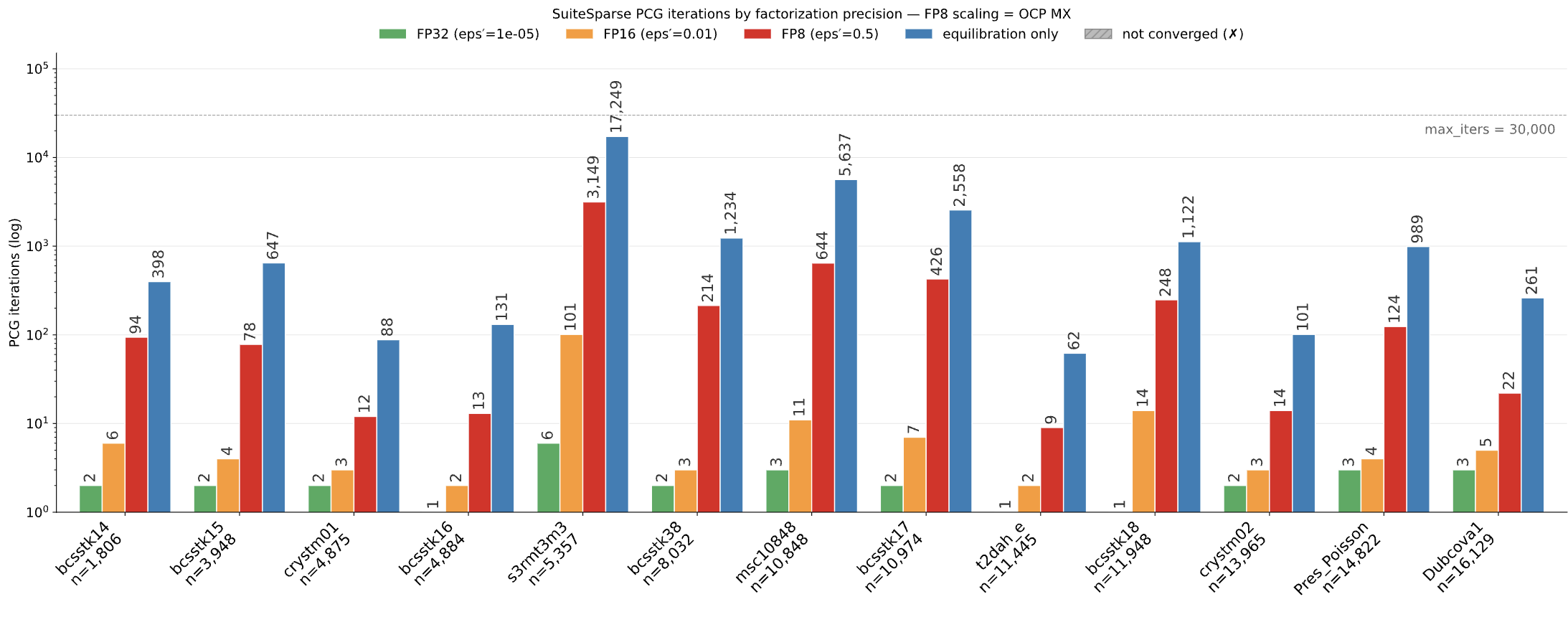}
  \caption{Number of PCG iterations for SuiteSparse matrices. \svknote{Change font of matrix names and numbers}}
  \label{fig:pcg-suitesparse}
\end{figure*}

For SuiteSparse, we test on 13 matrices with $n < 16384$ using different $\epsilon'$. The list of matrices is bcsstk14, bcsstk15, crystm01, bcsstk16, s3rmt3m3, bcsstk38, msc10848, bcsstk17, t2dah\_e, bcsstk18, crystm02, Pres\_Poisson, Dubcova1. Similarly to the Ozaki 1 case, the SuiteSparse matrices do not show any interesting distribution of precisions, and the factorizations choose mainly one precision for the whole factorization. For SuiteSparse, $\epsilon'$ falls into three possible ranges: $\epsilon' \le 10^{-3}$ where all GEMMs are computed in single precision, $10^{-3} < \epsilon' < 0.125$ where all GEMMs are in FP16, and $\epsilon' \ge 0.125$ where all GEMMs are in FP8. So, rather than presenting the number of iterations for each $\epsilon'$, we present the number of iterations for each precision. We find that the combination of MXFP8 and GMW81 is able to form an effective preconditioner for many matrices. Figure~\ref{fig:pcg-suitesparse} displays the number of PCG iterations required for each precision to reach double precision on all the SuiteSparse matrices mentioned above. As seen in Figure \ref{fig:pcg-suitesparse}, using low-precision GEMMs along with diagonal perturbation can form an effective preconditioner for Conjugate Gradient on many SuiteSparse matrices. 

In our experiments, Algorithm \ref{alg:adaptive-FP} uses the GMW81 path described in Section~\ref{sec:failures} since it outperformed diagonal clamping. We use the "optimal" $\beta^{\mathrm{GMW}}$ suggested in \cite{GillMurrayWright1981}. GMW81 outperformed the diagonal For synthetic matrices, we set $\delta^{\mathrm{GMW}}_{\mathrm{FP8}}$ and $\delta^{\mathrm{GMW}}_{\mathrm{FP16}}$ to the unit roundoff of each format ($2^{-4}$ and $2^{-11}$). Doing this for FP32 led to some failures for the ill-conditioned geometric synthetic matrices. However, increasing  $\delta^{\mathrm{GMW}}_{\mathrm{FP32}}$ to $2^{-18}$ seems to fix this. We use the same $\delta^{\mathrm{GMW}}_{\mathrm{FP32}}$ for SuiteSparse. To choose $\delta^{\mathrm{GMW}}_{FP8}$ for SuiteSparse, we sweep [1.0, 0.5, 0.3, 0.1] and find that $\delta^{\mathrm{GMW}}_{FP8} = 0.5$ works best. For $\delta^{\mathrm{GMW}}_{\mathrm{FP16}}$, we sweep $[10^{-3},\, 3\times 10^{-3},\, 10^{-2},\, 3\times 10^{-2}]$ and find that $3\times 10^{-3}$ works best across our test set. Note that the best value of $\delta^{\mathrm{GMW}}$ can depend on the matrix. For example, $\delta^{\mathrm{GMW}}_{\mathrm{FP8}} = 0.3$ gives a better preconditioner than $\delta^{\mathrm{GMW}}_{\mathrm{FP8}} = 0.5$ on many of the tested SuiteSparse matrices, but on bcsstk38 it is worse than CG on the equilibrated A.

\clearpage
\section{Conclusions and Future Work}
We have presented an analysis of floating-point mixed-precision Cholesky and a Cholesky performed via Ozaki 1 fixed-point emulation. We have also presented an algorithm to dynamically reduce precision/slices throughout the Cholesky factorization. For emulation, we show up to 29\% speedup 
\jwdnotenewnew{Abstract refers to
a 20\% speedup.}
over using 7 fixed slices throughout the factorization, and up to 69\% speedup over cuBLAS ADP. We also use the adaptive Cholesky, with a diagonal perturbation to handle factorization failures, to construct preconditioners for PCG: on the 13 SuiteSparse matrices we tested, an FP8 preconditioner reduces the number of PCG iterations by up to $12\times$ over CG on our test matrices. Future work includes constructing stress tests to see if the error bounds in Section~\ref{sec:mp-chol} are attainable. Tighter bounds could let us use low precision more often, thus giving us more speedup. Developing a high performance implementation and adapting our algorithm to the left-looking case is also of interest. Applying our algorithms to CholeskyQR, or developing similar algorithms for LU also seem like promising directions for future work.

\section{AI usage disclosure}
Claude Opus 5 and Claude Opus 5.5 were used to improve the writing in some places, and get rid of notational collisions throughout the paper. The LLMs were also able to find an error in the proof of Theorem 4.2 which allowed the authors to remove an extra $(1 + \epsilon_h)$ factor from the bound. They were also used to unify the proofs of Theorem 4.2 and 4.3 given an initial proof sketch of how to unify the old versions of 4.2 and 4.3 (which treated both MX and no MX cases separately). LLMs were also used to set up and run the experiments on the Alps cluster. ChatGPT 6 Astra-Pro and Claude Opus 5.5 were also used to review the paper. All code and proofs were checked by the authors.

\begin{acks}
\svknote{TODO: Add Alps Acknowledgment from Laura}

The experiments performed in this work on Alps were supported through a grant from the Swiss National Supercomputing Centre (CSCS) under project RandESC. \svknotenew{Acknowledgement spelling error}
\end{acks}

\section{Appendix}

\subsection{Adaptive precision with Ozaki 2}
In the appendix, we derive an adaptive algorithm for Ozaki 2 using the same methods we used in Section \ref{sec:adaptive_alg}. Results from our experiments are presented in Section~\ref{sec:oz2-exp}. The resulting algorithm is Algorithm~\ref{alg:adaptive-oz2}. For the GEMM error bound, we use only the dominant $O(t)$ terms from Equation (16) in \cite{uchino2026erroranalysismatrixmultiplication}. The paper uses, $t = \frac{1}{\sqrt{2^5(\prod_{\ell=1}^{N} p_\ell - 1)}}$, where $p_\ell$ is the $\ell$-th largest coprime in the list of coprime numbers from Equation (10) of \cite{uchino2026erroranalysismatrixmultiplication}, and $N$ being the number of coprimes used by the GEMM. We are interested in controlling how many co-prime numbers we need at any given iteration of Algorithm~\ref{A_chol}. So, rather than just $t$, we will use the notation $t(N) = \frac{1}{\sqrt{2^5(\prod_{\ell=1}^{N} p_\ell - 1)}}$. Now, to understand the error analysis of a Cholesky with Ozaki 2, it is sufficient to see what the error on the diagonal, $|e_{ii}| = |a_{ii} - \sum_{k=1}^{i}\hat{l}_{ik}^2|$ looks like. The rest of the proof will follow the exact same steps as the proof of Theorem~\ref{thm:chol-emu}. The dominant term of the bound from \cite[Eq.~(16)]{uchino2026erroranalysismatrixmultiplication} is $t(N)\bigl(|A|\,v\,(2^{\beta'})^T + 2^{\alpha'}v^T|B|\bigr)$, with $v = (1, \ldots, 1)^T$. Where, all terms are defined as in \cite{uchino2026erroranalysismatrixmultiplication}.  Since we are concerned with GEMMs of the form $\hat{L}_{IK}\hat{L}_{IK}^T$, $A = B^T = \hat{L}_{IK}$, $\alpha = \beta$. For $A \in \mathbb{R}^{m \times k}$ and $A^T = B$
\[
\begin{gathered}
    \alpha'_i = \alpha_i + \tfrac{1}{2}q_i, \qquad \alpha_i = \Bigl\lfloor \log_2 \max_{1 \le h \le k}|a_{ih}| \Bigr\rfloor, \qquad q_i = \log_2 \max_{1 \le j \le n} \bar{c}_{ij}, \\
    \bar{C} = \bigl\lceil \mathrm{diag}(2^{5 - \alpha})|A| \bigr\rceil \bigl\lceil |A^T|\,\mathrm{diag}(2^{5 - \alpha}) \bigr\rceil,
\end{gathered}
\]
where $q$ is the vector $e'$ of \cite{uchino2026erroranalysismatrixmultiplication}. So, we can write the inequality on the block diagonal as
\begin{equation} \label{eq:oz2-block-diag}
    \Bigl|A_{II} - \sum_{K=1}^{I}\hat{L}_{IK}\hat{L}_{IK}^T\Bigr| \le \sum_{K=1}^{I-1} t(\pi_g(K))\Bigl(|\hat{L}_{IK}|\,v\,\bigl(2^{\alpha'^{(K)}}\bigr)^T + 2^{\alpha'^{(K)}}v^T|\hat{L}_{IK}|^T\Bigr) + n_b\epsilon_h \Omega_{II},
\end{equation}
where $\alpha'^{(K)} = \alpha^{(K)} + \frac{1}{2}q^{(K)}$ is the vector $\alpha'$ of \cite{uchino2026erroranalysismatrixmultiplication} for the GEMM with panel $K$, and $q^{(K)}$ denotes the vector $e'$ $(= f')$ of \cite{uchino2026erroranalysismatrixmultiplication}. For row $i$ of block $I$, the $i$-th entry of both $|\hat{L}_{IK}|\,v$ and $v^T|\hat{L}_{IK}|^T$ is $r_i^{(K)} = \sum_{k=(K-1)n_b+1}^{Kn_b}|\hat{l}_{ik}|$, so the $(i,i)$ entry of each of the two terms in \eqref{eq:oz2-block-diag} is $t(\pi_g(K))\,2^{\alpha'^{(K)}_i} r_i^{(K)}$. On the diagonal element, we then get
\begin{equation} \label{eq:oz2-eii}
    |e_{ii}| \le \sum_{K=1}^{I-1} 2\,t(\pi_g(K))\,2^{\alpha'^{(K)}_i} r_i^{(K)} + n_b\epsilon_h a^u_{ii} = \sum_{K=1}^{I-1} 2\,t(\pi_g(K))\,2^{\alpha_i^{(K)}}\,2^{q_i^{(K)}/2}\, r_i^{(K)} + n_b\epsilon_h a^u_{ii}.
\end{equation}
Getting an upper bound on $2^{\alpha^{(K)}_i}$ is straightforward, since $\alpha^{(K)}_i = \bigl\lfloor \log_2 \max_{(K-1)n_b < k \le Kn_b}|\hat{l}_{ik}| \bigr\rfloor \implies 2^{\alpha^{(K)}_i} \le g_i^{(K)}$, where $g_i^{(K)} = \max_{(K-1)n_b < k \le Kn_b}|\hat{l}_{ik}|$. Now, $q^{(K)}$ is computed by first scaling $|\hat{L}_{IK}|$ by the diagonal matrix $\mathrm{diag}(2^{5 - \alpha^{(K)}})$, and computing $\bar{C} = \lceil \mathrm{diag}(2^{5 - \alpha^{(K)}})|\hat{L}_{IK}|\rceil\lceil \mathrm{diag}(2^{5 - \alpha^{(K)}})|\hat{L}_{IK}|\rceil^T$. We then compute $q^{(K)}_i = \log_2 \max_{h} \bar{c}_{ih}$. Since $2^{5-\alpha^{(K)}_i}|\hat{l}_{ik}| < 2^6$, every entry of $\lceil \mathrm{diag}(2^{5 - \alpha^{(K)}})|\hat{L}_{IK}|\rceil$ lies in $[0, 2^6]$, so $\bar{c}_{hh} \le 2^{12}n_b$ for every $h$. Since $\bar{C}$ is positive semidefinite, $\bar{c}_{ih} \le \sqrt{\bar{c}_{ii}\bar{c}_{hh}} \le 2^{12}n_b$. Then, we get
\begin{equation} \label{eq:oz2-alpha-prime}
    2^{q^{(K)}_i/2} \le 2^6\sqrt{n_b} \implies 2^{\alpha'^{(K)}_i} \le 2^6\sqrt{n_b}\,g_i^{(K)}.
\end{equation}
Finally, since every $|\hat{l}_{ik}| \le g_i^{(K)}$ in panel $K$, $r_i^{(K)} \le n_b\,g_i^{(K)}$. Substituting this and \eqref{eq:oz2-alpha-prime} into each term of \eqref{eq:oz2-eii} gives
\begin{equation} \label{eq:oz2-gi2}
    2\,t(\pi_g(K))\,2^{\alpha'^{(K)}_i} r_i^{(K)} \le 2\,t(\pi_g(K))\cdot 2^6\sqrt{n_b}\,g_i^{(K)}\cdot n_b\,g_i^{(K)} = 2^7 n_b\sqrt{n_b}\,t(\pi_g(K))\bigl(g_i^{(K)}\bigr)^2.
\end{equation}
The rest follows from the proof of Theorem~\ref{thm:chol-emu}, using $\bigl(g_i^{(K)}\bigr)^2 \le \sum_{k=(K-1)n_b+1}^{Kn_b}\hat{l}_{ik}^2$, with $\gamma_{\mathrm{oz2}}(\psi, \Psi, n, n_b) = 2^7 n_b\sqrt{n_b}\,t(\psi)$ in place of $\gamma_{\mathrm{emu}}$.

\begin{algorithm}
\caption{Blocked Cholesky Factorization with Adaptive Moduli Ozaki 2 GEMM}
\label{alg:adaptive-oz2}
\begin{algorithmic}[1]
\Require Matrix $A \in \mathbb{R}^{n \times n}$ (SPD) in some high precision $h$, block size $n_b$, tolerance $\epsilon'$, max moduli count $\Psi'$
\Ensure Lower Cholesky factor $L$ (in FP64) such that $A = LL^T$ (in-place)
\State $\mathrm{diag}_i \gets a_{ii}$ for $i = 1,\ldots,n$ \Comment{$\mathrm{diag} \in \mathbb{R}^n$: original diagonal of $A$ in $h$}
\State $\mathrm{diag}^u_i \gets a_{ii}$ for $i = 1, \ldots , n$
\Comment{$\mathrm{diag}^u \in \mathbb{R}^n$, in $h$}
\State $m_i \gets 0$ for $i = 1,\ldots,n$ \Comment{$m \in \mathbb{R}^n$ in $h$}
\State $g_i \gets 0$ for $i = 1,\ldots,n$ \Comment{$g \in \mathbb{R}^n$ in $h$}
\State $\psi \gets \Psi'$
\For{$I = 1, \ldots, n/n_b$}
    \State $o \gets (I{-}1)n_b$ 
    \State $w \gets \min(n_b,\, n - o)$ 
    \State $A_{II} \gets A[o{+}1{:}o{+}w,\; o{+}1{:}o{+}w]$ \Comment{diagonal block}
    \State $A_{II} \gets \textsc{Potrf}(A_{II})$ \Comment{$A_{II} \leftarrow L_{II} L_{II}^T$ in $h$}
    \If{$o + w < n$}
        \State $A_{I{+}1,I} \gets A[o{+}w{+}1{:}n,\; o{+}1{:}o{+}w]$ \Comment{off-diagonal panel}
        \State $A_{I{+}1,I} \gets A_{I{+}1,I} \cdot A_{II}^{-T}$  \Comment{TRSM is emulated using $\psi$ moduli}
        \For{$i = o{+}w{+}1, \ldots, n$} \Comment{one pass over row $i$ of the panel}
            \State $g_i \gets \max_{o<k\le o+w} |l_{ik}|$
            \State $\mathrm{diag}^u_i \gets \mathrm{diag}^u_i - \sum_{k=o+1}^{o+w} l_{ik}^2$
        \EndFor
        \For{$p_1 = 1, \ldots , \Psi'$}
        \For{$p_2 = 1, \ldots , p_1$}
            \If{$2^7n_b\sqrt{n_b}\,t(p_2)\, g_i^2 + \frac{\gamma_{\mathrm{oz2}}(p_1, \Psi', n - In_b, n_b)}{1 - \gamma_{\mathrm{oz2}}(p_1, \Psi', n - In_b, n_b)}a_{ii}^u \le \epsilon' a_{ii} - m_i$ for $o+w < i \le n$}
            \State $\pi_g(I) \gets p_2, \psi \gets p_1$
            \State goto line~\ref{line:oz2-after-select}
            \EndIf
        \EndFor
        \EndFor
        \For{$i = o{+}w{+}1, \ldots, n$}\label{line:oz2-after-select}
            \State $m_i \gets m_i + 2^7 n_b\sqrt{n_b}\,t(\pi_g(I))\, g_i^2$
        \EndFor
        \State $A_{I{+}1,I{+}1} \gets A[o{+}w{+}1{:}n,\; o{+}w{+}1{:}n]$ \Comment{trailing submatrix}
        \State $A_{I{+}1,I{+}1} \gets A_{I{+}1,I{+}1} - A_{I{+}1,I}\, A_{I{+}1,I}^T$
        \Comment{GEMM/SYRK is emulated using $\pi_g(I)$ moduli}
    \EndIf
\EndFor
\end{algorithmic}
\end{algorithm}

\subsection{Ozaki 2 Emulation on synthetic matrices}\label{sec:oz2-exp}
\begin{figure*}[tp]
  \centering
  \includegraphics[width=0.85\linewidth]{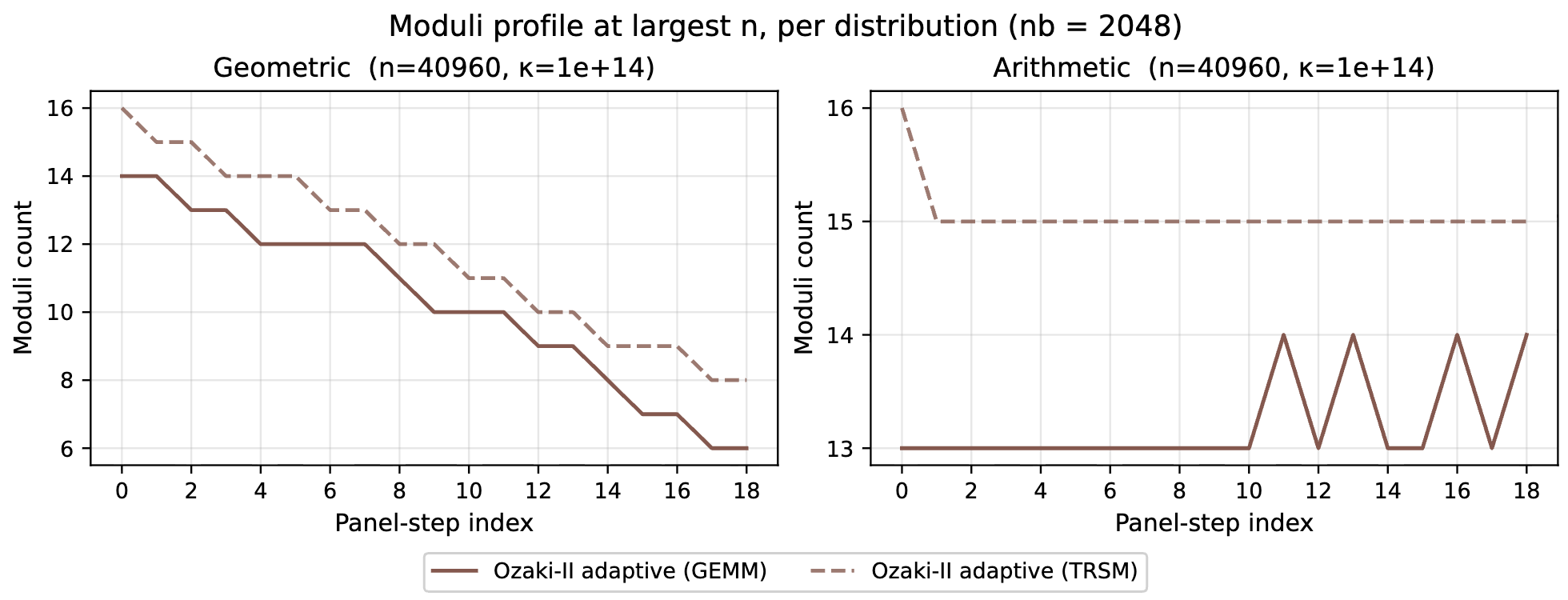}
  \caption{Moduli profile per panel for adaptive Ozaki 2 on the Geometric and Arithmetic distributions.}
  \label{fig:oz2-moduli}
\end{figure*}

Ozaki 2 is also a grade B GEMM so we expect it to have similar properties to Ozaki 1 in terms of accuracy, and thus slice selection. 
We use the implementation from \cite{ozaki2025ozakischemeiigemmoriented} for our experiments. We compare our adaptive algorithm (Algorithm~\ref{alg:adaptive-oz2}) with fixed 12, 13 and 14 moduli implementations. Similar to Ozaki 1, Figure \ref{fig:oz2-moduli} shows that the moduli profile is decreasing for geometric and pretty much constant for arithmetic, with GEMMs being less precise than TRSM. We also stay within the $n\epsilon_{64}$ bound despite only using the $O(t)$ terms to write Algorithm~\ref{alg:adaptive-oz2}.

\begin{figure*}[tp]
  \centering
  \includegraphics[width=\linewidth]{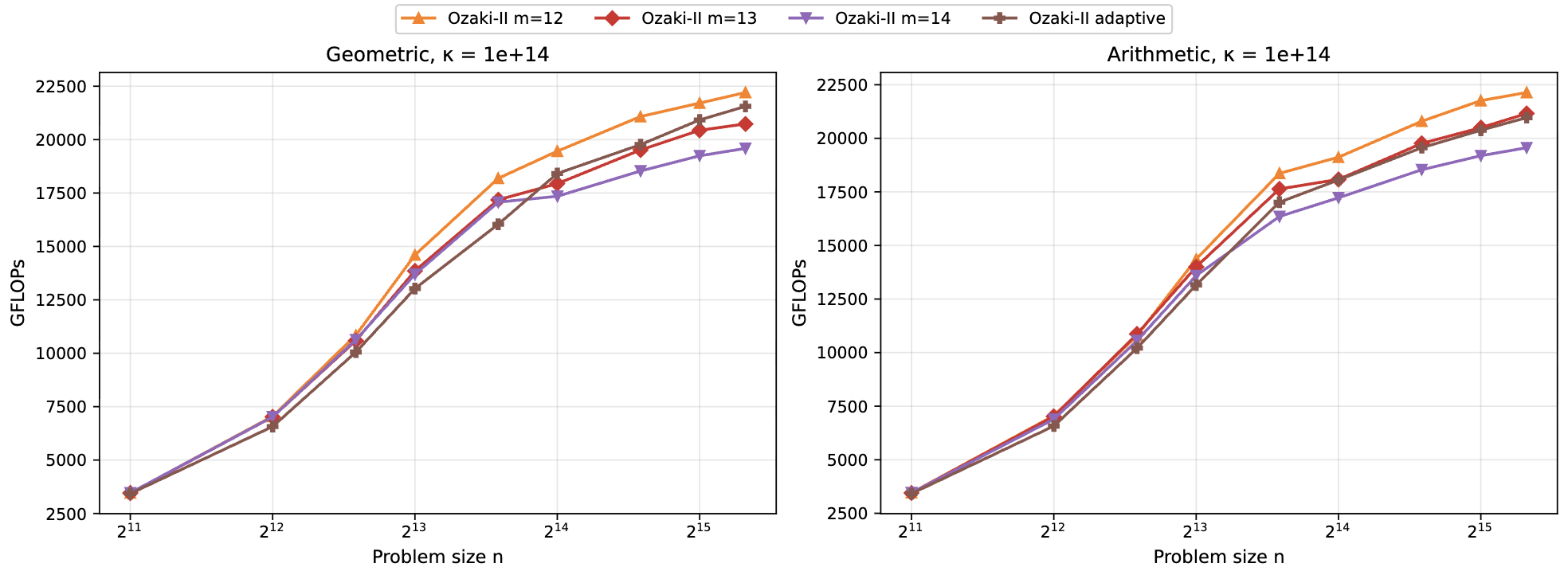}
  \caption{Performance of the Ozaki 2 emulation schemes on the synthetic distributions.}
  \label{fig:oz2-perf}
\end{figure*}

\begin{figure*}[tp]
  \centering
  \includegraphics[width=\linewidth]{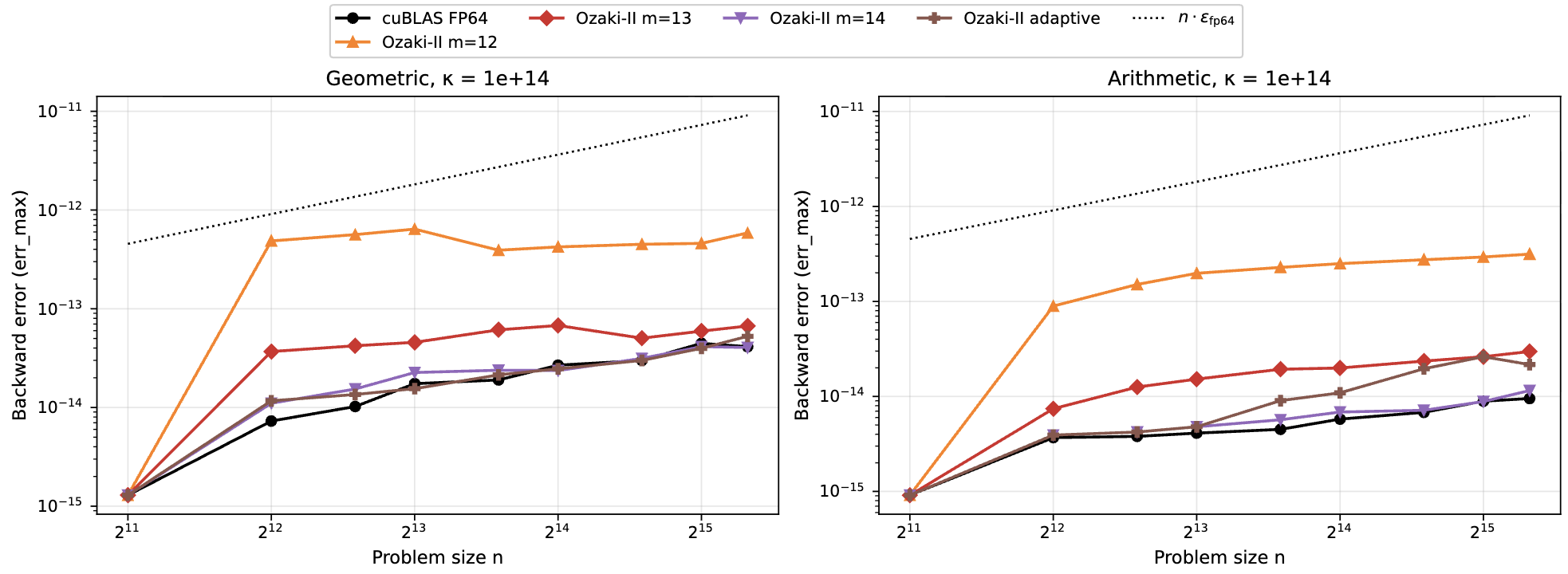}
  \caption{Backward error of the Ozaki 2 emulation schemes on the synthetic distributions.}
  \label{fig:oz2-acc}
\end{figure*}

\clearpage
\printbibliography

@inproceedings{haidar-mixed-prec,
author = {Haidar, Azzam and Wu, Panruo and Tomov, Stanimire and Dongarra, Jack},
title = {Investigating half precision arithmetic to accelerate dense linear system solvers},
year = {2017},
isbn = {9781450351256},
publisher = {Association for Computing Machinery},
address = {New York, NY, USA},
url = {https://doi.org/10.1145/3148226.3148237},
% doi = {10.1145/3148226.3148237},
abstract = {The use of low-precision arithmetic in mixed-precision computing methods has been a powerful tool to accelerate numerous scientific computing applications. Artificial intelligence (AI) in particular has pushed this to current extremes, making use of half-precision floating-point arithmetic (FP16) in approaches based on neural networks. The appeal of FP16 is in the high performance that can be achieved using it on today's powerful manycore GPU accelerators, e.g., like the NVIDIA V100, that can provide 120 TeraFLOPS alone in FP16. We present an investigation showing that other HPC applications can harness this power too, and in particular, the general HPC problem of solving Ax = b, where A is a large dense matrix, and the solution is needed in FP32 or FP64 accuracy. Our approach is based on the mixed-precision iterative refinement technique - we generalize and extend prior advances into a framework, for which we develop architecture-specific algorithms and highly-tuned implementations that resolve the main computational challenges of efficiently parallelizing, scaling, and using FP16 arithmetic in the approach on high-end GPUs. Subsequently, we show for a first time how the use of FP16 arithmetic can significantly accelerate, as well as make more energy efficient, FP32 or FP64-precision Ax = b solvers. Our results are reproducible and the developments will be made available through the MAGMA library. We quantify in practice the performance, and limitations of the approach.},
booktitle = {Proceedings of the 8th Workshop on Latest Advances in Scalable Algorithms for Large-Scale Systems},
articleno = {10},
numpages = {8},
keywords = {mixed-precision iterative refinement, half precision, HPC, GPGPU},
location = {Denver, Colorado},
series = {ScalA '17}
}

@techreport{demmel1989cholesky,
  author    = {J. W. Demmel},
  title     = {On Floating-Point Errors in {C}holesky},
  institution = {University of Tennesse at Knoxville},
  number    = {CS-89-87},
  note = {{LAPACK Working Note 14}},
  month     = {October},
  year      = {1989}
}

@article{five-prec-GMRES,
author = {Amestoy, Patrick and Buttari, Alfredo and Higham, Nicholas J. and L’Excellent, Jean-Yves and Mary, Theo and Vieubl\'{e}, Bastien},
title = {Five-Precision GMRES-Based Iterative Refinement},
journal = {SIAM Journal on Matrix Analysis and Applications},
volume = {45},
number = {1},
pages = {529-552},
year = {2024},
% doi = {10.1137/23M1549079},

% URL = {https://doi.org/10.1137/23M1549079},
eprint = { 
    
        https://doi.org/10.1137/23M1549079
    
    

}
,
    abstract = { Abstract. GMRES-based iterative refinement in three precisions (GMRES-IR3), proposed by Carson and Higham in 2018, uses a low precision LU factorization to accelerate the solution of a linear system without compromising numerical stability or robustness. GMRES-IR3 solves the update equation of iterative refinement using GMRES preconditioned by the LU factors, where all operations within GMRES are carried out in the working precision \(u\), except for the matrix–vector products and the application of the preconditioner, which require the use of extra precision \(u^2\). The use of extra precision can be expensive, and is especially unattractive if it is not available in hardware; for this reason, existing implementations have not used extra precision, despite the absence of an error analysis for this approach. In this article, we propose to relax the requirements on the precisions used within GMRES, allowing the use of arbitrary precisions \(u\_p\) for applying the preconditioned matrix–vector product and \(u\_g\) for the rest of the operations. We obtain the five-precision GMRES-based iterative refinement (GMRES-IR5) algorithm which has the potential to solve relatively badly conditioned problems in less time and memory than GMRES-IR3. We develop a rounding error analysis that generalizes that of GMRES-IR3, obtaining conditions under which the forward and backward errors converge to their limiting values. Our analysis makes use of a new result on the backward stability of MGS-GMRES in two precisions. On hardware where three or more arithmetics are available, which is becoming very common, the number of possible combinations of precisions in GMRES-IR5 is extremely large. We provide an analysis of our theoretical results that identifies a relatively small subset of relevant combinations. By choosing from within this subset one can achieve different levels of tradeoff between cost and robustness, which allows for a finer choice of precisions depending on the problem difficulty and the available hardware. We carry out numerical experiments on random dense and SuiteSparse matrices to validate our theoretical analysis and discuss the complexity of GMRES-IR5. }
}

@article{squeezing_fp16,
author = {Higham, Nicholas J. and Pranesh, Srikara and Zounon, Mawussi},
title = {Squeezing a Matrix into Half Precision, with an Application to Solving Linear Systems},
journal = {SIAM Journal on Scientific Computing},
volume = {41},
number = {4},
pages = {A2536-A2551},
year = {2019},
% doi = {10.1137/18M1229511},

% URL = {https://doi.org/10.1137/18M1229511},
eprint = { 
    
        https://doi.org/10.1137/18M1229511
    
    

}
,
    abstract = { Motivated by the demand in machine learning, modern computer hardware is increasingly supporting reduced precision floating-point arithmetic, which provides advantages in speed, energy, and memory usage over single and double precision. Given the availability of such hardware, mixed precision algorithms that work in single or double precision but carry out part of a computation in half precision are now of great interest for general scientific computing tasks. Because of the limited range of half precision arithmetic, in which positive numbers lie between \$6\times 10^{-8}\$ and \$7\times 10^4\$, a straightforward rounding of single or double precision data into half precision can lead to overflow, underflow, or subnormal numbers being generated, all of which are undesirable. We develop an algorithm for converting a matrix from single or double precision to half precision. It first applies two-sided diagonal scaling in order to equilibrate the matrix (that is, to ensure that every row and column has \$\infty\$-norm \$1\$), then multiplies by a scalar to bring the largest element within a factor \$\theta \le 1\$ of the overflow level, and finally rounds to half precision. The second step ensures that full use is made of the limited range of half precision arithmetic, and \$\theta\$ must be chosen to allow sufficient headroom for subsequent computations. We apply the new algorithm to GMRES-based iterative refinement (GMRES-IR), which solves a linear system \$Ax = b\$ with single or double precision data by LU factorizing \$A\$ in half precision and carrying out iterative refinement with the correction equations solved by GMRES preconditioned with the low precision LU factors. Previous implementations of this algorithm have used a crude conversion to half precision that our experiments show can cause slow convergence of GMRES-IR for badly scaled matrices or failure to converge at all. The new conversion algorithm computes \$\infty\$-norms of rows and columns of the matrix and its cost is negligible in the context of LU factorization. We show that it leads to faster convergence of GMRES-IR for badly scaled matrices and thereby allows a much wider class of problems to be solved. }
}

@misc{micikevicius2022fp8formatsdeeplearning,
      title={{FP8 Formats for Deep Learning}}, 
      author={Paulius Micikevicius and Dusan Stosic and Neil Burgess and Marius Cornea and Pradeep Dubey and Richard Grisenthwaite and Sangwon Ha and Alexander Heinecke and Patrick Judd and John Kamalu and Naveen Mellempudi and Stuart Oberman and Mohammad Shoeybi and Michael Siu and Hao Wu},
      year={2022},
%       eprint={2209.05433},
%       archivePrefix={arXiv},
      primaryClass={cs.LG},
      url={https://arxiv.org/abs/2209.05433} 
}

@misc{ocp2023mx,
  title        = {{OCP Microscaling Formats (MX) Specification Version 1.0}},
  author       = {{Open Compute Project}},
  organization = {Open Compute Project Foundation},
  year         = {2023},
  month        = {Sep},
  url          = {https://www.opencompute.org/documents/ocp-microscaling-formats-mx-v1-0-spec-final-pdf},
  note         = {Accessed: 2024-12-02}
}

@misc{p3109_public_repo,
  author       = {P3109},
  title        = {Public Repository},
  year         = {2024},
  url          = {https://github.com/P3109/Public/tree/main},
  note         = {Accessed: 2024-12-02},
  howpublished = {GitHub repository, main branch}
}

@book{Higham_main_book,
author = {Higham, Nicholas J.},
title = {Accuracy and Stability of Numerical Algorithms},
publisher = {Society for Industrial and Applied Mathematics},
year = {2002},
% doi = {10.1137/1.9780898718027},
address = {},
edition   = {Second},
% URL = {https://epubs.siam.org/doi/abs/10.1137/1.9780898718027},

eprint = {https://epubs.siam.org/doi/pdf/10.1137/1.9780898718027}
}

@techreport{noune2022numerical,
  author       = {B. Noune and P. Jones and D. Justus and D. Masters and C. Luschi},
  title        = {8-bit numerical formats for deep neural networks},
  year         = {2022},
  institution  = {arXiv cs.LG},
  url          = {https://arxiv.org/abs/2206.02915}
}

@techreport{tesla2023dojo,
  author       = {{Tesla, Inc.}},
  title        = {Tesla {Dojo} {Technology}: {A}
                         guide to {T}esla's configurable floating point formats and arithmetic},
  year         = {2023},
  institution  = {Tesla, Inc.},
  url          = {https://web.archive.org/web/20230503235751/https://tesla-cdn.thron.com/static/MXMU3S_tesla-dojo-technology_1WDVZN.pdf}
}

@article{McSweeney2020,
  author={P. {McSweeney}},
  title={Modified {Cholesky} Decomposition and Applications},
  year={2020},
  url={https://mcsweeney90.github.io/files/modified-cholesky-decomposition-and-applications.pdf}
}

@inproceedings{Schwarz_2026, series={SCA/HPCAsia 2026},
   title={{Guaranteed DGEMM Accuracy While Using Reduced Precision Tensor Cores Through Extensions of the Ozaki Scheme}},
   url={http://dx.doi.org/10.1145/3773656.3773670},
   DOI={10.1145/3773656.3773670},
   booktitle={Proceedings of the Supercomputing Asia and International Conference on High Performance Computing in Asia Pacific Region},
   publisher={ACM},
   author={Schwarz, Angelika and Anders, Anton and Brower, Cole and Bayraktar, Harun and Gunnels, John and Clark, Kate and Xu, RuQing G. and Rodriguez, Samuel and Cayrols, Sebastien and Tabaszewski, Pawel and Podlozhnyuk, Victor},
   year={2026},
   month=jan, pages={91–101},
   collection={SCA/HPCAsia 2026} }

@misc{ootomo2024dgemmintegermatrixmultiplication,
      title={DGEMM on Integer Matrix Multiplication Unit}, 
      author={Hiroyuki Ootomo and Katsuhisa Ozaki and Rio Yokota},
      year={2024},
      eprint={2306.11975},
      archivePrefix={arXiv},
      primaryClass={cs.DC},
      url={https://arxiv.org/abs/2306.11975}, 
}

@techreport{Ootomo_Ozaki_Yokota_2024,
author={H. Ootomo and K. Ozaki and R. Yokota},
title={{DGEMM} on Integer Matrix Multiplication Unit},
institution = {arXiv},
year = {2024},
number = {2306.11975v4}}

@Article{Sahraneshinsamani2024,
author={Sahraneshinsamani, Nima
and Catal{\'a}n, Sandra
and Herrero, Jos{\'e} R.},
title={{Mixed-precision pre-pivoting strategy for the LU factorization}},
journal={The Journal of Supercomputing},
year={2024},
month={Oct},
day={28},
volume={81},
number={1},
pages={87},
abstract={This paper investigates the efficient application of half-precision floating-point (FP16) arithmetic on GPUs for boosting LU decompositions in double (FP64) precision. Addressing the motivation to enhance computational efficiency, we introduce two novel algorithms: Pre-Pivoted LU (PRP) and Mixed-precision Panel Factorization (MPF). Deployed in both hybrid CPU-GPU setups and native GPU-only configurations, PRP identifies pivot lists through LU decomposition computed in reduced precision and subsequently reorders matrix rows in FP64 precision before executing LU decomposition without pivoting. Two variants of PRP, namely hPRP and xPRP, are introduced, differing in their computation of pivot lists in full half-precision or mixed half-single precision. The MPF algorithm generates FP64 LU factorization while internally utilizing hPRP for panel factorization, showcasing accuracy on par with standard DGETRF but with superior speed. The study further explores auxiliary functions required for the native mode implementation of PRP variants and MPF.},
% issn={1573-0484},
% doi={10.1007/s11227-024-06523-w},

url={https://doi.org/10.1007/s11227-024-06523-w}
}

@INPROCEEDINGS{mixed-prec-chol,
  author={Abdulah, Sameh and Ltaief, Hatem and Sun, Ying and Genton, Marc G. and Keyes, David E.},
  booktitle={2019 IEEE 26th International Conference on High Performance Computing, Data, and Analytics (HiPC)}, 
  title={{Geostatistical Modeling and Prediction Using Mixed Precision Tiled Cholesky Factorization}}, 
  year={2019},
  volume={},
  number={},
  pages={152-162},
  doi={10.1109/HiPC.2019.00028}}

@article{Amestoy2022Mixed,
  author    = {Amestoy, Patrick and Boiteau, Olivier and Buttari, Alfredo 
               and Gerest, Matthieu and J{\'e}z{\'e}quel, Fabienne and et al.},
  title     = {{Mixed Precision Low Rank Approximations and their Application 
               to Block Low Rank LU Factorization}},
  journal   = {IMA Journal of Numerical Analysis},
  year      = {2022},
  doi       = {10.1093/imanum/drac037},
  note      = {hal-03251738v3}
}

@misc{Ren2024Cholesky,
  author    = {Ren, Jie and Ltaief, Hatem and Abdulah, Sameh and Keyes, David E.},
  title     = {{Accelerating Mixed-Precision Out-of-Core Cholesky Factorization with Static Task Scheduling}},
% year      = {2024},
%  eprint    = {2410.09819},
%  archivePrefix = {arXiv},
  primaryClass = {cs.MS},
  note      = {Published on October 13, 2024},
  url       = {https://arxiv.org/pdf/2410.09819}
}

@INPROCEEDINGS{emulatedLU,
  author={Luszczek, Piotr and Gadepally, Vijay and Anderson, LaToya and Arcand, William and Bestor, David and Bergeron, William and Bonn, Alex and Burrill, Daniel J. and Byun, Chansup and Houle, Michael and Hubbell, Matthew and Jones, Michael and Michaleas, Peter and Morales, Guillermo and Mullen, Julia and Prout, Andrew and Reuther, Albert and Rosa, Antonio and Yee, Charles and Kepner, Jeremy},
  booktitle={2025 IEEE High Performance Extreme Computing Conference (HPEC)}, 
  title={Performance and Numerical Aspects of Decompositional Factorizations with FP64 Floating-Point Emulation in INT8}, 
  year={2025},
  volume={},
  number={},
  pages={1-7},
  doi={10.1109/HPEC67600.2025.11196192}}

@misc{ozaki2025ozakischemeiigemmoriented,
      title={{Ozaki Scheme II: A GEMM-oriented emulation of floating-point matrix multiplication using an integer modular technique}}, 
      author={Katsuhisa Ozaki and Yuki Uchino and Toshiyuki Imamura},
      year={2025},
      eprint={2504.08009},
      archivePrefix={arXiv},
      primaryClass={cs.MS},
      url={https://arxiv.org/abs/2504.08009}, 
}

@misc{luo2024ascendhifloat8formatdeep,
      title={Ascend HiFloat8 Format for Deep Learning}, 
      author={Yuanyong Luo and Zhongxing Zhang and Richard Wu and Hu Liu and Ying Jin and Kai Zheng and Minmin Wang and Zhanying He and Guipeng Hu and Luyao Chen and Tianchi Hu and Junsong Wang and Minqi Chen and Mikhaylov Dmitry and Korviakov Vladimir and Bobrin Maxim and Yuhao Hu and Guanfu Chen and Zeyi Huang},
      year={2024},
      eprint={2409.16626},
      archivePrefix={arXiv},
      primaryClass={cs.LG},
      url={https://arxiv.org/abs/2409.16626}, 
}

@article{10.14529/jsfi170206,
author = {Gustafson and Yonemoto},
title = {Beating Floating Point at its Own Game: Posit Arithmetic},
year = {2017},
issue_date = {June 2017},
publisher = {South Ural State University},
address = {Chelyabinsk, RUS},
volume = {4},
number = {2},
issn = {2409-6008},
url = {https://doi.org/10.14529/jsfi170206},
doi = {10.14529/jsfi170206},
journal = {Supercomput. Front. Innov.: Int. J.},
month = jun,
pages = {71–86},
numpages = {16}
}

@article{Uchino_2025,
   title={{Performance enhancement of the Ozaki Scheme on integer matrix multiplication unit}},
   volume={39},
   ISSN={1741-2846},
   url={http://dx.doi.org/10.1177/10943420241313064},
   DOI={10.1177/10943420241313064},
   number={3},
   journal={The International Journal of High Performance Computing Applications},
   publisher={SAGE Publications},
   author={Uchino, Yuki and Ozaki, Katsuhisa and Imamura, Toshiyuki},
   year={2025},
   month=Jan, pages={462–476} }

@article{Kolodziej2019,
  doi = {10.21105/joss.01244},
  url = {https://doi.org/10.21105/joss.01244},
  year = {2019},
  publisher = {The Open Journal},
  volume = {4},
  number = {35},
  pages = {1244},
  author = {Scott P. Kolodziej and Mohsen Aznaveh and Matthew Bullock and Jarrett David and Timothy A. Davis and Matthew Henderson and Yifan Hu and Read Sandstrom},
  title = {The SuiteSparse Matrix Collection Website Interface},
  journal = {Journal of Open Source Software}
}

@article{SuiteSparse2011,
author = {Davis, Timothy A. and Hu, Yifan},
title = {{The University of Florida Sparse Matrix Collection}},
year = {2011},
issue_date = {November 2011},
publisher = {Association for Computing Machinery},
address = {New York, NY, USA},
volume = {38},
number = {1},
issn = {0098-3500},
url = {https://doi.org/10.1145/2049662.2049663},
doi = {10.1145/2049662.2049663},
journal = {ACM Trans. Math. Softw.},
month = {dec},
articleno = {1},
numpages = {25}
}

@misc{uchino2026erroranalysismatrixmultiplication,
      title={{Error Analysis of Matrix Multiplication Emulation Using Ozaki-II Scheme}}, 
      author={Yuki Uchino and Katsuhisa Ozaki and Toshiyuki Imamura},
      year={2026},
      eprint={2602.02549},
      archivePrefix={arXiv},
      primaryClass={math.NA},
      url={https://arxiv.org/abs/2602.02549},
}

@manual{IEEE754-2019,
author = {{IEEE}},
title = {{IEEE} Standard for Floating-Point Arithmetic},
type = {Standard},
number = {IEEE Std 754-2019 (Revision of IEEE Std 754-2008)},
organization = {Institute of Electrical and Electronics Engineers},
address = {New York, NY, USA},
year = {2019},
month = {7},
doi = {10.1109/IEEESTD.2019.8766229},
isbn = {978-1-5044-5924-2},
url = {https://doi.org/10.1109/IEEESTD.2019.8766229},
}

@misc{demmel2026gradeaccuracyblas,
      title={How to grade the accuracy of the {BLAS}}, 
      author={James Demmel and Greg Henry and Igor Kozachenko and Julien Langou and Xiaoye Sherry Li and Jason Riedy and Jackson Vanover},
      year={2026},
      eprint={2609.12307},
      archivePrefix={arXiv},
      primaryClass={math.NA},
      url={https://arxiv.org/abs/2609.12307}, 
}

@book{GillMurrayWright1981,
  author    = {Gill, Philip E. and Murray, Walter and Wright, Margaret H.},
  title     = {Practical Optimization},
  publisher = {Academic Press},
  address   = {London},
  year      = {1981}
}

@misc{LoFloat,
  author       = {Sudhanva Kulkarni, James Demmel, Laura Grigori and David Chen},
  title        = {{LoFloat}: A {C++} Library for Low-Precision Floating-Point Experimentation},
  howpublished = {\url{https://github.com/SudhanvaKulkarni123/LoFloat}},
  year         = {2025},
  note         = {GitHub repository. Accessed: 2026-10-13}
}

@misc{Kouya2026OzakiLU,
      title={{Accelerating Multiple-Precision LU Decomposition with Ozaki Scheme II}},
      author={Tomonori Kouya},
      year={2026},
      eprint={2609.29554},
      archivePrefix={arXiv},
      primaryClass={math.NA},
      url={https://arxiv.org/abs/2609.29554},
}

@misc{Mukunoki2025FP8Ozaki,
      title={{DGEMM without FP64 Arithmetic -- Using FP64 Emulation and FP8 Tensor Cores with Ozaki Scheme}},
      author={Daichi Mukunoki},
      year={2025},
      eprint={2508.00441},
      archivePrefix={arXiv},
      primaryClass={cs.PF},
      url={https://arxiv.org/abs/2508.00441},
}

@misc{Hayashi2026FP4Ozaki,
      title={{DGEMM with Ozaki Scheme I/II on FP4 Tensor Cores: A Base-13 E2M1 Limb Representation}},
      author={Shun-ichiro Hayashi and Daichi Mukunoki and Tetsuya Hoshino and Takahiro Katagiri},
      year={2026},
      eprint={2608.06812},
      archivePrefix={arXiv},
      primaryClass={cs.DC},
      url={https://arxiv.org/abs/2608.06812},
}

@misc{Bake2025MixedPCG,
      title={{Forward and backward error bounds for a mixed precision preconditioned conjugate gradient algorithm}},
      author={Thomas Bake and Erin Carson and Yuxin Ma},
      year={2025},
      eprint={2510.11379},
      archivePrefix={arXiv},
      primaryClass={math.NA},
      url={https://arxiv.org/abs/2510.11379},
}

@misc{Kashi2024MixedPrecSurvey,
      title={{Mixed-precision numerics in scientific applications: survey and perspectives}},
      author={Aditya Kashi and Hao Lu and Wesley Brewer and David Rogers and Michael Matheson and Mallikarjun Shankar and Feiyi Wang},
      year={2024},
      eprint={2412.19322},
      archivePrefix={arXiv},
      primaryClass={cs.CE},
      url={https://arxiv.org/abs/2412.19322},
}

@misc{Cao2026RightPrecision,
      title={{Right Precision, Right Place: Runtime Adaptations for Mixed-Precision Cholesky Factorization}},
      author={Qinglei Cao and Theo Mary and Rabab Alomairy and Hatem Ltaief and David E. Keyes},
      year={2026},
      howpublished={HAL preprint hal-05793453},
      url={https://hal.science/hal-05793453v1},
}

@article{MXCholesky,
  TITLE = {{Error-Controlled Mixed-Precision Cholesky Factorization with MX Microscaling Formats}},
  AUTHOR = {Abdurakhmanov, Abdurauf and Mary, Theo and Ltaief, Hatem and Keyes, David E and Fahmy, Suhaib A.},
  URL = {https://hal.science/hal-05641535},
  JOURNAL = {{ACM Transactions on Mathematical Software}},
  PUBLISHER = {{Association for Computing Machinery}},
  YEAR = {2026},
  HAL_ID = {hal-05641535},
  HAL_VERSION = {v2},
}

\end{document}